\documentclass[12pt,twoside]{article}

\usepackage{graphicx}
\usepackage{color}
\usepackage[utf8]{inputenc}
\usepackage{natbib}

\usepackage{enumitem}
\usepackage{latexsym}
\usepackage{amssymb}
\usepackage{epsfig}
\usepackage{hyperref}
\usepackage{xcolor}

\usepackage[paper=a4paper, top=2cm, bottom=3cm, left=2.5cm, right=3.0cm]{geometry} 

\usepackage{amsmath}
\usepackage{float}
\usepackage{dsfont}
\usepackage{multirow}
\usepackage{bm}

\allowdisplaybreaks 

\def\bX{\mathbf{X}}
\def\dX{d_{\mathbf{X}}}
\def\bx{\mathbf{x}}
\def\bK{\mathbf{K}}

\newtheorem{theo}{Theorem}[section]
\newtheorem{lemma}[theo]{Lemma}
\newtheorem{rem}[theo]{Remark}

\makeatletter

\@addtoreset{equation}{section}
\makeatother

\begin{document}

	\title{\bf Estimation of the Transformation Function in Fully Nonparametric Transformation Models with Heteroscedasticity
	}
	
	\author{{\sc\Large Nick Kloodt}\footnote{Bundesstra\ss{}e 55, 20146 Hamburg, Nick.Kloodt@uni-hamburg.de, 040 428387167, ORCID 0000-0001-8998-1703}\\ Department of Mathematics, University of Hamburg}

	\maketitle

	\begin{abstract}
	Completely nonparametric transformation models with heteroscedastic errors are considered. Despite their flexibility, such models have rarely been used so far, since estimators of the model components have been missing and even identification of such models has not been clear until very recently. The results of \citet{Klo2020} are used to construct the first two estimators of the transformation function in these models. While the first estimator converges to the true transformation function at a parametric rate, the second estimator can be obtained by an explicit formula and is less computationally demanding. Finally, a simulation study is followed by some concluding remarks. Assumptions and proofs can be found in the appendix.
	\end{abstract}

	
	\section{Introduction}
	Consider two dependent random variables $Y$ and $\bX$, where $Y$ is univariate and the regressor $\bX$ is allowed to be multivariate. Analysing this dependency is often done by regression models. Nevertheless, with the aim to simplify the relationship between $Y$ and $X$ before fitting a regression model, more general transformation models have been applied more frequently during the last years. These models can be summarized as satisfying the model equation
	\begin{equation}\label{modeleq}
	h(Y)=g(\bX)+\sigma(\bX)\varepsilon,
	\end{equation}
	where $h$ is a strictly increasing transformation function, $g$ is the regression function and $\sigma^2$ is the variance function. The error $\varepsilon$ is usually assumed to be centred and independent of $\bX$ or at least some components of $\bX$.
	
	Since the beginnings of completely parametric modelling by \citet{BC1964} there has grown a large variety of transformation models. The parametric class of \citet{BC1964} was enlarged by \citet{YJ2000}, but there are meanwhile various classes of transformation functions, see \citet{ZR1969}, \citet{JD1980}, \citet{BD1981} or \citet{JP2009} for further examples. In their seminal article, \citet{LSvK2008} provided a parametric estimator of the transformation function in semiparametric models with nonparametric regression functions. Their results were extended by \citet{NNvK2016} to models with heteroscedastic errors. Nonparametric transformation functions have been considered by \citet{Hor1996} and \citet{EHN2004}. The arguably most general estimation results so far were provided by \citet{CKC2015} and \citet{VvK2018}, who considered general regression functions and homoscedastic errors, but allowed endogenous regressors. Although all of the approaches mentioned above fit to the framework of equation (\ref{modeleq}), none of these provides a method to estimate the transformation function nonparametrically in models with nonparametric regression and variance functions. In the following, this gap in the theory of transformation models is filled, by using the identification results of \citet{Klo2020} to define an estimator of the transformation function $h$ in the fully nonparametric transformation model (\ref{modeleq}) with heteroscedastic errors.
	
	The remainder is organized as follows. First, two nonparametric estimators are constructed from the identification result of \citet{Klo2020}. Afterwards, a parametric convergence rate is obtained for one of these estimators and the asymptotic distributions are given, before a short simulation study is followed by some concluding remarks. The proofs and assumptions are given in the Appendix.

	\section{The Estimator}
	Before the estimator can be defined, some notations are needed. Let $(Y,\bX),$ $(Y_1,\bX_1),...,$ $(Y_n,\bX_n)$ be independent and identically distributed $\mathbb{R}\times\mathbb{R}^{d_{\bX}}$-valued random variables fulfilling model (\ref{modeleq}). The densities of $\bX$ and $\varepsilon$ are denoted by $f=f_{\bX}$ and $f_{\varepsilon}$, respectively. As common in statistics, cumulated distribution functions will be denoted by a capital letter $F$. Consequently, let $F_{Y|\bX}$ denote the conditional cumulative distribution function of $Y$ given $\bX$. Further, let $f_{Y|\bX}:=\frac{f_{Y,\bX}}{f_{\bX}}$ be the corresponding conditional density. Moreover, define for arbitrary $y\in\mathbb{R}$, $\bx=(x_1,...,x_{d_{\bX}})\in\mathbb{R}^{\dX}$ and an appropriate index $j\in\{1,...,d_{\bX}\}$
	\begin{align*}
	p(y,\bx)&:=\int_{-\infty}^{y}f_{Y,\bX}(u,\bx)\,du,&p_y(y,\bx)&:=f_{Y,\bX}(y,\bx),
	\\[0,2cm]p_{x}(y,\bx)&:=\int_{-\infty}^{y}\frac{\partial}{\partial \bx_j}f_{Y,\bX}(u,\bx)\,du,&f_x(\bx)&:=\frac{\partial}{\partial \bx_j}f(\bx)
	\end{align*}
	as well as
	$$\Phi(y,\bx)=\frac{p(y,\bx)}{f_{\bX}(\bx)},\quad\Phi_y(y,\bx)=\frac{p_y(y,\bx)}{f_{\bX}(\bx)},\quad\Phi_{x}(y,\bx)=\frac{p_{x}(y,\bx)}{f_{\bX}(\bx)}-\frac{p(y,\bx)f_x(\bx)}{f_{\bX}^2(\bx)}.$$
	Here and in the following, the convention $f=f_{\bX}$ is used sometimes to denote the density of $\bX$ in order to make it better distinguishable from its derivative $f_x$. Define for some appropriate weight function $v$
	\begin{equation}\label{deflambda}
	\lambda(y):=\int v(\bx)\frac{\Phi_{x}(y,\bx)}{\Phi_y(y,\bx)}\,d\bx.
	\end{equation}
	Assumptions on $i,v$ and $F_{Y,\bX}$ are given in Appendix \ref{assumptions}. \citet{Klo2020} showed that $\lambda$ can be written as
	\begin{equation}\label{diffeq}
	\lambda(y)=-\frac{A+Bh(y)}{\frac{\partial}{\partial y}h(y)},
	\end{equation}
	where $A,B$ can be found in \ref{A4} and $B$ is an identified parameter, which is caused by the heteroskedasticity of the error in (\ref{modeleq}). Therefore, the transformation function can no longer be obtained by integrating $\frac{1}{\lambda}$ as was basically done by \citet{Hor2009} or \citet{CKC2015}. In contrast to this, \citet{Klo2020} found an explicit expression for $h$ based on the solution to the differential equation in (\ref{diffeq}). To formulate this expression, define $y_0:=\lambda^{-1}(0)$ and let $y_1,y_2\in\mathbb{R}$ be some fixed values with $y_1>y_0$ and $y_2<y_0$. The conditions
	\begin{equation}\label{idconstraints}
	h(y_0)=0,\quad\textup{and}\quad h(y_1)=1
	\end{equation}
	are used to identify the model. \citet{Klo2020} showed that the transformation function $h$ in (\ref{modeleq}) is identified under Assumptions \ref{A1}--\ref{A4} and can be expressed as
	\begin{equation}\label{exprh}
	h(y)=\left\{\begin{array}{ll}
	\exp\Big(-B\int_{y_1}^y\frac{1}{\lambda(u)}\,du\Big)&y>y_0\\
	0&y=y_0\\
	\alpha_2\exp\Big(-B\int_{y_2}^y\frac{1}{\lambda(u)}\,du\Big)&y<y_0\end{array}\right.,\vspace{0,3cm}
	\end{equation}
	where $\alpha_2$ is uniquely determined by requiring $\underset{y\searrow y_0}{\lim}\,h'(y)=\underset{y\nearrow _0}{\lim}\,h'(y)=h'(y_0)$ as
	\begin{equation}\label{defalpha2}
	\alpha_2=-\underset{t\rightarrow0}{\lim}\,\exp\bigg(B\bigg(\int_{y_2}^{y_0-t}\frac{1}{\lambda(u)}\,du-\int_{y_1}^{y_0+t}\frac{1}{\lambda(u)}\,du\bigg)\bigg).
	\end{equation}
	Due to the explicit expression of the transformation function in (\ref{exprh}), it seems natural to first find estimators of the unknown components $\lambda,B,y_0$ and $\alpha_2$ to use these estimators in a second step to define a plug in estimator of $h$.
	
	\subsection{Estimation of $\lambda$}
	Let $K:\mathbb{R}\rightarrow\mathbb{R}$ and $\bK:\mathbb{R}^{\dX}\rightarrow\mathbb{R}$ be kernel functions and let $h_x,h_y$ be bandwidths satisfying Assumption \ref{B3}. Let $K_{h_y}(y):=\frac{1}{h_y}K\big(\frac{y}{h_y}\big)$ as well as $\bK_{h_x}(y):=\frac{1}{h_x^{\dX}}\bK\big(\frac{\bx}{h_x}\big)$ and $\mathcal{K}_{h_y}(y):=\int_{-\infty}^yK_{h_y}(u)\,du$. Define an estimator
	$$\hat{f}_{\bX}(\bx):=\frac{1}{n}\sum_{i=1}^n\bK_{h_x}(\bx-\bX_i)$$
	of the marginal density $f_{\bX}$ as well as an estimator
	$$\hat{f}_{Y,\bX}(y|\bx):=\frac{1}{n}\sum_{i=1}^nK_{h_y}(y-Y_i)\bK_{h_x}(\bx-\bX_i)$$
	of the joint density $f_{Y,\bX}$ and its integrated version
	$$\hat{p}(y,\bx):=\frac{1}{n}\sum_{i=1}^n\mathcal{K}_{h_y}(y-Y_i)\bK_{h_x}(\bx-\bX_i).$$
	With these definitions, an estimator of the conditional distribution function of $Y$ conditional on $\bX$ is given by $\hat{\Phi}(y,\bx):=\hat{F}_{Y|\bX}(y|\bX):=\frac{\hat{p}(y,\bx)}{\hat{f}_{\bx}(\bx)}$, so that a plug in estimator of $\lambda$ can be defined as
	\begin{equation}\label{defestlambda}
	\hat{\lambda}(y):=\int v(\bx)\frac{\hat{\Phi}_{x}(y,\bx)}{\hat{\Phi}_y(y,\bx)}\,d\bx
	\end{equation}
	with $\hat{\Phi}_y(y,\bx):=\frac{\partial}{\partial y}\hat{\Phi}(y,\bx)$ and $\hat{\Phi}_x(y,\bx):=\frac{\partial}{\partial \bx_j}\hat{\Phi}(y,\bx)$.
	
	\subsection{Estimation of $y_0$ and $\alpha_2$}
	The estimation of $y_0$ is fairly straightforward. Using the estimator $\hat{\lambda}$ from equation (\ref{defestlambda}), we can simply define
	\begin{equation}\label{defesty0hat}
	\hat{y}_0:=\underset{y\,:\,\hat{\lambda}(y)=0}{\arg\min}\,|y|,
	\end{equation}
	which is equivalent to defining $\hat{y}_0=\hat{\lambda}^{-1}(0)$ as long as the inverse of zero is unique. Constructing an estimator of $\alpha_2$ is slightly more involved. Let $\tilde{B}$ be an estimator of $B$, e.g. the one given in equation (\ref{defestBtilde}) below. Further let $t_n\searrow0$ be an appropriate null sequence. Then, $\alpha_2$ can be estimated by
	\begin{equation}\label{defestalpha2hat}
	\hat{\alpha}_2=-\exp\bigg(\tilde{B}\bigg(\int_{y_2}^{\hat{y}_0-t_n}\frac{1}{\hat{\lambda}(u)}\,du-\int_{y_1}^{\hat{y}_0+t_n}\frac{1}{\hat{\lambda}(u)}\,du\bigg)\bigg).
	\end{equation}
	Here, $\hat{\lambda}$ and $\hat{y}_0$ are the estimators defined in (\ref{defestlambda}) and (\ref{defesty0hat}), respectively.

	\subsection{Estimation of $B$}\label{estB}
	There are two ways to estimate the unknown, but also identified parameter $B$. While the first approach presented here is rather simple and thus easy to calculate, the second estimator converges to $B$ at a faster rate. As was shown in \citet{Klo2020}, the function $\lambda$ can be expressed under the identification constraints (\ref{idconstraints}) as
	$$\lambda(y)=-\frac{Bh(y)}{\frac{\partial}{\partial y}h(y)},$$
	that is,
	$$\frac{\partial}{\partial y}\lambda(y)=-B\frac{\big(\frac{\partial}{\partial y}h(y)\big)^2-h(y)\frac{\partial^2}{\partial y^2}h(y)}{\big(\frac{\partial}{\partial y}h(y)\big)^2}\overset{y=y_0}{=}-B.$$
	Consequently, $B$ can be estimated via
	\begin{equation}\label{defestBtilde}
	\tilde{B}:=-\frac{\partial}{\partial y}\hat{\lambda}(y)\Big|_{y=\hat{y}_0},
	\end{equation}
	where $\hat{y}_0$ denotes the estimator from equation (\ref{defesty0hat}).
	
	The second approach is more sophisticated and uses the idea of the Mean-Square-Distance-from-Independence estimator of \citet{LSvK2008}. For two random variables $U,V$, where $U$ is real valued, and some $\tau,\beta\in(0,1),\tau\neq\beta,$ denote the $\tau$-quantile of $U$ conditional on $V=v$ by
	$$F_{U|V}^{-1}(\tau|v):=\inf\,\{u\in\mathbb{R}:F_{U|V}(u|v)\geq\tau\}$$
	and define
	$$\tilde{\varepsilon}:=\frac{h(Y)-h(F_{Y|\bX}^{-1}(\tau|\bX))}{F_{\sigma(\bX)(\varepsilon-F_{\varepsilon}^{-1}(\tau))|\bX}^{-1}(\beta|\bX)}=\frac{\varepsilon-F_{\varepsilon}^{-1}(\tau)}{F_{\varepsilon}^{-1}(\beta)-F_{\varepsilon}^{-1}(\tau)}.$$
	To be able to take account of the variability of the estimated errors which are caused by the uncertainty of the estimated model components, let $\mathfrak{h},f_{m_{\tau}},f_{m_{\beta}}$ belong to some function sets specified in (\ref{defHcal}) and (\ref{defHcaltilde}) in Appendix \ref{auxiliary} and define
	\begin{equation}\label{defhc}
	h_c(y)=\exp\bigg(-c\int_{y_1}^{y}\frac{1}{\lambda(u)}\,du\bigg)\quad\textup{and}\quad \hat{h}_c(y)=\exp\bigg(-c\int_{y_1}^{y}\frac{1}{\hat{\lambda}(u)}\,du\bigg)
	\end{equation}
	for $y>y_0$ as well as
	$$s:=(\mathfrak{h},f_{m_{\tau}},f_{m_{\beta}}),\quad s_0:=(h_1,F_{Y|\bX}^{-1}(\tau|\cdot),F_{Y|\bX}^{-1}(\beta|\cdot))$$
	and
	$$\tilde{\varepsilon}_c(s)=\frac{\mathfrak{h}_c(Y)-\mathfrak{h}_c(f_{m_{\tau}}(\bX))}{\mathfrak{h}_c(f_{m_{\beta}}(\bX))-\mathfrak{h}_c(f_{m_{\tau}}(\bX))}.$$
	Therefore, one has $\tilde{\varepsilon}=\tilde{\varepsilon}(s_0)$ with the convention $\mathfrak{h}_c(y)=\operatorname{sign}(\mathfrak{h}(y))|\mathfrak{h}(y)|^{c}$ for all $c>0$. For $\hat{s}:=(\hat{h}_1,\hat{F}_{Y|\bX}^{-1}(\tau|\cdot),\hat{F}_{Y|\bX}^{-1}(\beta|\cdot))$ the errors $\varepsilon_1,...,\varepsilon_n$ can be estimated by
	$$\hat{\varepsilon}_{c,i}:=\tilde{\varepsilon}_{c,i}(\hat{s})=\frac{\hat{h}_c(Y_i)-\hat{h}_c(\hat{F}_{Y|\bX}^{-1}(\tau|\bX_i))}{\hat{h}_c(\hat{F}_{Y|\bX}^{-1}(\beta|\bX_i))-\hat{h}_c(\hat{F}_{Y|\bX}^{-1}(\tau|\bX_i))},\quad i=1,...,n.$$
	Similar to \citet{LSvK2008}, the idea behind the second estimator is based on the fact that $\tilde{\varepsilon}_c(s_0)$ is independent of $\bX$ if and only if $c=B$ \citep{Klo2019}. This can be transferred to an $\mathcal{L}^2$-criterion by defining
	\begin{equation}\label{defA}
	A(c,s):=\sqrt{\int_{M_{\bX}}\int_{[e_a,e_b]}G_{MD}(c,s)(\bx,e)^2\,de\,d\bx}=:||G_{MD}(c,s)||_2,
	\end{equation}
	where $M_{\bX}$ and $[e_a,e_b]$ are some appropriate compact intervals of $\mathbb{R}^{\dX}$ and $\mathbb{R}$, respectively, and $G_{MD}(c,s)$ is defined as
	\begin{align}
	G_{MD}(c,s)(\bx,e)&=P\big(\bX\leq \bx,\tilde{\varepsilon}_c(s)\leq e|\bX\in M_{\bX}\big)\nonumber
	\\[0,2cm]&\quad-P\big(\bX\leq \bx|\bX\in M_{\bX}\big)P\big(\tilde{\varepsilon}_c(s)\leq e|\bX\in M_{\bX}\big).\label{defGMD}
	\end{align}
	Assumptions on the sets $M_{\bX}$ and $[e_a,e_b]$ are given in Appendix \ref{assumptions}. Then, $c=B$ is equivalent to $||G_{MD}(c,s_0)||_2=0$. In order to use this property to construct an estimator of $B$, we write
	\begin{align}
	\hat{P}\big(\bX\leq \bx,\tilde{\varepsilon}_c(s)\leq e|\bX\in M_{\bX}\big)&=\frac{\frac{1}{n}\sum_{i=1}^nI_{\{\tilde{\varepsilon}_{c,i}(s)\leq e\}}I_{\{\bX_i\leq \bx\}}I_{\{\bX_i\in M_{\bX}\}}}{\frac{1}{n}\sum_{i=1}^nI_{\{\bX_i\in M_{\bX}\}}},\nonumber
	\\[0,2cm]\hat{P}\big(\bX\leq \bx|\bX\in M_{\bX}\big)&=\frac{\frac{1}{n}\sum_{i=1}^nI_{\{\bX_i\leq \bx\}}I_{\{\bX_i\in M_{\bX}\}}}{\frac{1}{n}\sum_{i=1}^nI_{\{\bX_i\in M_{\bX}\}}},\nonumber
	\\[0,2cm]\hat{P}\big(\tilde{\varepsilon}_c(s)\leq e|\bX\in M_{\bX}\big)&=\frac{\frac{1}{n}\sum_{i=1}^nI_{\{\tilde{\varepsilon}_{c,i}(s)\leq e\}}I_{\{\bX_i\in M_{\bX}\}}}{\frac{1}{n}\sum_{i=1}^nI_{\{\bX_i\in M_{\bX}\}}},\nonumber
	\\[0,2cm]G_{nMD}(c,s)(\bx,e)&=\hat{P}\big(\bX\leq \bx,\tilde{\varepsilon}_c(s)\leq e|\bX\in M_{\bX}\big)\nonumber
	\\[0,2cm]&\quad-\hat{P}\big(\bX\leq \bx|\bX\in M_{\bX}\big)\hat{P}\big(\tilde{\varepsilon}_c(s)\leq e|\bX\in M_{\bX}\big).\label{defGnMD}
	\end{align}
	Finally, define $\hat{A}(c,s):=||G_{nMD}(c,s)||_2$ and
	\begin{equation}\label{defestBhat}
	\hat{B}:=\underset{c\in[B_1,B_2]}{\arg\min}\,\hat{A}(c,\hat{s})
	\end{equation}
	for some interval $[B_1,B_2]\subseteq(0,\infty)$ which contains $B$.
	
	\subsection{Putting Things together}
	Once we have estimated all of the unknown components in expression (\ref{exprh}), we can simply insert them to obtain a plug in estimator. Note that we forced the transformation function $h$ by condition (\ref{idconstraints}) to be zero at $y_0$. Consequently, the integral $\int_{y_1}^{y}\frac{1}{\lambda(u)}\,du$ diverges for $y\searrow y_0$, which might cause problems when calculating $\int_{y_1}^{y}\frac{1}{\hat{\lambda}(u)}\,du$ in a \mbox{neighbourhood} of $y_0$ and therefore complicates the estimation of $h$ for such values. As a solution to this problem, a null sequence $t_n\searrow0$ can be used similarly to the estimation of $\alpha_2$ to define linearised estimators
	\begin{equation}\label{defesthglobalhat}
	\hat{h}(y)=\left\{\begin{array}{ll}
	\exp\Big(-\hat{B}\int_{y_1}^y\frac{1}{\hat{\lambda}(u)}\,du\Big),&y\geq\hat{y}_0+t_n,\\[0,2cm]
	\frac{y-\hat{y}_0}{t_n}\hat{h}(\hat{y}_0+t_n),&y\in(\hat{y}_0,\hat{y}_0+t_n),\\[0,2cm]
	0,&y=\hat{y}_0,\\[0,2cm]
	\frac{\hat{y}_0-y}{t_n}\hat{h}(\hat{y}_0-t_n),&y\in(\hat{y}_0-t_n,\hat{y}_0),\\[0,2cm]
	\hat{\lambda}_2\exp\Big(-\hat{B}\int_{y_2}^y\frac{1}{\hat{\lambda}(u)}\,du\Big),&y\leq\hat{y}_0-t_n\end{array}\right.
	\end{equation}
	and
	\begin{equation}\label{defesthglobaltilde}
	\tilde{h}(y)=\left\{\begin{array}{ll}
	\exp\Big(-\tilde{B}\int_{y_1}^y\frac{1}{\hat{\lambda}(u)}\,du\Big),&y\geq\hat{y}_0+t_n,\\[0,2cm]
	\frac{y-\hat{y}_0}{t_n}\tilde{h}(\hat{y}_0+t_n),&y\in(\hat{y}_0,\hat{y}_0+t_n),\\[0,2cm]
	0,&y=\hat{y}_0,\\[0,2cm]
	\frac{\hat{y}_0-y}{t_n}\tilde{h}(\hat{y}_0-t_n),&y\in(\hat{y}_0-t_n,\hat{y}_0),\\[0,2cm]
	\tilde{\lambda}_2\exp\Big(-\tilde{B}\int_{y_2}^y\frac{1}{\hat{\lambda}(u)}\,du\Big),&y\leq\hat{y}_0-t_n.\end{array}\right.
	\end{equation}
	Here, the notations of $\hat{h}$ and $\tilde{h}$ are related to those of $\hat{B}$ and $\tilde{B}$, respectively.

	\section{Asymptotic Results}
	In this section, asymptotic results for the estimators given in (\ref{defesthglobalhat}) and (\ref{defesthglobaltilde}) are presented. Due to the plug in type of $\hat{h}$ and $\tilde{h}$, the asymptotic behaviour of these estimators can be obtained from those of the single components. Apart from those for $\hat{h}$ and $\tilde{h}$, the main focus will especially lie on convergence results for the estimators $\hat{\lambda},\hat{B}$ and $\tilde{B}$.
	
	In order to provide an asymptotic expression for $\hat{\lambda}-\lambda$, let
	\begin{align*}
	D_{p,0}&:=-\frac{f_x}{\Phi_yf^2},&D_{p,y}&:=-\frac{\Phi_x}{\Phi_y^2f},&D_{p,x}&:=\frac{1}{\Phi_yf},
	\\[0,2cm]D_{f,0}&:=\frac{2pf_x}{\Phi_yf^3}-\frac{p_x}{\Phi_yf^2}+\frac{p_y\Phi_x}{\Phi_y^2f^2},&D_{f,x}&:=-\frac{p}{\Phi_yf^2}&&
	\end{align*}
	\begin{lemma}\label{asympintegral}
		Assume \ref{A1}--\ref{A4} and \ref{B1}--\ref{B5} and let $\mathcal{K}\subseteq\mathbb{R}$ be compact. Then, with $\hat{\lambda}$ as in (\ref{defestlambda}) one has
		\begin{align}
		\hat{\lambda}(u)-\lambda(u)&=\frac{1}{n}\sum_{i=1}^n\bigg(v(\bX_i){D}_{p,0}(u,\bX_i)\mathcal{K}_{h_y}(u-Y_i)-\frac{\partial\big(v(\bX_i){D}_{p,x}(u,\bX_i)\big)}{\partial \bx_j}\mathcal{K}_{h_y}(u-Y_i)\nonumber
		\\[0,2cm]&\quad+v(\bX_i){D}_{p,y}(u,\bX_i)K_{h_y}(u-Y_i)+v(\bX_i){D}_{f,0}(u,\bX_i)\nonumber
		\\[0,2cm]&\quad-\frac{\partial\big(v(\bX_i){D}_{f,x}(u,\bX_i)\big)}{\partial \bx_j}\bigg)+o_p\bigg(\frac{1}{\sqrt{n}}\bigg)\nonumber
		\\[0,2cm]&=\mathcal{O}_p\Bigg(\sqrt{\frac{\log(n)}{nh_y}}\Bigg)\label{uniconvlambda}
		\end{align}
		uniformly in $y\in\mathcal{K}$. Furthermore, the process $(Z_n(y))_{y\in\mathcal{K}}$ defined by
		$$Z_n(y):=\sqrt{n}\int_{y_1}^y\bigg(\frac{1}{\hat{\lambda}(u)}-\frac{1}{\lambda(u)}\bigg)\,du$$
		converges weakly to a centred Gaussian process $Z_{\lambda}$ with a covariance function which can be found in \citet{Klo2019}.
	\end{lemma}\noindent
	The proof can be found in Section \ref{proofasympintegral}. Similar techniques as in the article of \citet{CKC2015} are applied.
	
	The following result can be shown for the estimators $\hat{y}_0$ and $\hat{\alpha}_2$.
	\begin{lemma}\label{asympy0alpha2}
		\ref{A1}--\ref{A4} and \ref{B1}--\ref{B5} in Appendix \ref{assumptions}. Moreover, let $\hat{y}_0$ and $\hat{\alpha}_2$ be defined as in (\ref{defesty0hat}) and (\ref{defestalpha2hat}), respectively, where the sequence $t_n>0$ fulfils $t_n\sim\big(\frac{\log(n)^2}{nh_y}\big)^{\frac{1}{4}}$. Then,
		$$\sqrt{nh_y}(\hat{y}_0-y_0)\overset{\mathcal{D}}{\rightarrow}\mathcal{N}(0,\sigma_{y_0}^2)$$
		with
		$$\sigma_{y_0}^2=\frac{\int K(z)^2\,dz}{B^2}\int v(\mathbf{w})^2{D}_{p,y}(y_0,\mathbf{w})^2f_{Y,\bX}(y_0,\mathbf{w})\,d\mathbf{w}$$
		and
		$$\hat{\alpha}_2-\alpha_2=\frac{\alpha_2\frac{\partial^2}{\partial y^2}h(y_0)t_n}{\frac{\partial}{\partial y}h(y_0)}+o_p(t_n).$$
	\end{lemma}\noindent
	Since especially the second part of the proof of Lemma \ref{asympy0alpha2} is rather technical and does not really contribute to the overall picture of the article, the proof of Lemma \ref{asympy0alpha2} is omitted here. It can be found in \citet{Klo2019}.
	
	\begin{lemma}\label{asympB}
		Let $d_X=1$ and assume \ref{A1}--\ref{A4}, \ref{B1}--\ref{B5}and \ref{M1}--\ref{M5} from Appendix \ref{assumptions} for some compact interval $[z_a,z_b]\subseteq(y_0,\infty)$. Moreover, consider the estimator $\hat{s}=(\hat{h}_1,\hat{F}_{Y|\bX}^{-1}(\tau|\cdot),\hat{F}_{Y|\bX}^{-1}(\beta|\cdot))$ of $s_0$.
		Then,
		$$\sqrt{n}(\hat{B}-B)\overset{\mathcal{D}}{\rightarrow}\mathcal{N}\bigg(0,\frac{\sigma_A^2}{||\Gamma_{1}(B,s_0)||_2^4}\bigg),$$
		where $\Gamma_1(B,s_0)(\bx,e):=\frac{\partial}{\partial c}G_{MD}(c,s_0)(\bx,e)\big|_{c=B}$ and $\sigma_A>0$ is defined in \ref{C6} in Appendix \ref{auxiliary}.
	\end{lemma}\noindent
	The proof can be found in Subsection \ref{proofasympB}. The assumptions \ref{B2} and \ref{B4} are chosen in a way such that the derivatives of $\hat{\lambda}$ and $\hat{F}_{Y|\bX}(\tau|\cdot)$ up to order two converge uniformly to those of $\lambda$ and $F_{Y|\bX}(\tau|\cdot)$, respectively, on compact sets. Hence, Lemma \ref{lemmac4} in Appendix \ref{auxiliary} and consequently Lemma \ref{asympB} only hold for $d_X=1$. It is conjectured that the assumptions can be adjusted such that the convergence results can be extended to higher order derivatives which would generalize these Lemmas to arbitrary $\dX>1$. Nevertheless, the case $\dX=1$ is considered for simplicity in the following.
	
	Another approach presented in \ref{estB} consisted in estimating $B$ via $\tilde{B}$ from (\ref{defestBtilde}).
	\begin{lemma}\label{asympaltB}
		Assume \ref{A1}--\ref{A4} and \ref{B1}--\ref{B5} from Appendix \ref{assumptions}. Then,
		\begin{equation}\label{convBtilde}
		\sqrt{nh_y^3}(\tilde{B}-B)\overset{\mathcal{D}}{\rightarrow}\mathcal{N}(0,\sigma_{\tilde{B}}^2),
		\end{equation}
		where
		$$\sigma_{\tilde{B}}^2=\bigg(\int\bigg(\frac{\partial}{\partial z}K(z)\bigg)^2\,dz\bigg)\bigg(\int v(\mathbf{w})^2D_{p,y}(y_0,\mathbf{w})^2f_{Y,\bX}(y_0,\mathbf{w})\,d\mathbf{w}\bigg).$$
	\end{lemma}
	The proof is given in Subsection \ref{proofasympaltB}.
	
	It remains to combine the previous lemmas to the main result of this article, a convergence theorem for $\hat{h}$ and $\tilde{h}$, respectively. When doing so, the rate of uniform convergence depends on the set over which the convergence is considered.
	\begin{theo}\label{asymph}
		Let $\hat{s}=(\hat{h}_1,\hat{F}_{Y|\bX}^{-1}(\tau|\cdot),\hat{F}_{Y|\bX}^{-1}(\beta|\cdot))$ be the estimator from Section \ref{estB} and let $\mathcal{K}\subseteq(y_0,\infty),\tilde{\mathcal{K}}\subseteq\mathbb{R}$ be compact sets.
		\begin{enumerate}[label=(\roman*)]
			\item\label{asymphi} Under the assumptions of Lemma \ref{asympB}, the stochastic process $(H_n(y))_{y\in\mathcal{K}}$ defined by
			$$H_n(y):=\sqrt{n}(\hat{h}(y)-h(y))$$
			converges weakly to a centred Gaussian process $(Z_h(y))_{y\in\mathcal{K}}$ with the covariance function
			\begin{align*}
			\kappa_h(u,v)&=h(u)h(v)E\bigg[\bigg(B\eta_1(u)+\int_{y_1}^{u}\frac{1}{\lambda(y)}\,dy\ \psi_{\Gamma_2}(Y_1,\bX_1)\bigg)
			\\[0,2cm]&\quad\quad\bigg(B\eta_1(v)+\int_{y_1}^{v}\frac{1}{\lambda(y)}\,dy\ \psi_{\Gamma_2}(Y_1,\bX_1)\bigg)\bigg]
			\end{align*}
			with $\eta_1$ and $\psi_{\Gamma_2}$ as in (\ref{defeta}) and \ref{C6} in Appendix \ref{auxiliary}.
			\item\label{asymphii} Under the assumptions of Lemma \ref{asympaltB}, the process $(\tilde{H}_n(y))_{y\in\mathcal{K}}$ defined by $\tilde{H}_n(y)=\sqrt{nh_y^3}(\tilde{h}(y)-h(y))$ with $\tilde{h}$ as in (\ref{defesthglobaltilde}) converges weakly to the centred Gaussian process
			$$(Z_{\tilde{h}}(y))_{y\in\mathcal{K}}=\bigg(h(y)\int_{y_1}^{y}\frac{1}{\lambda(u)}\,du\ Z_{\tilde{B}}\bigg)_{y\in\mathcal{K}}$$
			with $Z_{\tilde{B}}$ from Lemma \ref{asympaltB}. Furthermore, if $t_n\sim\big(\frac{\log(n)^2}{nh_y}\big)^{\frac{1}{4}}$, it holds that
			$$\underset{y\in\tilde{\mathcal{K}}}{\sup}\,|\tilde{h}(y)-h(y)|=\mathcal{O}_p\Bigg(\bigg(\frac{\log(n)^2}{nh_y}\bigg)^{\frac{1}{4}}\Bigg).$$
		\end{enumerate}
	\end{theo}\noindent
	The proof can be found in Section \ref{proofasymph}. Theorem \ref{asymph} yields uniform consistency of both of the estimators $\hat{h}$ and $\tilde{h}$ on compact sets. These estimators of the transformation function consist of single estimators of the unknown components in (\ref{exprh}). Therefore, the rate of convergence to the true transformation function $h$ corresponds to the rate of these estimators. The unknown components $B$ and $\int_{y_1}^y\frac{1}{\lambda(u)}\,du$ can be estimated at a $n^{-\frac{1}{2}}$-rate uniformly on compacts sets $y\in\mathcal{K}\subseteq(y_0,\infty)$, so can $h$. When estimating $h$ for values $y$, which are below or close to $y_0$, the rate worsens, since the estimators of $y_0$ and $\alpha_2$ do not converge at a $n^{-\frac{1}{2}}$-rate. Although $\hat{B}$ converges to $B$ at a faster rate than $\tilde{B}$ does, the calculation of $\tilde{B}$ is due to the explicit formula in (\ref{defestBtilde}) much easier and less computationally demanding than that of $\hat{B}$.
	\begin{rem}
		\citet{CvK2018} adjusted the estimator of \cite{CKC2015} by first transforming $Y$ with the empirical distribution function $\hat{F}_Y$. In the context of the heteroscedastic model presented here, such a pretransformation is conceivable as well and it is conjectured that the influence on the asymptotic distribution will be similar to that in \citet{CvK2018}.
	\end{rem}

	\section{Simulation Study}
	The purpose of this section consists in providing an estimation approach, which works and is practically applicable. Moreover, the behaviour of the estimators of $y_0,B$ and $h$ given in (\ref{defesty0hat}), (\ref{defestBtilde}) and (\ref{defesthglobaltilde}), respectively, for finite sample sizes is examined. The estimator $\tilde{B}$ is chosen for practical reasons, since calculating $\hat{B}$ is accompanied with a rather complex and computationally demanding optimization problem, which requires the choice of additional bandwidths, kernels and sets, whereas equation (\ref{defestBtilde}) can be applied fast and easily.
	
	Independent observations of real valued random variables $\bX\sim\mathcal{U}([0,1])$ and $\varepsilon\sim\mathcal{U}([-1,1])$ are generated. Afterwards, $Y$ is defined by
	$$Y=\frac{\big(1+\bX+\frac{(1+\bX)^2}{2}\varepsilon\big)^3}{8}+\frac{7\big(1+\bX+\frac{(1+\bX)^2}{2}\varepsilon\big)}{8},$$
	that is, model (\ref{modeleq}) is fulfilled with
	$$h^{-1}(y)=\frac{y^3}{8}+\frac{7y}{8},\quad g(\bx)=1+\bx\quad\textup{and}\quad\sigma(\bx)=\frac{(1+\bx)^2}{2}.$$
	The transformation function $h$ is chosen such that it is strictly monotonic. Furthermore, it fulfils the identification conditions $h(0)=0$ and $h(1)=1$ and thus needs to be linearly transformed later when comparing it to the estimator $\tilde{h}$.
	
	The simulations are conducted with the language \textit{R} (\citet{R2013}). Some of the already implemented commands such as \textit{integrate} and \textit{h.select} are applied and an interface for C++ is used to reduce the computation time. The weight function $v$ is chosen to be the indicator function of $[0,1]$. Similarly to \citet{CvK2018}, the mean of $N_x=100$ evaluations of $$\bx\mapsto v(\bx)\frac{\frac{\partial \hat{F}_{Y|\bX}(y|\bx)}{\partial \bx_j}}{\frac{\partial \hat{F}_{Y|\bX}(y|\bx)}{\partial y}}$$
	at equidistant points between the minimum and the maximum of the observations of $\bX$ is taken instead of integrating the quotient as in (\ref{defestlambda}). To calculate the bandwidths $h_y$ and $h_x$, cross validation and the normal reference rule, respectively, have been applied (\citet{Sil1986}). The kernel $K$ is chosen to be the Epanechnikov kernel.
	
\begin{figure}[htbp]
		\centering
		\includegraphics[width=15cm,height=15cm]{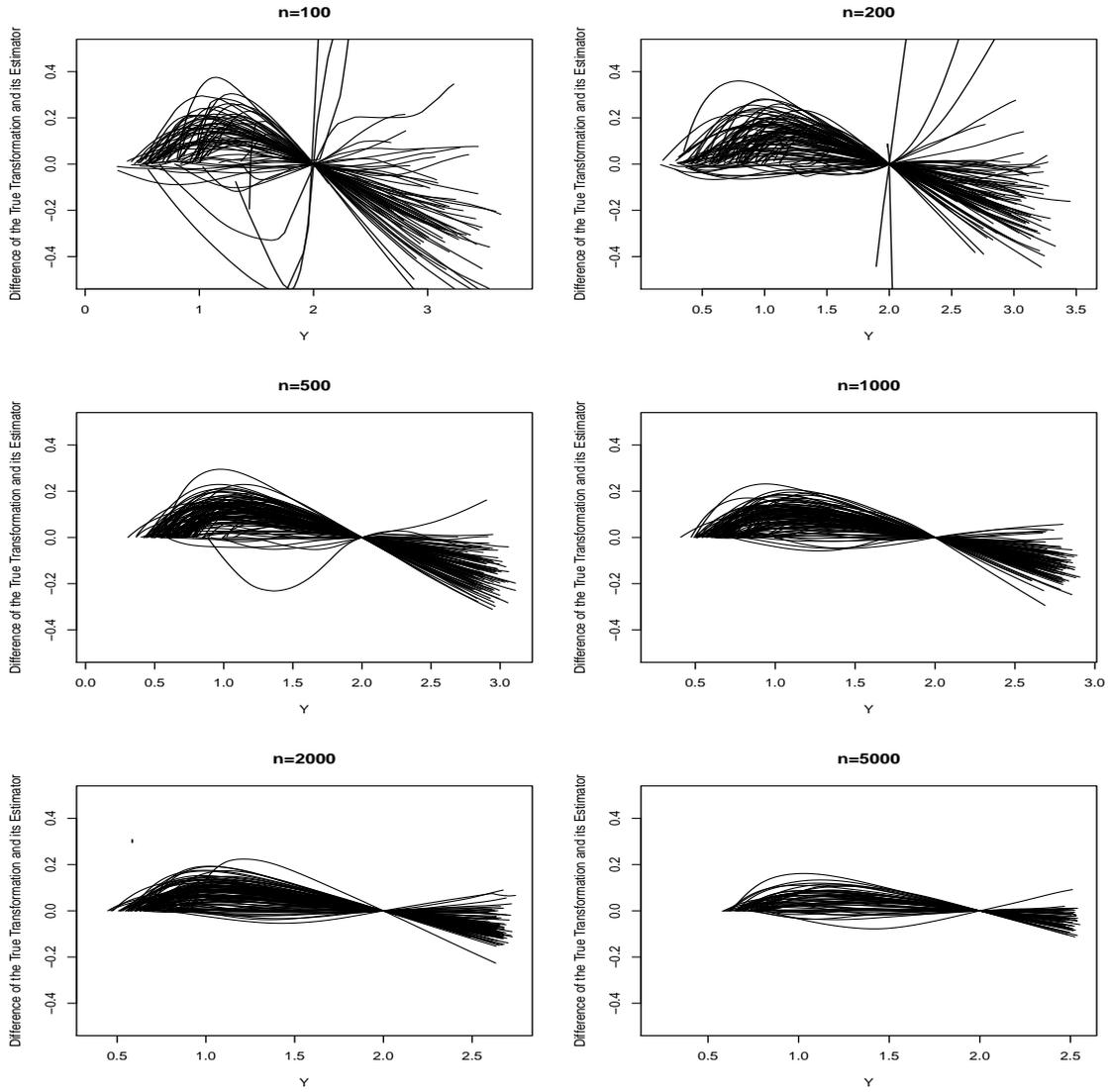}
		\caption[Difference of Nonparametric Estimators and the True Transformation]{The difference of the true transformation function and its estimator under the same identification conditions is shown for the sample sizes of $n=100,n=200,n=500,n=1000,n=2000,n=5000,n=10000$.}
		\label{est_differences_trafo}
	\end{figure}

	It can be shown that $B=\log(4)\approx1.39$ and $y_0=\frac{1}{8\log(4)^3}+\frac{7}{8\log(4)}\approx0.68.$
	Observations are simulated for sample sizes of $n\in\{100,200,500,1000,2000,5000,10000\}$. For computational reasons, the number $m$ of simulation runs for each of the scenarios decreases with the sample size and can be found in Table \ref{_repetition_samplesize}.
	\begin{table}[htbp]
		\centering
		\begin{tabular}{|l|c|c|c|c|c|c|c|}
			\hline
			Sample Size $n$&100&200&500&1000&2000&5000&10000\\
			\hline
			Number of Simulation Runs $m$&500&500&200&200&100&50&20\\
			\hline
		\end{tabular}
		\caption[Number of Simulation Runs Depending on the Sample Size]{The sample sizes and the corresponding number of the simulation runs.}
		\label{_repetition_samplesize}
	\end{table}\noindent
	\begin{figure}[htbp]
		\centering
		\includegraphics[width=10cm,height=10cm]{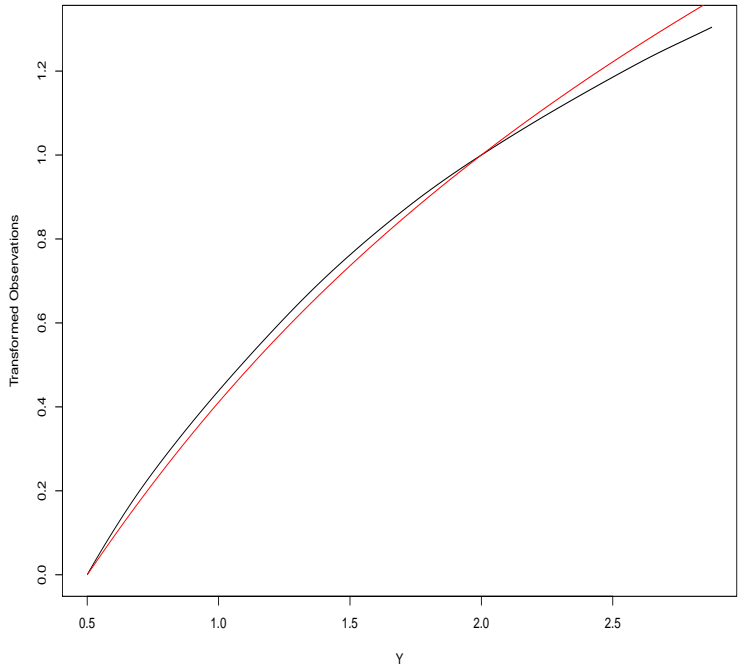}
		\caption[A Nonparametric Estimate of the Heteroscedastic Transformation]{One realization of the estimated transformation function (black curve) and the true transformation function (red curve) are shown for $n=500$.}
		\label{hetmodel_estimated_transformation}
	\end{figure}\noindent
	Figure \ref{hetmodel_estimated_transformation} shows a realization of the estimator $\tilde{h}$ in (\ref{defesthglobaltilde}) which is based on $n=500$ observations in black and the true transformation function $h$ in red, both for $y>\hat{y}_0$. Here and in the following, the identification constraints in (\ref{idconstraints}) for $\hat{y}_0$ and $y_1=2$ are used, that is $\tilde{h}(\hat{y}_0)=0=h(\hat{y}_0)$ and $\tilde{h}(2)=1=h(2)$. Therefore, both functions have to intersect at least in $\hat{y}_0$ and $y_1=2$. The approximation in Figure \ref{hetmodel_estimated_transformation} seems to be quite good, although the estimator for values below $y_1=2$ slightly overestimates the true transformation function, whereas the opposite holds for values above $y_1=2$. As can be seen in Graphic \ref{est_differences_trafo}, this phenomenon carries over to all of the simulated scenarios. There, the difference $\tilde{h}-h$ of the estimator and the true transformation function, again based on the same identification conditions, for different sample sizes is shown. Table \ref{By0MISE} indicates that this bias is caused by an underestimation of $B$, since it holds that $h(y)\in(0,1)$ if and only if $y\in(y_0,2)$.
	\noindent
	\begin{table}[htbp]
		\centering
		\begin{tabular}{|l|c|c|c|}
			\hline
			Sample Size&Mean of $\hat{y}_0$&Mean of $\tilde{B}$&Est. MISE of $\tilde{h}$\\
			\hline
			$n=100$&1.14&0.80&33.19\\
			$n=200$&0.88&0.76&12.90\\
			$n=500$&0.64&0.81&2.38\\
			$n=1000$&0.65&0.85&2.25\\
			$n=2000$&0.66&0.99&2.19\\
			$n=5000$&0.71&1.10&2.30\\
			$n=10000$&0.66&1.16&1.92\\
			True Values&0.68&1.39&\\
			\hline
		\end{tabular}
		\caption[Means of the Estimators of some Model Components]{Means of the estimators $\hat{y}_0$ and $\tilde{B}$ as well as the estimated MISE of the estimated transformation function for the sample sizes of $n=100,n=200,n=500,n=1000,n=2000,n=5000,n=10000$.}
		\label{By0MISE}
	\end{table}\noindent
	Whereas $\hat{y}_0$ already seems to be unbiased for $n=500$, the value of $\tilde{B}$ is even for $n=10000$ below the true value of $B=1.39$, although the gap between $\tilde{B}$ and $B$ decreases with a growing sample size.
	
\begin{figure}[htbp]
		\centering
		\includegraphics[width=10cm,height=10cm]{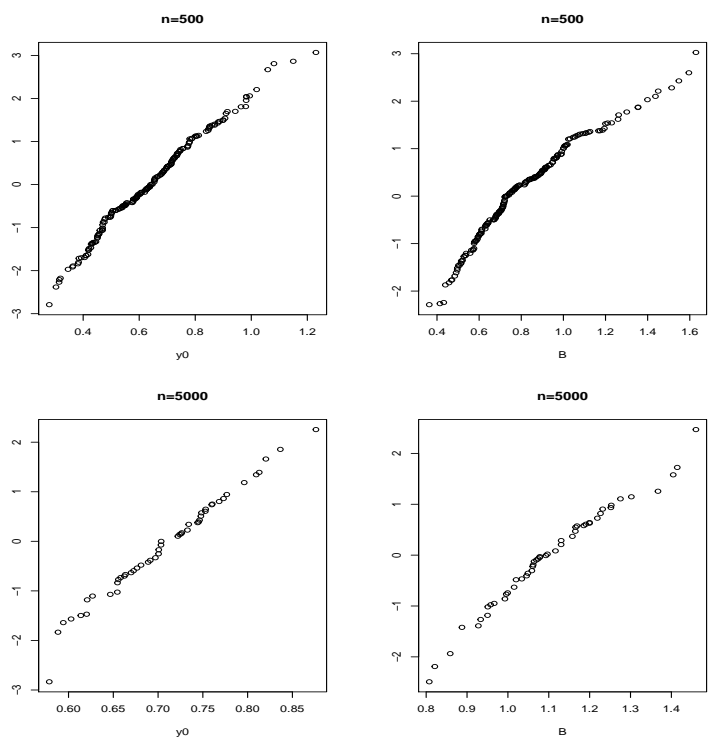}
		\caption[QQ-Plots of the Estimators of some Model Components]{Normal-QQ-Plots of the estimators $\hat{y}_0$ and $\tilde{B}$ for the sample sizes of $n=500$ and $n=5000$.}
		\label{est_qqplots_y0B}
	\end{figure}

	Finally, some QQ-plots for $\hat{y}_0$ and $\tilde{B}$ are given in Figure \ref{est_qqplots_y0B}. There, the empirical quantiles of the estimators are compared to those of standard normally distributed random variables. While the distribution of $\hat{y}_0$ seems to be almost normal already for a sample size of $n=500$, the corresponding curve for $\tilde{B}$ has a small bump for $n=500$, but at least seems to be linear for $n=5000$.

	\section{Conclusion}
	Based on the results of \citet{Klo2020}, the so far most general approach for estimating the transformation function in the heteroscedastic model (\ref{modeleq}) has been developed. Depending on the chosen approach for estimating $B$, two estimators $\hat{h}$ and $\tilde{h}$ related to that of \citet{CKC2015} have been provided. Consistency as well as weak convergence results for the proposed estimators and its components have been proven. Moreover, the estimators and convergence results have been promoted by a simulation study.
	
	Since this has been the first step in the context of estimating the transformation function in fully nonparametric and heteroskedastic models like (\ref{modeleq}), there is potential for further adjustments and examinations. Future research could consist in simplifying the estimation of at least some components of $\hat{h}$ or $\tilde{h}$. Additionally, an examination of the behaviour of $\hat{B}$ for finite sample sizes would be worthwhile, since the usage of $\tilde{B}$ seems to be accompanied with a small bias. Another research aspect might consist in a comparison to more restrictive estimators like those of \citet{CKC2015} or \citet{ZLJ2009}.

	\begin{appendix}

		\section{Assumptions}\label{assumptions}
		
		\subsection{Assumptions Ensuring Identification of the Model}
		The following assumptions are necessary to identify model (\ref{modeleq}) and are taken from \citet{Klo2020}.
		\begin{enumerate}[label=(\textbf{A\arabic{*}})]
			\item\label{A1} Let $Y,\varepsilon$ and $\bX$ be real valued and $\mathbb{R}^{\dX}$-valued random variables, respectively, with
			$$h(Y)=g(\bX)+\sigma(\bX)\varepsilon$$
			for some transformation, regression and variance functions $h,g$ and $\sigma^2$.
			\item\label{A2} $\varepsilon$ is a centred random variable independent of $\bX$ with $E[\varepsilon]=0$ and $\operatorname{Var}(\varepsilon)=1$.
			\item\label{A3} Let the density $f_{\varepsilon}$ of $\varepsilon$ be continuous and let $h,g$ and $\sigma$ from \ref{A1} be continuously differentiable.
			\item\label{A4} The conditional cumulative distribution function $(y,\bx)\mapsto F_{Y|\bX}(y|\bx)$ is continuously differentiable with respect to $y$ and $\bx$. Let $j\in\{1,...,\dX\},v\geq0$ be an index and a weight function with support $\operatorname{supp}(v)$ such that $\frac{\partial}{\partial y}F_{Y|\bX}(y|\bx)>0$ for all $y\in\mathbb{R},\bx\in\operatorname{supp}(v)$ and such that (with $g$ and $\sigma$ from \ref{A1})
			$$A:=\int v(\bx)\left(\frac{\sigma(\bx)\frac{\partial g(\bx)}{\partial \bx_j}-g(\bx)\frac{\partial\sigma(\bx)}{\partial \bx_j}}{\sigma(\bx)}\right)\,d\bx\quad\textup{and}\quad B:=\int v(\bx)\frac{\frac{\partial\sigma(\bx)}{\partial \bx_j}}{\sigma(\bx)}\,d\bx$$
			are well defined with $B\neq0$.
		\end{enumerate}

		\subsection{Assumptions for the Convergence Results}
		Let $m\in\mathbb{N}$ and let $v$ be a weight function with a compact support.
		\begin{enumerate}[label=(\textbf{B\arabic{*}})]
			\item\label{B1} Let $(Y,X),(Y_1,X_1),...,(Y_n,X_n)$ be independent and identically distributed observations from model (\ref{modeleq}). Let the density $f_{Y,X}$ of the joint distribution of $(Y,X)$ be $(m+1)$-times continuously differentiable. Assume $f_{Y,X}$ to be bounded and $f_X$ to be bounded away from zero on the support of $v$.
			\item\label{B2} Let $K$ be a twice continuously differentiable kernel of order $m$ with compact support.
			\item\label{B3} Let $\sqrt{n}h_y^m\rightarrow0,\sqrt{n}h_x^m\rightarrow0,\frac{nh_y^5h_x}{\log(n)}\rightarrow\infty,\frac{nh_x^{\dX+4}}{\log(n)}\rightarrow\infty$ and $\frac{nh_x^3h_y^3}{\log(n)}\rightarrow\infty$.
			\item\label{B4} Let $v$ be $(m+1)$-times continuously differentiable.
			\item\label{B5} Let there exist some $c>0$ such that the function $\bx\mapsto\frac{\sigma(\bx)}{g(\bx)}$ is not almost surely constant on $M_{>c}=\{\bx:g(\bx)>c\}\,\cap\,\operatorname{supp}(v)$ or $M_{<-c}=\{\bx:g(\bx)<-c\}\,\cap\,\operatorname{supp}(v)$.
		\end{enumerate}

		\subsection{Assumptions Required for the Set $M_{\bX}$}
		Let $[z_a,z_b]\subseteq(y_0,\infty),[e_a,e_b]\subseteq\mathbb{R}$ be fixed and let $\tau<\beta$. The following restrictions of the compact set $M_{\bX}$ are needed to justify some technical calculations in the proof. It can be shown that a set satisfying these assumptions can be constructed from the data using some of the previous assumptions \ref{A1}--\ref{A4}, \ref{B1}--\ref{B5} \citep{Klo2019}.
		\begin{enumerate}[label=(\textbf{M\arabic{*}})]
			\item\label{M1} $M_{\bX}\subseteq\operatorname{supp}(v)$ and $f_{\bX}(\bx)>0$ for all $\bx\in M_{\bX}$,
			\item\label{M2} $\bx\mapsto\frac{g(\bx)}{\sigma(\bx)}$ is not almost surely constant on $M_{\bX}$,
			\item\label{M3} $F_{Y|\bX}^{-1}(\tau|\bx),F_{Y|\bX}^{-1}(\beta|\bx)\in(z_a,z_b)$ for all $\bx\in M_{\bX}$,
			\item\label{M4} $\underset{\bx\in M_{\bX},e\in[e_a,e_b],c\in[B_1,B_2]}{\sup}\,h_c(F_{Y|\bX}^{-1}(\tau|\bx))+e(h_c(F_{Y|\bX}^{-1}(\beta|\bx))-h_c(F_{Y|\bX}^{-1}(\tau|\bx)))<h_c(z_b)$
			\item\label{M5} $\underset{\bx\in M_{\bX},e\in[e_a,e_b],c\in[B_1,B_2]}{\inf}\,h_c(F_{Y|\bX}^{-1}(\tau|\bx))+e(h_c(F_{Y|\bX}^{-1}(\beta|\bx))-h_c(F_{Y|\bX}^{-1}(\tau|\bx)))>h_c(z_a)$.
		\end{enumerate}
		Since $M_{\bX}$ is an interval, the boundary of $M_{\bX}$ has Lebesgue-measure equal to zero.

		\section{Proof of the Main Result}
		In the following, the proofs of Lemma \ref{asympintegral} and Theorem \ref{asymph} are given, while the proofs of the auxiliary Lemmas can be found in Appendix \ref{auxiliary}.

		\subsection{Proof of Theorem \ref{asymph}}\label{proofasymph}
		(i) is proven first. In Appendix \ref{auxiliary}, it will be shown that \ref{C1}--\ref{C6} hold. Especially, one has
		$$\sqrt{n}(\hat{B}-B)=\frac{1}{\sqrt{n}}\sum_{i=1}^n\psi_{\Gamma_2}(Y_i,\bX_i)+o_p\bigg(\frac{1}{\sqrt{n}}\bigg).$$
		Further, apply a Taylor expansion to (\ref{defhc}) to obtain
		\begin{align*}
		\hat{h}_c(y)-h_c(y)&=\exp\bigg(-c\int_{y_1}^{y}\frac{1}{\hat{\lambda}(u)}\,du\bigg)-\exp\bigg(-c\int_{y_1}^{y}\frac{1}{\lambda(u)}\,du\bigg)
		\\[0,2cm]&=-c\exp\bigg(-c\int_{y_1}^{y}\frac{1}{\lambda(u)}\,du\bigg)\bigg(\int_{y_1}^{y}\frac{1}{\hat{\lambda}(u)}\,du-\int_{y_1}^{y}\frac{1}{\lambda(u)}\,du\bigg)
		\\[0,2cm]&\quad+o_p\bigg(\int_{y_1}^{y}\frac{1}{\hat{\lambda}(u)}\,du-\int_{y_1}^{y}\frac{1}{\lambda(u)}\,du\bigg)
		\\[0,2cm]&=-\frac{ch_c(y)}{n}\sum_{i=1}^n\eta_i(y)+o_p\bigg(\frac{1}{\sqrt{n}}\bigg)
		\end{align*}
		with $\eta_i,i=1,...,n,$ as in (\ref{defeta}), where the last equation is shown in the proof of Lemma \ref{asympintegral}. Moreover,
		\begin{align*}
		\underset{y\in\mathcal{K}}{\sup}\,|\hat{h}_1(y)-h_1(y)|^2&=o_p\bigg(\frac{1}{\sqrt{n}}\bigg),
		\\[0,2cm]\underset{y\in\mathcal{K}}{\sup}\,|\log(\hat{h}_1(y))-\log(h_1(y))|^2&=o_p\bigg(\frac{1}{\sqrt{n}}\bigg).
		\end{align*}
		Replacing $c$ by the estimator $\hat{B}$ and applying \ref{C6} results in
		\begin{align}
		H_n(y)&=\sqrt{n}(\hat{h}_{\hat{B}}(y)-h_{\hat{B}}(y)+h_{\hat{B}}(y)-h_B(y))\nonumber
		\\[0,2cm]&=-\frac{Bh(y)}{\sqrt{n}}\sum_{i=1}^n\eta_i(y)+\sqrt{n}\Big(\exp\big(\hat{B}\log(h_1(y))\big)-\exp\big(B\log(h_1(y))\big)\Big)+o_p(1)\nonumber
		\\[0,2cm]&=-\frac{Bh(y)}{\sqrt{n}}\sum_{i=1}^n\eta_i(y)+\exp\big(B\log(h_1(y))\big)\log(h_1(y))\sqrt{n}(\hat{B}-B)+o_p(1)\nonumber
		\\[0,2cm]&=\frac{h(y)}{\sqrt{n}}\sum_{i=1}^n\bigg(B\eta_i(y)+\log(h_1(y))\psi_{\Gamma_2}(Y_i,\bX_i)\bigg).\label{exprHn}
		\end{align}
		Convergence of the finite dimensional distributions follows from the Central Limit Theorem. Since $h_1$ is continuous and bounded away from zero on $\mathcal{K}$, asymptotic equicontinuity as defined in condition (2.1.8) of \citet{vdVW1996} is implied by that of $\big(\sqrt{n}(\hat{h}_1(y)-h(y))\big)_{y\in\mathcal{K}}$. Hence, Corollary 2.3.12 of \citet{vdVW1996} yields
		$$(H_n(y))_{y\in\mathcal{K}}\rightsquigarrow(Z_h(y))_{y\in\mathcal{K}},$$
		where the stated expression for the covariance function directly follows from (\ref{exprHn}).
		
		To prove (ii), write by the same reasoning as above
		\begin{align*}
		\tilde{h}(y)-h(y)&=Bh(y)\frac{\tilde{h}_1(y)-h_1(y)}{h_1(y)}+\exp\big(B\log(h_1(y))\big)\log(h_1(y))(\tilde{B}-B)
		\\[0,2cm]&\quad+o_p\bigg(\frac{1}{\sqrt{n}}+|\tilde{B}-B|\bigg)
		\\[0,2cm]&=-h(y)\int_{y_1}^{y}\frac{1}{\lambda(u)}\,du(\tilde{B}-B)+o_p\Bigg(\frac{1}{\sqrt{nh_y^3}}\Bigg)
		\end{align*}
		uniformly on compact sets $\mathcal{K}\subseteq(y_0,\infty)$. The weak convergence follows from Theorem \ref{asympaltB}. The proof of the second part of (ii) can be found in \citet{Klo2019}.

		\section{Proof of the Auxiliary Lemmas}\label{auxiliary}
		Here, the proofs of the auxiliary Lemmas \ref{asympintegral}, \ref{asympB} and \ref{asympaltB} are given.
		
		\subsection{Proof of Lemma \ref{asympintegral}}\label{proofasympintegral}
		First, the expansion in (\ref{uniconvlambda}) is derived before the weak convergence result is proven.\\[0,2cm]
		\textbf{Proof of the first assertion:}
		
		Uniform convergence results like
		$$\underset{\bx\in\operatorname{supp}(v)}{\sup}\,v(\bx)|\hat{f}(\bx)-f(\bx)|=o_p\big(n^{-\frac{1}{4}}\big)$$
		and
		$$\underset{y\in\mathcal{K},\,\bx\in\operatorname{supp}(v)}{\sup}\,|\hat{p}_y(y,\bx)-p_y(y,\bx)|=o_p\big(n^{-\frac{1}{4}}\big)$$
		can be shown similarly to \citet{Han2008}. Therefore, the same reasoning as in \citet{CKC2015} leads to
		\begin{align*}
		\frac{\hat{\Phi}_{x}(y,\bx)}{\hat{\Phi}_y(y,\bx)}-\frac{\Phi_{x}(y,\bx)}{\Phi_y(y,\bx)}&=(\hat{p}(y,\bx)-p(y,\bx))D_{p,0}(y,\bx)+(\hat{p}_x(y,\bx)-p_x(y,\bx))D_{p,x}(y,\bx)
		\\[0,2cm]&\quad+(\hat{p}_y(y,\bx)-p_y(y,\bx))D_{p,y}(y,\bx)+(\hat{f}(\bx)-f(\bx))D_{f,0}(y,\bx)
		\\[0,2cm]&\quad+(\hat{f}_x(\bx)-f_x(\bx))D_{f,x}(y,\bx)+o_p\bigg(\frac{1}{\sqrt{n}}\bigg)
		\end{align*}
		so that one obtains
		\begin{align*}
		&\hat{\lambda}(u)-\lambda(u)
		\\[0,2cm]&=\frac{1}{n}\sum_{i=1}^n\int\bigg({D}_{p,0}(u,\bx)\mathcal{K}_{h_y}(u-Y_i)\bK_{h_x}(\bx-\bX_i)
		\\[0,2cm]&\quad+{D}_{p,x}(u,\bx)\mathcal{K}_{h_y}(u-Y_i)\frac{\partial \bK_{h_x}(\bx-\bX_i)}{\partial \bx_j}+{D}_{p,y}(u,\bx)K_{h_y}(u-Y_i)\bK_{h_x}(\bx-\bX_i)
		\\[0,2cm]&\quad+{D}_{f,0}(u,\bx)\bK_{h_x}(\bx-\bX_i)+{D}_{f,x}(u,\bx)\frac{\partial \bK_{h_x}(\bx-\bX_i)}{\partial \bx_j}\bigg)v(\bx)\,d\bx+o_p\bigg(\frac{1}{\sqrt{n}}\bigg)
		\end{align*}
		by inserting the definition of $\hat{p},\hat{p}_x,\hat{p}_y,\hat{f},\hat{f}_x$. Due to the assumptions \ref{B2} and \ref{B3}, a Taylor expansion similarly to \citet{CvK2018} leads to
		$$\int l(\bx)\bK_{h_x}(\bx-\bX_i)\,d\bx=\int l(\bX_i+h_x\bx)\bK(\bx)\,d\bx=l(\bX_i)+o\bigg(\frac{1}{\sqrt{n}}\bigg)$$
		for every $m$-times continuously differentiable function $l$ with bounded support. Moreover, integration by parts yields
		\begin{align*}
		-\int l(\bx)\frac{\partial}{\partial \bx_j}\bK_{h_x}(\bx-\bX_i)\,d\bx&=\int \bK_{h_x}(\bx-\bX_i)\frac{\partial}{\partial \bx_j}l(\bx)\,d\bx
		\\[0,2cm]&=\int \bK(\bx)\frac{\partial}{\partial \bx_j}l(\bx)\bigg|_{\bx=\bX_i+h_x\bx}\,d\bx
		\\[0,2cm]&=\frac{\partial}{\partial \bx_j}l(\bx)\bigg|_{\bx=\bX_i}+o\bigg(\frac{1}{\sqrt{n}}\bigg)
		\end{align*}
		for every $(m+1)$-times continuously differentiable function $l$ with bounded support.  Due to the compactness of $\operatorname{supp}(v)$ and $\mathcal{K}$, all derivatives of $D_{p,0},...,D_{f,x}$ are bounded, so that
		\begin{align*}
		\\[0,2cm]\hat{\lambda}(u)-\lambda(u)&=\frac{1}{n}\sum_{i=1}^n\bigg(v(\bX_i){D}_{p,0}(u,\bX_i)\mathcal{K}_{h_y}(u-Y_i)-\frac{\partial v(\bX_i){D}_{p,x}(u,\bX_i)}{\partial \bx_j}\mathcal{K}_{h_y}(u-Y_i)
		\\[0,2cm]&\quad+v(\bX_i){D}_{p,y}(u,\bX_i)K_{h_y}(u-Y_i)+v(\bX_i){D}_{f,0}(u,\bX_i)
		\\[0,2cm]&\quad-\frac{\partial v(\bX_i){D}_{f,x}(u,\bX_i)}{\partial \bx_j}\bigg)+o_p\bigg(\frac{1}{\sqrt{n}}\bigg).
		\end{align*}
		Finally,
		$$\underset{y\in\mathcal{K}}{\sup}\,|\hat{\lambda}(u)-\lambda(u)|=\mathcal{O}_p\Bigg(\sqrt{\frac{\log(n)}{nh_y}}\Bigg)$$
		follows again as in \citet{Han2008}.\\[0,2cm]
		\textbf{Proof of the weak convergence:}
		
		The main idea to prove the second assertion is to find an expression
		\begin{equation}\label{asympeta}
		\int_{y_1}^y\bigg(\frac{1}{\hat{\lambda}(u)}-\frac{1}{\lambda(u)}\bigg)\,du=\frac{1}{n}\sum_{i=1}^n(\eta_i(y)-E[\eta_i(y)])+o_p\bigg(\frac{1}{\sqrt{n}}\bigg)
		\end{equation}
		($\eta_i$ will be defined later) for which some weak convergence results can be applied.
		
		First, remark that since $\mathcal{K}\subseteq(y_0,\infty)$ is compact, $u\mapsto\frac{1}{\lambda(u)}$ is bounded and bounded away from zero on $\mathcal{K}$. Hence, one has
		\begin{align*}
		\int_{y_1}^y\bigg(\frac{1}{\hat{\lambda}(u)}-\frac{1}{\lambda(u)}\bigg)\,du&=\int_{y_1}^y\frac{\lambda(u)-\hat{\lambda}(u)}{\hat{\lambda}(u)\lambda(u)}\,du
		\\[0,2cm]&=\int_{y_1}^y\frac{\lambda(u)-\hat{\lambda}(u)}{\lambda(u)^2}\bigg(1-\frac{\hat{\lambda}(u)-\lambda(u)}{\hat{\lambda}(u)}\bigg)\,du.
		\end{align*}
		Possibly, extend $\mathcal{K}$ such that $y_1$ is included. Due to assumption \ref{B3}, the expansion in (\ref{uniconvlambda}) leads to
		\begin{align*}
		Z_n(u)&\,\ =\sqrt{n}\Bigg(\int_{y_1}^y\frac{\lambda(u)-\hat{\lambda}(u)}{\lambda(u)^2}\,du+\mathcal{O}_p\Big(\underset{u\in\mathcal{K}}{\sup}\,|\hat{\lambda}(u)-\lambda(u)|^2\Big)\Bigg)
		\\[0,2cm]&\overset{(\ref{uniconvlambda})}{=}\sqrt{n}\bigg(\int_{y_1}^y\frac{\lambda(u)-\hat{\lambda}(u)}{\lambda(u)^2}\,du\bigg)+o_p(1)
		\end{align*}
		as well as
		\begin{align*}
		\int_{y_1}^y\frac{\lambda(u)-\hat{\lambda}(u)}{\lambda(u)^2}\,du&=\frac{1}{n}\sum_{i=1}^n\int_{y_1}^y\frac{-1}{\lambda(u)^2}\bigg(v(\bX_i){D}_{p,0}(u,\bX_i)\mathcal{K}_{h_y}(u-Y_i)
		\\[0,2cm]&\quad-\frac{\partial v(\bX_i){D}_{p,x}(u,\bX_i)}{\partial \bx_j}\mathcal{K}_{h_y}(u-Y_i)+v(\bX_i){D}_{p,y}(u,\bX_i)K_{h_y}(u-Y_i)
		\\[0,2cm]&\quad+v(\bX_i){D}_{f,0}(u,\bX_i)-\frac{\partial v(\bX_i){D}_{f,x}(u,\bX_i)}{\partial \bx_j}\bigg)\,du+o_p\bigg(\frac{1}{\sqrt{n}}\bigg)
		\\[0,2cm]&=:\frac{1}{n}\sum_{i=1}^n\tilde{\eta}_i(y)+o_p\bigg(\frac{1}{\sqrt{n}}\bigg).
		\end{align*}
		The following lemma is similar to Proposition 2 of \citet{CvK2018}. The proof is omitted here and can be found in \citet{Klo2019}.
		\begin{lemma}\label{lemmacolling}
			Let $\mathcal{K}\subseteq(y_0,\infty)$ be compact, $l:\mathbb{R}\times\mathbb{R}^{\dX}\rightarrow\mathbb{R},(u,\bx)\mapsto l(u,\bx),$ be bounded on compact sets and let $l$ have a compact support with respect to the $\bx$-component, which will be denoted by $\operatorname{supp}_\bx(l)$ in the following. Then, under the conditions of Lemma \ref{asympintegral} one has
			$$\frac{1}{n}\sum_{i=1}^n\int_{y_1}^yl(u,\bX_i)\big(\mathcal{K}_{h_y}(u-Y_i)-I_{\{Y_i\leq u\}}\big)\,du=o_p\bigg(\frac{1}{\sqrt{n}}\bigg)$$
			uniformly in $y\in\mathcal{K}$.
		\end{lemma}
		
		Define
		\begin{align}
		\eta_i(y)&:=\int_{y_1}^y\frac{-1}{\lambda(u)^2}\bigg(v(\bX_i){D}_{p,0}(u,\bX_i)-\frac{\partial v(\bX_i){D}_{p,x}(u,\bX_i)}{\partial \bx_j}\bigg)I_{\{u\geq Y_i\}}\,du\nonumber
		\\[0,2cm]&\quad-\frac{v(\bX_i){D}_{p,y}(Y_i,\bX_i)}{\lambda(Y_i)^2}\big(I_{\{Y_i\leq y\}}-I_{\{Y_i\leq y_1\}}\big)\nonumber
		\\[0,2cm]&\quad+\int_{y_1}^y\frac{-1}{\lambda(u)^2}\bigg(v(\bX_i){D}_{f,0}(u,\bX_i)-\frac{\partial v(\bX_i){D}_{f,x}(u,\bX_i)}{\partial \bx_j}\bigg)\,du.\label{defeta}
		\end{align}
		Then, Lemma \ref{lemmacolling} leads to
		\begin{align*}
		\frac{1}{n}\sum_{i=1}^n\tilde{\eta}_i(y)&=\frac{1}{n}\sum_{i=1}^n\int_{y_1}^y\frac{-1}{\lambda(u)^2}\bigg(v(\bX_i){D}_{p,0}(u,\bX_i)I_{\{u\geq Y_i\}}-\frac{\partial v(\bX_i){D}_{p,x}(u,\bX_i)}{\partial \bx_j}I_{\{u\geq Y_i\}}
		\\[0,2cm]&\quad+v(\bX_i){D}_{p,y}(u,\bX_i)K_{h_y}(u-Y_i)+v(\bX_i){D}_{f,0}(u,\bX_i)
		\\[0,2cm]&\quad-\frac{\partial v(\bX_i){D}_{f,x}(u,\bX_i)}{\partial \bx_j}\bigg)\,du+o_p\bigg(\frac{1}{\sqrt{n}}\bigg)
		\\[0,2cm]&=\frac{1}{n}\sum_{i=1}^n{\eta}_i(y)+o_p\bigg(\frac{1}{\sqrt{n}}\bigg),
		\end{align*}
		where
		\begin{align*}
		&\int_{y_1}^{y}\frac{-1}{\lambda(u)^2}v(\bX_i){D}_{p,y}(u,\bX_i)K_{h_y}(u-Y_i)\,du
		\\[0,2cm]&=-\int\frac{1}{\lambda(Y_i+h_yu)^2}v(\bX_i)D_{p,y}(Y_i+h_yu)K(u)\big(I_{\{Y_i\leq y-h_yu\}}-I_{\{Y_i\leq y_1-h_yu\}}\big)\,du
		\\[0,2cm]&=-\frac{v(\bX_i){D}_{p,y}(Y_i,\bX_i)}{\lambda(Y_i)^2}\big(I_{\{Y_i\leq y\}}-I_{\{Y_i\leq y_1\}}\big)+o_p\bigg(\frac{1}{\sqrt{n}}\bigg)
		\end{align*}
		can be shown similarly to Lemma \ref{lemmacolling}. Straightforward calculations lead to $E[\eta_i(y)]=o_p\Big(\frac{1}{\sqrt{n}}\Big)$ uniformly in $y\in\mathcal{K}$, so that equation (\ref{asympeta}) is valid. So far, the asymptotic representation
		$$\int_{y_1}^y\bigg(\frac{1}{\hat{\lambda}(u)}-\frac{1}{\lambda(u)}\bigg)\,du=\frac{1}{n}\sum_{i=1}^n(\eta_i(y)-E[\eta_i(y)])+o_p\bigg(\frac{1}{\sqrt{n}}\bigg)$$
		was proven. It remains to show weak convergence of the corresponding process to an appropriate Gaussian process. For this purpose, define
		\begin{align*}
		\eta_{z,\bx}^a(y)&:=\int_{y_1}^y\frac{-1}{\lambda(u)^2}\bigg(v(\bx){D}_{p,0}(u,\bx)-\frac{\partial v(\bx){D}_{p,x}(u,\bx)}{\partial \bx_j}\bigg)I_{\{u\geq z\}}\,du
		\\[0,2cm]&\quad+\int_{y_1}^y\frac{-1}{\lambda(u)^2}\bigg(v(\bx){D}_{f,0}(u,\bx)-\frac{\partial v(\bx){D}_{f,x}(u,\bx)}{\partial \bx_j}\bigg)\,du,
		\\[0,2cm]\eta_{z,\bx}^b(y)&:=-\bigg(\frac{v(\bx){D}_{p,y}(z,\bx)}{\lambda(z)^2}\bigg)_+\big(I_{\{z\leq y\}}-I_{\{z\leq y_1\}}\big),
		\\[0,2cm]\eta_{z,\bx}^c(y)&:=\bigg(\frac{v(\bx){D}_{p,y}(z,\bx)}{\lambda(z)^2}\bigg)_{-}\big(I_{\{z\leq y\}}-I_{\{z\leq y_1\}}\big),
		\end{align*}
		where for some value $a\in\mathbb{R}$ the terms $(a)_+$ and $(a)_-$ denote the positive and negative part of $a$, respectively. Hence,
		$$\eta_i(y)=\eta_{Y_i,\bX_i}^a(y)+\eta_{Y_i,\bX_i}^b(y)+\eta_{Y_i,\bX_i}^c(y).$$
		It can be easily seen that $\eta_{z,\bx}^a(y),\eta_{z,\bx}^b(y)$ and $\eta_{z,\bx}^c(y)$ are bounded by some constant $\tilde{C}>0$ uniformly in $y,\tilde{y}\in\mathcal{K}$. In the following, it will be proven, that the function classes
		$$\mathcal{F}^k:=\big\{(z,\bx)\mapsto\eta_{z,\bx}^k(y),y\in\mathcal{K}\big\},\quad k\in\{a,b,c\},$$
		are Donsker. Example 2.10.7 of \citet{vdVW1996} then implies that the class $\mathcal{F}=\{(z,\bx)\mapsto\eta_{z,\bx}(y),y\in\mathcal{K}\}$ is Donsker as well. While the Donsker property of $\mathcal{F}^b$ and $\mathcal{F}^c$ can be shown by standard arguments as for indicator functions, one has
		\begin{align*}
		|\eta_{z,\bx}^a(y)-\eta_{z,\bx}^a(\tilde{y})|&=\bigg|\int_{\tilde{y}}^y\frac{-1}{\lambda(u)^2}\bigg(v(\bx){D}_{p,0}(u,\bx)+\frac{\partial v(\bx){D}_{p,x}(u,\bx)}{\partial \bx_j}\bigg)I_{\{u\geq z\}}\,du
		\\[0,2cm]&\quad+\int_{\tilde{y}}^y\frac{-1}{\lambda(u)^2}\bigg(v(\bx){D}_{f,0}(u,\bx)+\frac{\partial v(\bx){D}_{f,x}(u,\bx)}{\partial \bx_j}\bigg)\,du\bigg|
		\\[0,2cm]&\leq C|y-\tilde{y}|
		\end{align*}
		for all $y,\tilde{y}\in\mathcal{K}$ and an appropriate constant $C>0$, so that
		\begin{align*}
		\sqrt{E[(\eta_{Z_1,\bX_1}^a(y)-\eta_{Z_1,\bX_1}^a(\tilde{y}))^2]}\leq C|\tilde{y}-y|.
		\end{align*}
		Let $\xi>0$. Then, $\xi$-brackets $[l,u]$ for the function class $\mathcal{F}^a$ can be defined as
		$$l(z,\bx)=\eta_{z,\bx}^a(y_k^*)-\frac{\sqrt{\xi}}{C}\quad\textup{and}\quad u(z,\bx)=\eta_{z,\bx}^a(y_k^*)+\frac{\sqrt{\xi}}{C},\quad k=1,...,K,$$
		for some $K\in\mathbb{N}$ and appropriate values $y_1^*,...,y_K^*\in\mathcal{K}$. Consequently, the bracketing number can be deduced from that of $\mathcal{K}$ and for some constant $C$ the bracketing integral
		$$\int_0^{\infty}\sqrt{\log(\mathcal{N}_{[\,]}(\varepsilon,\mathcal{F},L_2(P^{Y,\bX})))}\,d\varepsilon=C\int_0^{\infty}\sqrt{\log\bigg(\max\bigg(\frac{1}{\varepsilon^2},1\bigg)\bigg)}\,d\varepsilon<\infty$$
		is finite. Theorem 2.5.6 of \citet{vdVW1996} ensures that $\mathcal{F}^a$ is Donsker, as long as the finite dimensional distributions converge, but this in turn (as for $\mathcal{F}^b,\mathcal{F}^c$ and $\mathcal{F}$) is implied by the multivariate Central Limit Theorem. After some rather technical computations for the indicator functions, the covariance function can be written as in \citet{KNVK2019}. Finally, the weak convergence
		$$(Z_n(y))_{y\in\mathcal{K}}\rightsquigarrow(Z(y))_{y\in\mathcal{K}}$$
		was proven, where $Z$ is a centred Gaussian process.\hfill$\square$

		\subsection{Proof of Lemma \ref{asympB}}\label{proofasympB}
		Before Lemma \ref{asympB} can be proven, some further notations are needed. These notations are used in the proof of some preliminary lemmas, which in turn will imply the assertion of Lemma \ref{asympB}.

		\subsubsection{Preliminary Notations}
		Denote the supremum norms of the functions $f_{m_{\tau}},f_{m_{\beta}}:\mathbb{R}^{\dX}\rightarrow\mathbb{R}$ and $h:\mathbb{R}\rightarrow\mathbb{R}$ on $M_{\bX}$ and $[z_a,z_b]$ by $||.||_{M_{\bX}}$ and $||.||_{[z_a,z_b]}$, respectively. Let $C>0$ such that $\underset{u\in[z_a,z_b]}{\sup}\,\big|\frac{\partial^2}{\partial u^2}h_1(u)\big|\leq\frac{C}{2}$ and define the set of functions
		\begin{align}
		\mathcal{H}&=\bigg\{s=(\mathfrak{h},f_{m_{\tau}},f_{m_{\beta}})^t:\mathfrak{h}\in\mathcal{C}^2([z_a,z_b]),f_{m_{\tau}},f_{m_{\beta}}\in\mathcal{C}^2(M_{\bX}),f_{m_{\tau}}(M_{\bX})\subseteq(z_a,z_b),\nonumber
		\\[0,2cm]&\quad\quad f_{m_{\beta}}(M_{\bX})\subseteq(z_a,z_b),\bigg|\frac{\partial^2}{\partial u^2}\mathfrak{h}(u)\bigg|\leq C,2\underset{u\in[z_a,z_b]}{\inf}\,\frac{\partial}{\partial u}\mathfrak{h}(u)>\underset{u\in[z_a,z_b]}{\inf}\,\frac{\partial}{\partial u}h_1(u)\bigg\}\label{defHcal}
		\end{align}
		endowed with the supremum norm
		$$||s||_{\mathcal{H}}=\max\big(||h||_{[z_a,z_b]},||f_{m_{\tau}}||_{M_{\bX}},||f_{m_{\beta}}||_{M_{\bX}}\big).$$
		Following Section 2.7.1 of \citet{vdVW1996}, consider for some $\gamma,R>0$ the (H\"older-)class $C_R^{\gamma}$ of all functions on $M_{\bX}$ such that all partial derivatives up to order $\lfloor\gamma\rfloor$ are uniformly bounded by $R$ and the partial derivatives of highest order are Lipschitz of order $\gamma-\lfloor\gamma\rfloor$. More precisely, define for any multi-index $l=(l_1,...,l_{\dX})$ the differential
		$$D_l=\frac{\partial^l}{\partial \bx_1^{l_1}...\partial \bx_d^{l_d}}$$
		as well as the norm
		$$||f||_{\gamma}=\underset{l\leq\lfloor\gamma\rfloor}{\max}\,\underset{\bx\in M_{\bX}}{\sup}\,|D_lf(\bx)|+\underset{l=\lfloor\gamma\rfloor}{\max}\,\underset{\bx\neq y\in M_{\bX}}{\sup}\,\frac{|D_lf(\bx)-D_lf(y)|}{||\bx-y||^{\gamma-\lfloor\gamma\rfloor}},$$
		where the inequality $l\leq\lfloor\gamma\rfloor$ has to be read in the sense of $\sum_{i=1}^{\dX}l_i\leq\lfloor\gamma\rfloor$ for every multi-index $l=(l_1,...,l_{\dX})$. In the case of $\dX=1$ the norm can be written as
		$$||f||_{\gamma}=\underset{l=1,...,\lfloor\gamma\rfloor}{\max}\,\underset{\bx\in M_{\bX}}{\sup}\,\bigg|\frac{\partial^l}{\partial \bx^l}f(\bx)\bigg|+\underset{\bx\neq y\in M_{\bX}}{\sup}\,\frac{\big|\frac{\partial^{\lfloor\gamma\rfloor}}{\partial \bx^{\lfloor\gamma\rfloor}}f(\bx)-\frac{\partial^{\lfloor\gamma\rfloor}}{\partial \bx^{\lfloor\gamma\rfloor}}f(y)\big|}{||\bx-y||^{\gamma-\lfloor\gamma\rfloor}}.$$
		Further, define for some $R>0$ the set $C_R^{\gamma}(M_{\bX})$ as the set of all (sufficiently often differentiable) functions $f$ with $||f||_{\gamma}\leq R$ and
		\begin{equation}\label{defHcaltilde}
		\tilde{\mathcal{H}}=\bigg\{s\in\mathcal{H}: h\in C_{R_h}^{\gamma_h}([z_a,z_b]),f_{m_{\tau}}\in C_{R_{f_{m_{\tau}}}}^{\gamma_{f_{m_{\tau}}}}(M_{\bX}),f_{m_{\beta}}\in C_{R_{f_{m_{\beta}}}}^{\gamma_{f_{m_{\beta}}}}(M_{\bX})\bigg\}
		\end{equation}
		for some constants $\gamma_h>1,\gamma_{f_{m_{\tau}}},\gamma_{f_{m_{\beta}}}>\dX$ and $R_h,R_{f_{m_{\tau}}},R_{f_{m_{\beta}}}<\infty$.
		
		Recall that $\Gamma_1(c,s_0)$ was defined as the ordinary derivative of $G_{MD}$ with respect to $c\in[B_1,B_2]$. For all $(\bx,e)\in M_{\bX}\times[e_a,e_b]$ let $\Gamma_2(c,s_0)(\bx,e)[s-s_0]$ denote the directional derivative of $G_{MD}(c,s_0)(\bx,e)$ with respect to $s$, that is
		$$\Gamma_2(c,s_0)(\bx,e)[s-s_0]:=\underset{t\rightarrow0}{\lim}\,\frac{G_{MD}(c,s_0+t(s-s_0))(\bx,e)-G_{MD}(c,s_0)(\bx,e)}{t}.$$
		Furthermore, let
		$$\begin{array}{c}
		D_hG_{MD}(c,s_0)(\bx,e)[\mathfrak{h}-h_1],
		\\[0,2cm]D_{f_{m_{\tau}}}G_{MD}(c,s_0)(\bx,e)\big[f_{m_{\tau}}-F_{Y|\bX}^{-1}(\tau|\cdot)\big],
		\\[0,2cm]D_{f_{m_{\beta}}}G_{MD}(c,s_0)(\bx,e)\big[f_{m_{\beta}}-F_{Y|\bX}^{-1}(\beta|\cdot)\big]
		\end{array}$$
		denote the directional derivatives with respect to $\mathfrak{h},f_{m_{\tau}}$ and $f_{m_{\beta}}$, respectively. Now, some further properties can be formulated. These will allow to proceed in the proof of Theorem \ref{asympB} below similarly to \citet{LSvK2008}. Let $G_{MD}$ and $G_{nMD}$ be defined as in (\ref{defGMD}) and (\ref{defGnMD}) and let $\hat{s}=(\hat{h}_1,\hat{F}_{Y|\bX}^{-1}(\tau|\cdot),\hat{F}_{Y|\bX}^{-1}(\beta|\cdot))$ be the estimator of $s_0$ from Section \ref{estB}. In the proof of Theorem \ref{asympB}, it will be shown first that the assumptions \ref{C1}--\ref{C6} below are already implied by \ref{A1}--\ref{A4}, \ref{B1}--\ref{B5} and \ref{M1}--\ref{M5} from Appendix \ref{assumptions}. In a second step, these assumptions are used to proof the assertion similarly to \citet{CLvK2003} and \citet{LSvK2008}.
		\begin{enumerate}[label=(\textbf{C\arabic{*}})]
			\item\label{C1} One has $G_{MD}(B,s_0)\equiv0$ and $\hat{B}-B=o_p(1)$.
			\item\label{C2} For all $(\bx,e)\in M_{\bX}\times[e_a,e_b]$ the ordinary derivative $\Gamma_1(c,s_0)(\bx,e)$ with respect to $c$ of $G_{MD}(c,s_0)(\bx,e)$ exists in a neighbourhood of $B$ and is continuous at $c=B$. $\Gamma_1(B,s_0)(\bx,e)$ is different from zero on a set with positive $\lambda_{M_{\bX}\times[e_a,e_b]}$-measure.
			\item\label{C3} For any $\delta>0$ let $B_{\delta}$ be the $\delta$-neighbourhood of $B$ in $[B_1,B_2]$ and $\tilde{\mathcal{H}}_{\delta}=\{s\in\tilde{\mathcal{H}}:||s-s_0||_{\mathcal{H}}<\delta\}$. With these notations, the directional derivative $\Gamma_2(c,s_0)(\bx,e)[s-s_0]$ of $G_{MD}(c,s_0)(\bx,e)$ with respect to $s$ exists for all $c\in B_{\delta},(\bx,e)\in M_{\bX}\times[e_a,e_b]$ and in all directions $[s-s_0]$ with $s\in\tilde{\mathcal{H}}$ and $\tilde{\mathcal{H}}$ as in (\ref{defHcaltilde}). Consider a positive sequence $\delta_n\rightarrow0$ and $(c,s)\in B_{\delta_n}\times\tilde{\mathcal{H}}_{\delta_n}$. Then,
			\begin{enumerate}[label=(\roman*)]
				\item for an appropriate constant $C\geq0$ (independent of $c$ and $s$) it holds that
				\begin{align*}
				&||G_{MD}(c,s)-G_{MD}(c,s_0)-\Gamma_2(c,s_0)[s-s_0]||_2
				\\[0,2cm]&\leq C\big(||h-h_1||_{[z_a,z_b]}^{\frac{3}{2}}+||f_{m_{\tau}}-F_{Y|\bX}^{-1}(\tau|\cdot)||_{M_{\bX}}^2+||f_{m_{\beta}}-F_{Y|\bX}^{-1}(\beta|\cdot)||_{M_{\bX}}^2\big).
				\end{align*}
				\item one has $||\Gamma_2(c,s_0)[\hat{s}-s_0]-\Gamma_2(B,s_0)[\hat{s}-s_0]||=o_p(|c-B|)+o_p\big(n^{-\frac{1}{2}}\big)$ uniformly in $c\in B_{\delta_n}$.
			\end{enumerate}
			\item\label{C4} It holds that $\hat{s}\in\tilde{\mathcal{H}}$ with probability converging to one, $||\hat{h}_1-h_1||_{[z_a,z_b]}^{\frac{3}{2}}=o_p(n^{-\frac{1}{2}})$ and
			$$||\hat{F}_{Y|\bX}^{-1}(\tau|\cdot)-F_{Y|\bX}^{-1}(\tau|\cdot)||_{M_{\bX}},||\hat{F}_{Y|\bX}^{-1}(\beta|\cdot)-F_{Y|\bX}^{-1}(\beta|\cdot)||_{M_{\bX}}=o_p(n^{-\frac{1}{4}}).$$
			\item\label{C5} $\underset{||c-B||\leq\delta_n,||s-s_0||\leq\delta_n}{\sup}\,||G_{nMD}(c,s)-G_{MD}(c,s)-G_{nMD}(B,s_0)||_2=o_p(n^{-\frac{1}{2}})$.
			\item\label{C6} There exists a real valued function $\psi_{\Gamma_2}$ with $E[\psi_{\Gamma_2}(Y,\bX)]=o\big(n^{-\frac{1}{2}}\big)$ and $\sigma_A^2:=E[\psi_{\Gamma_2}(Y,\bX)^2]\in(0,\infty)$ such that
			\begin{align*}
			&\sqrt{n}\int_{M_{\bX}}\int_{[e_a,e_b]}\Gamma_1(B,s_0)(\bx,e)\Big(G_{nMD}(B,s_0)(\bx,e)+\Gamma_2(B,s_0)(\bx,e)[\hat{s}-s_0]\Big)\,de\,d\bx
			\\[0,2cm]&=\quad\frac{1}{\sqrt{n}}\sum_{i=1}^n\psi_{\Gamma_2}(Y_i,\bX_i)+o_p(1)
			\\[0,2cm]&\overset{\mathcal{D}}{\rightarrow}\mathcal{N}(0,\sigma_A^2).
			\end{align*}
		\end{enumerate}
		Before the corresponding Lemmas are stated and proven define for $s\in\tilde{\mathcal{H}},c\in[B_1,B_2],$ $\bx\in M_{\bX},e\in[e_a,e_b]$
		\begin{equation}\label{defkc}
		k_c(s,\bx,e)=\frac{h_1\big(\mathfrak{h}^{-1}\big(\big(\mathfrak{h}^c(f_{m_{\tau}}(\bx))+e(\mathfrak{h}^c(f_{m_{\beta}}(\bx))-\mathfrak{h}^c(f_{m_{\tau}}(\bx)))\big)^{\frac{1}{c}}\big)\big)^B-g(\bx)}{\sigma(\bx)}
		\end{equation}
		with $h_1$ as in (\ref{defhc}), so that $h=h_1^B$ and
		\begin{align}
		&G_{MD}(c,s)(\bx,e)\nonumber
		\\[0,2cm]&=P(\bX\leq \bx,\tilde{\varepsilon}_c(s)\leq e|\bX\in M_{\bX})-P(\bX\leq \bx|\bX\in M_{\bX})P(\tilde{\varepsilon}_c(s)\leq e|\bX\in M_{\bX})\nonumber
		\\[0,2cm]&=P(\bX\leq \bx,\varepsilon\leq k_c(s,\bX,e)|\bX\in M_{\bX})\nonumber
		\\[0,2cm]&\quad-P(\bX\leq \bx|\bX\in M_{\bX})P(\varepsilon\leq k_c(s,\bX,e)|\bX\in M_{\bX})\nonumber
		\\[0,2cm]&=\frac{1}{P(\bX\in M_{\bX})}\int_{M_{\bX}}\big(I_{\{\mathbf{w}\leq \bx\}}-P(\bX\leq \bx|\bX\in M_{\bX})\big)f_\bX(\mathbf{w})F_{\varepsilon}(k_c(s,\mathbf{w},e))\,d\mathbf{w}
		.\label{exprGMDkc}
		\end{align}
		\hfill$\square$

		\subsubsection{Preliminary Lemmas}
		In the following, assume validity of \ref{A1}--\ref{A4}, \ref{B1}--\ref{B5} and \ref{M1}--\ref{M5} from Appendix \ref{assumptions} for a compact interval $[z_a,z_b]\subseteq(y_0,\infty)$.
		\begin{lemma}\label{lemmac4}
			If $\dX=1$, one has $\hat{s}\in\tilde{\mathcal{H}}$ with probability converging to one for $\gamma_h=\gamma_{f_{m_{\tau}}}=\gamma_{f_{m_{\beta}}}=2$ and some sufficiently large constants $R_h,R_{f_{m_{\tau}}},R_{f_{m_{\beta}}}>0$. Moreover, the second part of \ref{C4} is valid as well.
		\end{lemma}\noindent
		\textbf{Proof:} Recall the definition of $\tilde{\mathcal{H}}\subseteq\mathcal{H}$ (remember $\gamma_h=\gamma_{f_{m_{\tau}}}=\gamma_{f_{m_{\beta}}}=2$):
		$$\tilde{\mathcal{H}}=\bigg\{s\in\mathcal{H}: h\in C_{R_h}^{2}([z_a,z_b]),f_{m_{\tau}}\in C_{R_{f_{m_{\tau}}}}^{2}(M_{\bX}),f_{m_{\beta}}\in C_{R_{f_{m_{\beta}}}}^{2}(M_{\bX})\bigg\}.$$
		The convergence rate of $||\hat{h}_1-h_1||_{[z_a,z_b]}$ follows from Theorem \ref{asympintegral}, while those of $||\hat{F}_{Y|\bX}^{-1}(\tau|\cdot)-F_{Y|\bX}^{-1}(\tau|\cdot)||_{M_{\bX}}$ and $||\hat{F}_{Y|\bX}^{-1}(\beta|\cdot)-F_{Y|\bX}^{-1}(\beta|\cdot)||_{M_{\bX}}$ can be derived by writing
		\begin{align*}
		0&=\hat{F}_{Y|\bX}(\hat{F}_{Y|\bX}^{-1}(\tau|\bx)|\bx)-F_{Y|\bX}(F_{Y|\bX}^{-1}(\tau|\bx)|\bx)
		\\[0,2cm]&=\hat{F}_{Y|\bX}(F_{Y|\bX}^{-1}(\tau|\bx)|\bx)+\hat{f}_{Y|\bX}(\tilde{F}(\tau,\bx)|\bx)\big(\hat{F}_{Y|\bX}^{-1}(\tau|\bx)-F_{Y|\bX}^{-1}(\tau|\bx)\big)
		\\[0,2cm]&\quad-F_{Y|\bX}(F_{Y|\bX}^{-1}(\tau|\bx)|\bx)
		\\[0,2cm]&=\hat{F}_{Y|\bX}(F_{Y|\bX}^{-1}(\tau|\bx)|\bx)-F_{Y|\bX}(F_{Y|\bX}^{-1}(\tau|\bx)|\bx)
		\\[0,2cm]&\quad+f_{Y|\bX}(F_{Y|\bX}^{-1}(\tau|\bx)|\bx)\big(\hat{F}_{Y|\bX}^{-1}(\tau|\bx)-F_{Y|\bX}^{-1}(\tau|\bx)\big)+o_p\big(\hat{F}_{Y|\bX}^{-1}(\tau|\bx)-F_{Y|\bX}^{-1}(\tau|\bx)\big)
		\end{align*}
		uniformly in $\bx\in\operatorname{supp}(v)$, which in turn results in
		\begin{align}
		&\hat{F}_{Y|\bX}^{-1}(\tau|\bx)-F_{Y|\bX}^{-1}(\tau|\bx)\nonumber
		\\[0,2cm]&=-\frac{\hat{F}_{Y|\bX}(F_{Y|\bX}^{-1}(\tau|\bx)|\bx)-F_{Y|\bX}(F_{Y|\bX}^{-1}(\tau|\bx)|\bx)}{f_{Y|\bX}(F_{Y|\bX}^{-1}(\tau|\bx)|\bx)}(1+o_p(1))\nonumber
		\\[0,2cm]&=-\frac{(1+o_p(1))}{f_{Y|\bX}(F_{Y|\bX}^{-1}(\tau|\bx)|\bx)}\frac{1}{n}\sum_{i=1}^n\bigg(\frac{1}{f_{\bX}(\bx)}\mathcal{K}_{h_y}(F_{Y|\bX}^{-1}(\tau|\bx)-Y_i)\bK_{h_x}(\bx-\bX_i)\nonumber
		\\[0,2cm]&\quad-\frac{p(F_{Y|\bX}^{-1}(\tau|\bx),\bx)}{f_{\bX}(\bx)^2}\bK_{h_x}(\bx-\bX_i)\bigg)\label{FYXhatinverse}
		\end{align}
		uniformly in $\bx\in\operatorname{supp}(v)$. Now, similar techniques as in \citet{Han2008} can be applied to obtain the result.
		
		To prove $\hat{s}\in\tilde{\mathcal{H}}$ with probability converging to one it suffices to show uniform convergence of the functions $\hat{h}_1,\hat{F}_{Y|\bX}^{-1}(\tau|\cdot)$, $\hat{F}_{Y|\bX}^{-1}(\beta|\cdot)$ and their derivatives up to order two to $h_1,F_{Y|\bX}^{-1}(\tau|\cdot),F_{Y|\bX}^{-1}(\beta|\cdot)$ and the corresponding derivatives, respectively. Without loss of generality, when proving $\hat{F}_{Y|\bX}^{-1}(\tau|\cdot)\in C_{R_{f_{m_{\tau}}}}^{\gamma_{f_{m_{\tau}}}}(M_{\bX})$ and $\hat{F}_{Y|\bX}^{-1}(\beta|\cdot)\in C_{R_{f_{m_{\beta}}}}^{\gamma_{f_{m_{\beta}}}}(M_{\bX})$ only derivatives with respect to $\bx_1$ are considered since other derivatives can be treated analogously. For $\hat{h}_1$ this can be derived from Lemma \ref{asympintegral}, since
		$$\frac{\partial}{\partial y}\hat{h}_1(y)=\frac{\partial}{\partial y}\exp\bigg(-\int_{y_1}^{y}\frac{1}{\hat{\lambda}(u)}\,du\bigg)=-\frac{\hat{h}_1(y)}{\hat{\lambda}(y)}$$
		and
		$$\frac{\partial^2}{\partial y^2}\hat{h}_1(y)=-\frac{\hat{\lambda}(y)\frac{\partial}{\partial y}\hat{h}_1(y)-\hat{h}_1(y)\frac{\partial}{\partial y}\hat{\lambda}(y)}{\hat{\lambda}(y)^2}=\frac{\hat{h}_1(y)+\hat{h}_1(y)\frac{\partial}{\partial y}\hat{\lambda}(y)}{\hat{\lambda}(y)^2}.$$
		As will be seen in the following, the assertion for $\hat{F}_{Y|\bX}^{-1}(\tau|\cdot)$ and $\hat{F}_{Y|\bX}^{-1}(\beta|\cdot)$ follows from the corresponding assertion for $\hat{F}_{Y|\bX}(y|\cdot)$ and hence can be obtained again as in \citet{Han2008}. One has
		$$\tau=\hat{F}_{Y|\bX}(\hat{F}_{Y|\bX}^{-1}(\tau|\bx)|\bx),$$
		so that
		$$0=\frac{\partial}{\partial \bx_1}\hat{F}_{Y|\bX}(\hat{F}_{Y|\bX}^{-1}(\tau|\bx)|\bx)=\hat{F}_y(\hat{F}_{Y|\bX}^{-1}(\tau|\bx)|\bx)\frac{\partial}{\partial \bx_1}\hat{F}_{Y|\bX}^{-1}(\tau|\bx)+\hat{F}_x(\hat{F}_{Y|\bX}^{-1}(\tau|\bx)|\bx),$$
		where $\hat{F}_y$ and $\hat{F}_x$ denote the derivative of $(y,\bx)\mapsto\hat{F}_{Y|\bX}(y|\bx)$ with respect to $y$ and $\bx$, respectively. Note that $\dX=1$ was assumed. Therefore,
		$$\frac{\partial}{\partial \bx_1}\hat{F}_{Y|\bX}^{-1}(\tau|\bx)=-\frac{\hat{F}_x(\hat{F}_{Y|\bX}^{-1}(\tau|\bx)|\bx)}{\hat{F}_y(\hat{F}_{Y|\bX}^{-1}(\tau|\bx)|\bx)}$$
		and
		\begin{align*}
		&\frac{\partial^2}{\partial \bx_1^2}\hat{F}_{Y|\bX}^{-1}(\tau|\bx)
		\\[0,2cm]&=-\frac{\hat{F}_y(\hat{F}_{Y|\bX}^{-1}(\tau|\bx)|\bx)\frac{\partial}{\partial \bx_1}\hat{F}_x(\hat{F}_{Y|\bX}^{-1}(\tau|\bx)|\bx)-\hat{F}_x(\hat{F}_{Y|\bX}^{-1}(\tau|\bx)|\bx)\frac{\partial}{\partial \bx_1}\hat{F}_y(\hat{F}_{Y|\bX}^{-1}(\tau|\bx)|\bx)}{\hat{F}_y(\hat{F}_{Y|\bX}^{-1}(\tau|\bx)|\bx)^2}.
		\end{align*}
		Similar to before, let $\hat{F}_{yy},\hat{F}_{xy}$ and $\hat{F}_{xx}$ denote the partial derivatives of $(y,\bx)\mapsto\hat{F}_{Y|\bX}(y|\bx)$ of order two. Then, it holds that
		$$\frac{\partial}{\partial \bx_1}\hat{F}_y(\hat{F}_{Y|\bX}^{-1}(\tau|\bx)|\bx)=\hat{F}_{xy}(\hat{F}_{Y|\bX}^{-1}(\tau|\bx)|\bx)+\hat{F}_{yy}(\hat{F}_{Y|\bX}^{-1}(\tau|\bx)|\bx)\frac{\partial}{\partial \bx_1}\hat{F}_{Y|\bX}^{-1}(\tau|\bx)$$
		as well as
		$$\frac{\partial}{\partial \bx_1}\hat{F}_x(\hat{F}_{Y|\bX}^{-1}(\tau|\bx)|\bx)=\hat{F}_{xx}(\hat{F}_{Y|\bX}^{-1}(\tau|\bx)|\bx)+\hat{F}_{xy}(\hat{F}_{Y|\bX}^{-1}(\tau|\bx)|\bx)\frac{\partial}{\partial \bx_1}\hat{F}_{Y|\bX}^{-1}(\tau|\bx).$$
		Consequently, the desired results
		$$\underset{\bx\in M_{\bX}}{\sup}\,\bigg|\frac{\partial}{\partial \bx_1}\hat{F}_{Y|\bX}^{-1}(\tau|\bx)-\frac{\partial}{\partial \bx_1}F_{Y|\bX}^{-1}(\tau|\bx)\bigg|=o_p(1)$$
		and
		$$\underset{\bx\in M_{\bX}}{\sup}\,\bigg|\frac{\partial^2}{\partial \bx_1^2}\hat{F}_{Y|\bX}^{-1}(\tau|\bx)-\frac{\partial^2}{\partial \bx_1^2}F_{Y|\bX}^{-1}(\tau|\bx)\bigg|=o_p(1)$$
		can be shown as in \citet{Han2008}.
		\hfill$\square$

		\begin{lemma}\label{lemmac1}
			With $A$ as in (\ref{defA}), it holds that $A(B,s_0)=0$ and $\hat{B}-B=o_p(1)$, that is, \ref{C1} is valid.
		\end{lemma}\noindent
		\textbf{Proof:} The first part follows directly from the definition of $A$ in (\ref{defA}) in Section \ref{estB}. For the second part, consider the function classes
		$$\mathcal{F}=\{(\bX,\varepsilon)\mapsto I_{\{\bX\in M_{\bX}\}}I_{\{\varepsilon\leq k_c(\mathfrak{h},f_{m_{\tau}},f_{m_{\beta}},\bX,e)\}}:s\in\mathcal{H},c\in[B_1,B_2],e\in[e_a,e_b]\}$$
		and
		\begin{align*}
		\tilde{\mathcal{F}}&=\{(\bX,\varepsilon)\mapsto I_{\{\bX\in M_{\bX}\}}I_{\{\bX\leq \bx\}}I_{\{\varepsilon\leq k_c(\mathfrak{h},f_{m_{\tau}},f_{m_{\beta}},\bX,e)\}}:
		\\[0,2cm]&\qquad s\in\mathcal{H},c\in[B_1,B_2],\bx\in M_{\bX},e\in[e_a,e_b]\}.
		\end{align*}
		It will be shown in the proof of Lemma \ref{lemmac5} below that the classes $\mathcal{F}$ and $\tilde{\mathcal{F}}$ are Donsker with respect to $\mathcal{L}^2\big(P^{(\bX,\varepsilon)}\big)$. Hence,
		\begin{align*}
		\hat{P}(\bX\leq \bx,\tilde{\varepsilon}_c(s_0)\leq e|\bX\in M_{\bX})&=\hat{P}(\bX\leq \bx,\varepsilon\leq k_c(s_0,\bX,e)|\bX\in M_{\bX})
		\\[0,2cm]&=\frac{\frac{1}{n}\sum_{i=1}^nI_{\{\bX_i\leq \bx,\varepsilon_i\leq k_c(s_0,\bX_i,e)\}}I_{\{\bX_i\in M_{\bX}\}}}{\frac{1}{n}\sum_{i=1}^nI_{\{\bX_i\in M_{\bX}\}}}
		\\[0,2cm]&=P(\bX\leq \bx,\varepsilon\leq k_c(s_0,\bX,e)|\bX\in M_{\bX})+\mathcal{O}_p\bigg(\frac{1}{\sqrt{n}}\bigg).
		\end{align*}
		Lemma \ref{lemmac4} yields
		$$\hat{f}_{m_{\tau}}(\bx)-F_{Y|\bX}^{-1}(\tau|\bx)=o_p\big(n^{-\frac{1}{4}}\big),\quad\hat{f}_{m_{\beta}}(\bx)-F_{Y|\bX}^{-1}(\beta|\bx)=o_p\big(n^{-\frac{1}{4}}\big)$$
		and
		$$\bar{h}_1(y)-h_1(y)=o_p\big(n^{-\frac{1}{4}}\big)$$
		uniformly in $y\in[z_a,z_b]$ and $\bx\in M_{\bX}$. Consequently, it holds that
		$$\underset{c\in[B_1,B_2],\bx\in M_{\bX},e\in[e_a,e_b]}{\sup}\,|k_c(\hat{s},\bx,e)-k_c(s_0,\bx,e)|=o_p(\delta_n),$$
		where the sequence $(\delta_n)_{n\in\mathbb{N}}$ follows from standard arguments, see Lemma 1.5.1 in \citet{Klo2019} for details. Assumption \ref{C4} ensures $\hat{s}\in\tilde{\mathcal{H}}$ with probability converging to one, so that Corollary 2.3.12 of \citet{vdVW1996} leads to
		\begin{align*}
		&\underset{c\in[B_1,B_2],\bx\in M,e\in[e_a,e_b]}{\sup}\,\big|\hat{P}(\bX\leq \bx,\tilde{\varepsilon}_c(\hat{s})\leq e|\bX\in M_{\bX})-P(\bX\leq \bx,\tilde{\varepsilon}_c(s_0)\leq e|\bX\in M_{\bX})\big|
		\\[0,2cm]&=\underset{c\in[B_1,B_2],\bx\in M,e\in[e_a,e_b]}{\sup}\,\big|\hat{P}(\bX\leq \bx,\varepsilon\leq k_c(\hat{s},\bX,e)|\bX\in M_{\bX})
		\\[0,2cm]&\quad\quad-P(\bX\leq \bx,\varepsilon\leq k_c(s_0,\bX,e)|\bX\in M_{\bX})\big|
		\\[0,2cm]&=o_p(1).
		\end{align*}
		Analogous calculations can be done for $\hat{P}(\bX\leq \bx|\bX\in M_{\bX})$ and $\hat{P}(\tilde{\varepsilon}_c(s_0)\leq e|\bX\in M_{\bX})$. Therefore,
		\begin{align*}
		\hat{A}(c,\hat{s})&=\bigg(\int_M\int_{[e_a,e_b]}\big(\hat{P}(\bX\leq \bx,\tilde{\varepsilon}_c(\hat{s})\leq e|\bX\in M_{\bX})
		\\[0,2cm]&\quad\quad-\hat{P}(\bX\leq \bx|\bX\in M_{\bX})\hat{P}(\tilde{\varepsilon}_c(\hat{s})\leq e|\bX\in M_{\bX})\big)^2\,de\,d\bx\bigg)^{\frac{1}{2}}
		\\[0,2cm]&=\bigg(\int_M\int_{[e_a,e_b]}\big(P(\bX\leq \bx,\tilde{\varepsilon}_c(s_0)\leq e|\bX\in M_{\bX})
		\\[0,2cm]&\quad\quad-P(\bX\leq \bx|\bX\in M_{\bX})P(\tilde{\varepsilon}_c(s_0)\leq e|\bX\in M_{\bX})\big)^2\,de\,d\bx\bigg)^{\frac{1}{2}}+o_p(1)
		\\[0,2cm]&=A(c,s_0)+o_p(1)
		\end{align*}
		uniformly in $c\in[B_1,B_2]$. Since the map $c\mapsto A(c,s_0)$ is continuous and $c=B$ is the unique minimizer, it holds that
		$$\underset{c\in[B_1,B_2],|c-B|>\delta}{\inf}\,A(c,s_0)>0$$
		for all $\delta>0$ and thus, $\hat{B}=\underset{c\in[B_1,B_2]}{\arg\min}\,\hat{A}(c,\hat{s})=\underset{c\in[B_1,B_2]}{\arg\min}\,A(c,s_0)+o_p(1)=B+o_p(1)$.
		\hfill$\square$

		\begin{lemma}\label{lemmac2}
			The ordinary derivative $\Gamma_1(c,s_0)(\bx,e)$ of $G_{MD}(c,s_0)(\bx,e)$ (with respect to c) exists for all $(\bx,e)\in M_{\bX}\times[e_a,e_b]$ in a neighbourhood of $B$ and is continuous at $c=B$. $\Gamma_1(B,s_0)(\bx,e)$ is different from zero on a set with positive $\lambda_{M_{\bX}\times[e_a,e_b]}$-measure. Consequently, \ref{C2} holds true.
		\end{lemma}\noindent
		\textbf{Proof:} The proof can be divided into three steps namely the proof of the continuous differentiability of $c\mapsto k_c(s_0,\bx,e)$, the proof of continuous differentiability of $c\mapsto P(\bX\leq \bx,\varepsilon\leq k_c(s_0,\bX,e)|\bX\in M_{\bX})$ and $c\mapsto G_{MD}(c,s_0)$ (each for all $(\bx,e)\in M\times[e_a,e_b]$) and finally the proof of $\Gamma_1(B,s_0)(\bx,e)\neq0$, which is omitted here and can be found in \citet{Klo2019}. First, recall the definition of $k_c$ in (\ref{defkc}) and write for all $(\bx,e)\in M\times[e_a,e_b]$
		\begin{align}
		&\frac{\partial}{\partial c}k_c(s_0,\bx,e)\nonumber
		\\[0,2cm]&=\frac{\partial}{\partial c}\frac{\big(h_c(F_{Y|\bX}^{-1}(\tau|\bx))+e(h_c(F_{Y|\bX}^{-1}(\beta|\bx))-h_c(F_{Y|\bX}^{-1}(\tau|\bx)))\big)^{\frac{B}{c}}-g(\bx)}{\sigma(\bx)}\nonumber
		\\[0,2cm]&=\frac{1}{\sigma(\bx)}\Bigg[-\frac{B}{c^2}\bigg(h_1\big(F_{Y|\bX}^{-1}(\tau|\bx)\big)^c+e\Big(h_1\big(F_{Y|\bX}^{-1}(\beta|\bx)\big)^c-h_1\big(F_{Y|\bX}^{-1}(\tau|\bx)\big)^c\Big)\bigg)^{\frac{B}{c}}\nonumber
		\\[0,2cm]&\quad\quad\log\bigg(h_1\big(F_{Y|\bX}^{-1}(\tau|\bx)\big)^c+e\Big(h_1\big(F_{Y|\bX}^{-1}(\beta|\bx)\big)^c-h_1\big(F_{Y|\bX}^{-1}(\tau|\bx)\big)^c\Big)\bigg)\nonumber
		\\[0,2cm]&\quad+\frac{B}{c}\bigg(h_1\big(F_{Y|\bX}^{-1}(\tau|\bx)\big)^c+e\Big(h_1\big(F_{Y|\bX}^{-1}(\beta|\bx)\big)^c-h_1\big(F_{Y|\bX}^{-1}(\tau|\bx)\big)^c\Big)\bigg)^{\frac{B}{c}-1}\nonumber
		\\[0,2cm]&\quad\quad\bigg(\log\Big(h_1\big(F_{Y|\bX}^{-1}(\tau|\bx)\big)\Big)h_1\big(F_{Y|\bX}^{-1}(\tau|\bx)\big)^c+e\bigg(\log\Big(h_1\big(F_{Y|\bX}^{-1}(\beta|\bx)\big)\Big)h_1\big(F_{Y|\bX}^{-1}(\beta|\bx)\big)^c\nonumber
		\\[0,2cm]&\quad\quad-\log\Big(h_1\big(F_{Y|\bX}^{-1}(\tau|\bx)\big)\Big)h_1\big(F_{Y|\bX}^{-1}(\tau|\bx)\big)^c\bigg)\bigg)\Bigg].\label{kcderivative}
		\end{align}
		Due to $0<h_c(F_{Y|\bX}^{-1}(\tau|\bx)),h_c(F_{Y|\bX}^{-1}(\beta|\bx))$ as well as
		$$0<h_c(z_a)\leq h_c(F_{Y|\bX}^{-1}(\tau|\bx))+e(h_c(F_{Y|\bX}^{-1}(\beta|\bx))-h_c(F_{Y|\bX}^{-1}(\tau|\bx)))\leq h_c(z_b)$$
		for all $\bx\in M_{\bX},e\in[e_a,e_b]$ the function $(c,\bx,e)\mapsto\frac{\partial}{\partial c}k_c(s_0,\bx,e)$ is well defined, continuous and thus bounded on $[B_1,B_2]\times M_{\bX}\times[e_a,e_b]$. 
		
		Second, $P(\bX\leq \bx,\varepsilon\leq k_c(s_0,\bX,e)|\bX\in M_{\bX})$ can be written as
		\begin{align*}
		&P(\bX\leq \bx,\varepsilon\leq k_c(s_0,\bX,e)|\bX\in M_{\bX})
		\\[0,2cm]&=\frac{P(\bX\leq \bx,\varepsilon\leq k_c(s_0,\bX,e),\bX\in M_{\bX})}{P(\bX\in M_{\bX})}
		\\[0,2cm]&=\frac{1}{P(\bX\in M_{\bX})}\int_{M_{\bX}\cap(-\infty,\bx]}F_{\varepsilon}(k_c(s_0,\mathbf{v},e))f_{\bX}(\mathbf{v})\,d\mathbf{v}.
		\end{align*}
		Analogously,
		$$P(\varepsilon\leq k_c(s_0,\bX,e)|\bX\in M_{\bX})=\frac{1}{P(\bX\in M_{\bX})}\int_{M_{\bX}}F_{\varepsilon}(k_c(s_0,\mathbf{v},e))f_{\bX}(\mathbf{v})\,d\mathbf{v}.$$
		The Dominated Convergence Theorem leads to
		\begin{align*}
		&\frac{\partial}{\partial c}P(\bX\leq \bx,\varepsilon\leq k_c(s_0,\bX,e)|\bX\in M_{\bX})
		\\[0,2cm]&\quad=\frac{1}{P(\bX\in M_{\bX})}\int_{M_{\bX}\cap(-\infty,\bx]}\frac{\partial}{\partial c}F_{\varepsilon}(k_c(s_0,\mathbf{v},e))f_{\bX}(\mathbf{v})\,d\mathbf{v}
		\\[0,2cm]&\quad=\frac{1}{P(\bX\in M_{\bX})}\int_{M_{\bX}\cap(-\infty,\bx]}f_{\varepsilon}(k_c(s_0,\mathbf{v},e))\frac{\partial}{\partial c}k_c(s_0,\mathbf{v},e)f_{\bX}(\mathbf{v})\,d\mathbf{v},
		\end{align*}
		where the supremum of the integrand, which is continuous and evaluated on a compact set, can be taken as a majorant. Consequently
		\begin{align*}
		\frac{\partial}{\partial c}G_{MD}(c,s_0)(\bx,e)&=\frac{1}{P(\bX\in M_{\bX})}\int_{M_{\bX}}f_{\varepsilon}(k_c(s_0,\mathbf{v},e))
		\\[0,2cm]&\quad\quad\frac{\partial}{\partial c}k_c(s_0,\mathbf{v},e)(I_{(-\infty,\bx]}(\mathbf{v})-P(\bX\leq \bx|\bX\in M_{\bX}))f_{\bX}(\mathbf{v})\,d\mathbf{v}.
		\end{align*}
		\hfill$\square$

		\begin{lemma}\label{lemmac3}
			There exists a $\delta>0$ such that for all $c\in B_{\delta},(\bx,e)\in M_{\bX}\times[e_a,e_b]$ the directional derivative $\Gamma_2(c,s_0)(\bx,e)[s-s_0]$ of $G_{MD}(c,s_0)(\bx,e)$ with respect to $s$ exists in all directions $[s-s_0]$. Moreover, consider a positive sequence $\delta_n\rightarrow0$ and $(c,s)\in B_{\delta_n}\times\tilde{\mathcal{H}}_{\delta_n}$. Then,
			\begin{enumerate}[label=(\roman*)]
				\item for an appropriate constant $C\geq0$ one has
				\begin{align*}
				&||G_{MD}(c,s)-G_{MD}(c,s_0)-\Gamma_2(c,s_0)[s-s_0]||_2
				\\[0,2cm]&\leq C\big(||\mathfrak{h}-h_1||_{[z_a,z_b]}^{\frac{3}{2}}+||f_{m_{\tau}}-F_{Y|\bX}^{-1}(\tau|\cdot)||_{M_{\bX}}^2+||f_{m_{\beta}}-F_{Y|\bX}^{-1}(\beta|\cdot)||_{M_{\bX}}^2\big).
				\end{align*}
				\item one has $||\Gamma_2(c,s_0)[\hat{s}-s_0]-\Gamma_2(B,s_0)[\hat{s}-s_0]||=o_p(|c-B|)+o_p\big(n^{-\frac{1}{2}}\big)$.
			\end{enumerate}
			Therefore, \ref{C3} is valid.
		\end{lemma}\noindent
		\textbf{Proof:} First, existence of the directional derivatives is shown, before conditions $(i)$ and $(ii)$ are proven.\\[0,2cm]
		\textbf{Directional derivative with respect to} $\mathbf{h}$: Define for some fixed $c,\mathfrak{h},\bx,e$\label{deffhtzct}
		\begin{align}
		f_{h,t}&:=h_1+t(\mathfrak{h}-h_1),\label{deffht}
		\\[0,2cm]\psi(t,z)&:=f_{h,t}^{-1}(z),\nonumber
		\\[0,2cm]z_c(f_{h,t},f_{m_{\tau}},f_{m_{\beta}},\bx,e)&:=\big(f_{h,t}^c(f_{m_{\tau}}(\bx))+e(f_{h,t}^c(f_{m_{\beta}}(\bx))-f_{h,t}^c(f_{m_{\tau}}(\bx)))\big)^{\frac{1}{c}}.\label{defzc}
		\end{align}
		Mostly, the components $\bx,e$ will be omitted and $z_c(t)$ will be written as an abbreviation for $z_c(f_{h,t},f_{m_{\tau}},f_{m_{\beta}},\bx,e)$. Further, all derivatives with respect to $t$ are marked with a ``$\,\cdot\,$", those with respect to $y$ are marked with a ``$\,'\,$". Then, one has
		$$h_1'(z)=-\frac{h_1(z)}{\lambda(z)},$$
		$$\psi'(t,z)=\frac{1}{f_{h,t}'(f_{h,t}^{-1}(z_c(t)))}\overset{t=0}{=}-\frac{\lambda(h_1^{-1}(z_c(0)))}{z_c(0)}$$
		as well as
		\begin{align*}
		&\frac{\partial}{\partial t}\psi(t,f_{h,t}(h_1^{-1}(z_c(t))))
		\\[0,2cm]&=\dot{\psi}(t,f_{h,t}(h_1^{-1}(z_c(t))))
		\\[0,2cm]&\quad+\psi'(t,f_{h,t}(h_1^{-1}(z_c(t))))\bigg(\dot{f}_{h,t}(h_1^{-1}(z_c(t)))+f_{h,t}'(h_1^{-1}(z_c(t)))\frac{\frac{\partial}{\partial t}z_c(t)}{h_1'(h_1^{-1}(z_c(t)))}\bigg)
		\\[0,2cm]&\overset{t\rightarrow0}{\longrightarrow}\dot{\psi}(0,z_c(0))+\psi'(0,z_c(0))\bigg((\mathfrak{h}-h_1)(h_1^{-1}(z_c(0)))+\frac{\partial}{\partial t}z_c(t)\big|_{t=0}\bigg)
		\\[0,2cm]&=\dot{\psi}(0,z_c(0))-\frac{\lambda(h_1^{-1}(z_c(0)))}{z_c(0)}\bigg(\mathfrak{h}(h_1^{-1}(z_c(0)))-z_c(0)+\frac{\partial}{\partial t}z_c(t)\bigg|_{t=0}\bigg).
		\end{align*}
		Due to
		$$\frac{\partial}{\partial t}\psi(t,f_{h,t}(h_1^{-1}(z_c(t))))=\frac{\frac{\partial}{\partial t}z_c(t)}{h_1'(h_1^{-1}(z_c(t)))}\overset{t\rightarrow0}{\longrightarrow}-\frac{\lambda(h_1^{-1}(z_c(0)))\frac{\partial}{\partial t}z_c(t)\big|_{t=0}}{z_c(0)},$$
		it holds that
		$$\dot{\psi}(0,z_c(0))=\frac{\lambda(h_1^{-1}(z_c(0)))}{z_c(0)}(\mathfrak{h}(h_1^{-1}(z_c(0)))-z_c(0)),$$
		so that
		\begin{align*}
		\frac{\partial}{\partial t}\psi(t,z_c(t))&=\dot{\psi}(t,z_c(t))+\psi'(t,z_c(t))\frac{\partial}{\partial t}z_c(t)
		\\[0,2cm]&\overset{t\rightarrow0}{\longrightarrow}\dot{\psi}(0,z_c(0))+\psi'(0,z_c(0))\frac{\partial}{\partial t}z_c(t)\bigg|_{t=0}
		\\[0,2cm]&=\frac{\lambda(h_1^{-1}(z_c(0)))}{z_c(0)}\bigg(\mathfrak{h}(h_1^{-1}(z_c(0)))-z_c(0)-\frac{\partial}{\partial t}z_c(t)\bigg|_{t=0}\bigg).
		\end{align*}
		Additionally,
		\begin{align*}
		\frac{\partial}{\partial t}z_c(t)\bigg|_{t=0}&=\frac{\partial}{\partial t}\big(f_{h,t}^c(f_{m_{\tau}}(\bx))+e(f_{h,t}^c(f_{m_{\beta}}(\bx))-f_{h,t}^c(f_{m_{\tau}}(\bx)))\big)^{\frac{1}{c}}\bigg|_{t=0}
		\\[0,2cm]&=\frac{1}{c}\big(h_1^c(f_{m_{\tau}}(\bx))+e(h_1^c(f_{m_{\beta}}(\bx))-h_1^c(f_{m_{\tau}}(\bx)))\big)^{\frac{1}{c}-1}
		\\[0,2cm]&\quad\big(ch_1^{c-1}(f_{m_{\tau}}(\bx))(\mathfrak{h}(f_{m_{\tau}}(\bx))-h_1(f_{m_{\tau}}(\bx))
		)
		\\[0,2cm]&\quad+e(ch_1^{c-1}(f_{m_{\beta}}(\bx))(\mathfrak{h}(f_{m_{\beta}}(\bx))-h_1(f_{m_{\beta}}(\bx))
		\\[0,2cm]&\quad\quad-ch_1^{c-1}(f_{m_{\tau}}(\bx))(\mathfrak{h}(f_{m_{\tau}}(\bx))-h_1(f_{m_{\tau}}(\bx)))\big).
		\end{align*}
		This in turn results in (for the special case $f_{m_{\tau}}=F_{Y|\bX}^{-1}(\tau|\cdot),f_{m_{\beta}}=F_{Y|\bX}^{-1}(\beta|\cdot)$)
		\begin{align}
		D_hk_c(s_0,\bx,e)[\mathfrak{h}-h_1]&=\frac{\partial}{\partial t}k_c(f_{h,t},F_{Y|\bX}^{-1}(\tau|\cdot),F_{Y|\bX}^{-1}(\beta|\cdot),\bx,e)\bigg|_{t=0}\nonumber
		\\[0,2cm]&=\frac{\partial}{\partial t}\frac{h_1(\psi(t,z_c(t)))^B-g(\bx)}{\sigma(\bx)}\bigg|_{t=0}\nonumber
		\\[0,2cm]&=\frac{Bh_1(\psi(t,z_c(t)))^{B-1}h_1'(\psi(t,z_c(t)))\frac{\partial}{\partial t}\psi(t,z_c(t))}{\sigma(\bx)}\bigg|_{t=0}\nonumber
		\\[0,2cm]&=\frac{Bz_c(0)^{B-1}\big(z_c(0)-\mathfrak{h}(h_1^{-1}(z_c(0)))+\frac{\partial}{\partial t}z_c(t)\big|_{t=0}\big)}{\sigma(\bx)}\label{Dhkc}
		\end{align}
		and by applying the Dominated Convergence Theorem to equation (\ref{exprGMDkc})
		\begin{align*}
		D_hG_{MD}(c,s_0)(\bx,e)[\mathfrak{h}-h_1]&=\frac{1}{P(\bX\in M_{\bX})}\bigg(\int_{M_{\bX}}\big(I_{\{\mathbf{w}\leq \bx\}}-P(\bX\leq \bx|\bX\in M_{\bX})\big)
		\\[0,2cm]&\quad f_{\varepsilon}(k_c(s_0,\mathbf{w},e))D_hk_c(s_0,\mathbf{w},e)[\mathfrak{h}-h_1]f_{\bX}(\mathbf{w})\,d\mathbf{w}\bigg).
		\end{align*}
		\textbf{Directional derivative with respect to} $\mathbf{f_{m_{\tau}}}$\textbf{ and }$\mathbf{f_{m_{\beta}}}$: For $\mathfrak{h}=h_1$, $k_c$ simplifies to
		$$k_c(h_1,f_{m_{\tau}},f_{m_{\beta}},\bx,e)=\frac{z_c(h_1,f_{m_{\tau}},f_{m_{\beta}})^B-g(\bx)}{\sigma(\bx)}.$$
		Hence, with
		$$f_{m_{\tau},t}=F_{Y|\bX}^{-1}(\tau|\cdot)+t\big(f_{m_{\tau}}-F_{Y|\bX}^{-1}(\tau|\cdot)\big)\quad\textup{and}\quad f_{m_{\beta},t}=F_{Y|\bX}^{-1}(\beta|\cdot)+t\big(f_{m_{\beta}}-F_{Y|\bX}^{-1}(\beta|\cdot)\big)$$
		one has
		\begin{align*}
		&D_{f_{m_{\tau}}}k_c(h_1,F_{Y|\bX}^{-1}(\tau|\cdot),F_{Y|\bX}^{-1}(\beta|\cdot),\bx,e)\big[f_{m_{\tau}}-F_{Y|\bX}^{-1}(\tau|\cdot)\big]
		\\[0,2cm]&=\frac{\partial}{\partial t}\frac{z_c(h_1,f_{m_{\tau},t},F_{Y|\bX}^{-1}(\beta|\cdot))^B-g(\bx)}{\sigma(\bx)}\bigg|_{t=0}
		\\[0,2cm]&=\frac{\partial}{\partial t}\frac{\big(h_c(f_{m_{\tau},t}(\bx))+e(h_c(F_{Y|\bX}^{-1}(\beta|\bx))-h_c(f_{m_{\tau},t}(\bx)))\big)^{\frac{B}{c}}-g(\bx)}{\sigma(\bx)}\bigg|_{t=0}
		\\[0,2cm]&=-\frac{B(...)^{\frac{B}{c}-1}(1-e)h_{c}(F_{Y|\bX}^{-1}(\tau|\bx))(f_{m_{\tau}}(\bx)-F_{Y|\bX}^{-1}(\tau|\bx))}{\sigma(\bx)\lambda(F_{Y|\bX}^{-1}(\tau|\bx))}
		\end{align*}
		as well as
		\begin{align*}
		&D_{f_{\beta}}k_c(s_0,\bx,e)\big[f_{\beta}-F_{Y|\bX}^{-1}(\beta|\cdot)\big]
		\\[0,2cm]&=\frac{\partial}{\partial t}\frac{\big(h_c(F_{Y|\bX}^{-1}(\tau|\bx))+e(h_c(f_{m_{\beta},t}(\bx))-h_c(F_{Y|\bX}^{-1}(\tau|\bx)))\big)^{\frac{B}{c}}-g(\bx)}{\sigma(\bx)}\bigg|_{t=0}
		\\[0,2cm]&=-\frac{B(...)^{\frac{B}{c}-1}eh_{c}(F_{Y|\bX}^{-1}(\beta|\bx))(f_{m_{\beta}}(\bx)-F_{Y|\bX}^{-1}(\beta|\bx))}{\sigma(\bx)\lambda(F_{Y|\bX}^{-1}(\beta|\bx))}.
		\end{align*}
		The Dominated Convergence Theorem yields
		\begin{align*}
		&D_{f_{m_{\tau}}}G_{MD}(c,s_0)(\bx,e)\big[f_{m_{\tau}}-F_{Y|\bX}^{-1}(\tau|\cdot)\big]
		\\[0,2cm]&=\frac{1}{P(\bX\in M_{\bX})}\Big(\int_{M_{\bX}}f_{\varepsilon}(k_c(s_0,\mathbf{w},e))\big(I_{\{\mathbf{w}\leq \bx\}}-P(\bX\leq \bx|\bX\in M_{\bX})\big)
		\\[0,2cm]&\quad D_{f_{m_{\tau}}}k_c(s_0,\mathbf{w},e)\big[f_{m_{\tau}}-F_{Y|\bX}^{-1}(\tau|\cdot)\big]f_{\bX}(\mathbf{w})\,d\mathbf{w}\Big).
		\end{align*}
		\textbf{Directional derivative with respect to} $\mathbf{s}$: This results from the previous parts of the proof as follows: Define $s_t=(f_{h,t},f_{m_{\tau},t},f_{m_{\beta},t})$ and $t\mapsto\tilde{z}_c(t):=z_c(s_t,\bx,e)$. Then,
		$$k_c(s_t,\bx,e)=\frac{\partial}{\partial t}k_c(s_t,\bx,e)\bigg|_{t=0}=\frac{\partial}{\partial t}\frac{h_1(\psi(t,\tilde{z}_c(t)))^B-g(\bx)}{\sigma(\bx)}\bigg|_{t=0}$$
		only depends on $f_{m_{\tau},t}$ and $f_{m_{\beta},t}$ via $\tilde{z}_c(t)$, respectively.
		Due to $\tilde{z}_c(0)=z_c(s_0)$ ($=z_c(0)$ with the notation from before), one can proceed as for the derivative with respect to $h$ to obtain
		$$D_sk_c(s_0,\bx,e)[s-s_0]=\frac{B\tilde{z}_c(0)^{B-1}(\tilde{z}_c(0)-\mathfrak{h}(h_1^{-1}(\tilde{z}_c(0)))+\frac{\partial}{\partial t}\tilde{z}_c(t)\big|_{t=0})}{\sigma(\bx)},$$
		where $D_s$ denotes the derivative with respect to $s$. At the same time,
		\begin{align*}
		\frac{\partial}{\partial t}\tilde{z}_c(t)\big|_{t=0}&=\frac{\partial}{\partial t}z_c(f_{h,t},f_{m_{\tau},t},f_{m_{\beta},t})\bigg|_{t=0}
		\\[0,2cm]&=\left(\begin{array}{ccccc}D_hz_c(s_0)&,&D_{f_{m_{\tau}}}z_c(s_0)&,&D_{f_{m_{\beta}}}z_c(s_0)\end{array}\right)\left(\begin{array}{c}\mathfrak{h}-h_1\\f_{m_{\tau}}-F_{Y|\bX}^{-1}(\tau|\bx)\\f_{m_{\beta}}-F_{Y|\bX}^{-1}(\beta|\bx)\end{array}\right),
		\end{align*}
		which in total leads to
		\begin{align*}
		D_sk_c(s_0,\bx,e)[s-s_0]&=D_hk_c(s_0,\bx,e)[\mathfrak{h}-h_1]+D_{f_{m_{\tau}}}k_c(s_0,\bx,e)\big[f_{m_{\tau}}-F_{Y|\bX}^{-1}(\tau|\cdot)\big]
		\\[0,2cm]&\quad+D_{f_{m_{\beta}}}k_c(s_0,\bx,e)\big[f_{m_{\beta}}-F_{Y|\bX}^{-1}(\beta|\cdot)\big]
		\end{align*}
		and (after applying the Dominated Convergence Theorem)
		\begin{align}
		&D_sG_{MD}(c,s_0)(\bx,e)[s-s_0]\nonumber
		\\[0,2cm]&=D_hG_{MD}(c,s_0)(\bx,e)[\mathfrak{h}-h_1]+D_{f_{m_{\tau}}}G_{MD}(c,s_0)(\bx,e)\big[f_{m_{\tau}}-F_{Y|\bX}^{-1}(\tau|\cdot)\big]\nonumber
		\\[0,2cm]&\quad+D_{f_{m_{\beta}}}G_{MD}(c,s_0)(\bx,e)\big[f_{m_{\beta}}-F_{Y|\bX}^{-1}(\beta|\cdot)\big].\label{sepGamma2}
		\end{align}
		\textbf{Proof of }$\mathbf{(i)}$:
		First, the following Lemma can be shown by straightforward calculations, see \citet{Klo2019} for details.
		\begin{lemma}\label{lemmahdevinv}
			Let $\delta_n\searrow0,s=(\mathfrak{h},f_{m_{\tau}},f_{m_{\beta}})\in\tilde{\mathcal{H}}_{\delta_n}$ and $0<\eta<\frac{h_1(z_b)-h_1(z_a)}{2}$. Then,
			$$\underset{t\in[h_1(z_a)+\eta,h_1(z_b)-\eta]}{\sup}\,|\mathfrak{h}^{-1}(t)-h_1^{-1}(t)|=\mathcal{O}\big(||\mathfrak{h}-h_1||_{[z_a,z_b]}\big)$$
			and
			$$||\mathfrak{h}'-h_1'||_{[z_a,z_b]}=\mathcal{O}\Big(\sqrt{||\mathfrak{h}-h_1||_{[z_a,z_b]}}\Big)=\sqrt{\delta_n}.$$
		\end{lemma}
		
		Let $\delta_n\searrow0$ and $(c,s)\in B_{\delta_n}\times\tilde{\mathcal{H}}_{\delta_n}$. To apply the lemma from above, use (\ref{sepGamma2}) to split the norm into three parts
		\begin{align*}
		&||G_{MD}(c,s)-G_{MD}(c,s_0)-\Gamma_2(c,s_0)[s-s_0]||_2
		\\[0,2cm]&=||G_{MD}(c,\mathfrak{h},f_{m_{\tau}},f_{m_{\beta}})-G_{MD}(c,h_1,f_{m_{\tau}},f_{m_{\beta}})-D_hG_{MD}(c,s_0)[\mathfrak{h}-h_1]
		\\[0,2cm]&\quad+G_{MD}(c,h_1,f_{m_{\tau}},f_{m_{\beta}})-G_{MD}(c,h_1,F_{Y|\bX}^{-1}(\tau|\cdot),f_{m_{\beta}})
		\\[0,2cm]&\quad-D_{f_{m_{\tau}}}G_{MD}(c,s_0)\big[f_{m_{\tau}}-F_{Y|\bX}^{-1}(\tau|\cdot)\big]+G_{MD}(c,h_1,F_{Y|\bX}^{-1}(\tau|\cdot),f_{m_{\beta}})
		\\[0,2cm]&\quad-G_{MD}(c,h_1,F_{Y|\bX}^{-1}(\tau|\cdot),F_{Y|\bX}^{-1}(\beta|\cdot))-D_{f_{m_{\beta}}}G_{MD}(c,s_0)\big[f_{m_{\beta}}-F_{Y|\bX}^{-1}(\beta|\cdot)\big]||_2
		\\[0,2cm]&\leq||G_{MD}(c,\mathfrak{h},f_{m_{\tau}},f_{m_{\beta}})-G_{MD}(c,h_1,f_{m_{\tau}},f_{m_{\beta}})-D_hG_{MD}(c,s_0)[\mathfrak{h}-h_1]||_2
		\\[0,2cm]&\quad+||G_{MD}(c,h_1,f_{m_{\tau}},f_{m_{\beta}})-G_{MD}(c,h_1,F_{Y|\bX}^{-1}(\tau|\cdot),f_{m_{\beta}})
		\\[0,2cm]&\quad-D_{f_{m_{\tau}}}G_{MD}(c,s_0)\big[f_{m_{\tau}}-F_{Y|\bX}^{-1}(\tau|\cdot)\big]||_2+||G_{MD}(c,h_1,F_{Y|\bX}^{-1}(\tau|\cdot),f_{m_{\beta}})-G_{MD}(c,s_0)
		\\[0,2cm]&\quad-D_{f_{m_{\beta}}}G_{MD}(c,s_0)\big[f_{m_{\beta}}-F_{Y|\bX}^{-1}(\beta|\cdot)\big]||_2.
		\end{align*}
		Notice for the first summand that due to
		\begin{align*}
		&||G_{MD}(c,\mathfrak{h},f_{m_{\tau}},f_{m_{\beta}})-G_{MD}(c,h_1,f_{m_{\tau}},f_{m_{\beta}})-D_hG_{MD}(c,s_0)[\mathfrak{h}-h_1]||_2
		\\[0,2cm]&=\bigg(\int_{M_{\bX}}\int_{[e_a,e_b]}\Big(P(\bX\leq \bx,\varepsilon\leq k_c(\mathfrak{h},f_{m_{\tau}},f_{m_{\beta}},\bX,e)|\bX\in M_{\bX})
		\\[0,2cm]&\quad-P(\bX\leq \bx|\bX\in M_{\bX})P(\varepsilon\leq k_c(\mathfrak{h},f_{m_{\tau}},f_{m_{\beta}},\bX,e)|\bX\in M_{\bX})
		\\[0,2cm]&\quad-P(\bX\leq \bx,\varepsilon\leq k_c(h_1,f_{m_{\tau}},f_{m_{\beta}},\bX,e)|\bX\in M_{\bX})
		\\[0,2cm]&\quad+P(\bX\leq \bx|\bX\in M_{\bX})P(\varepsilon\leq k_c(h_1,f_{m_{\tau}},f_{m_{\beta}},\bX,e)|\bX\in M_{\bX})
		\\[0,2cm]&\quad-D_hG_{MD}(h_1,F_{Y|\bX}^{-1}(\tau|\cdot),F_{Y|\bX}^{-1}(\beta|\cdot),\bX,e)(\bx,e)[\mathfrak{h}-h_1]\Big)^2\,de\,d\bx\bigg)^{\frac{1}{2}}
		\\[0,2cm]&=\bigg(\int_{M_{\bX}}\int_{[e_a,e_b]}\bigg(\frac{1}{P(\bX\in M_{\bX})}\int_{M_{\bX}}\big(I_{\{v\leq \bx\}}-P(\bX\leq \bx|\bX\in M_{\bX})\big)
		\\[0,2cm]&\quad\quad\Big(F_{\varepsilon}(k_c(\mathfrak{h},f_{m_{\tau}},f_{m_{\beta}},v,e))-F_{\varepsilon}(k_c(h_1,f_{m_{\tau}},f_{m_{\beta}},v,e))-f_{\varepsilon}(k_c(h_1,f_{m_{\tau}},f_{m_{\beta}},v,e))
		\\[0,2cm]&\quad\quad D_hk_c(h_1,F_{Y|\bX}^{-1}(\tau|\cdot),F_{Y|\bX}^{-1}(\beta|\cdot),v,e)[\mathfrak{h}-h_1]\Big)f_{\bX}(v)\,dv\bigg)^2\,de\,d\bx\bigg)^{\frac{1}{2}}
		\end{align*}
		and
		\begin{align*}
		&F_{\varepsilon}(k_c(\mathfrak{h},f_{m_{\tau}},f_{m_{\beta}},v,e))-F_{\varepsilon}(k_c(h_1,f_{m_{\tau}},f_{m_{\beta}},v,e))
		\\[0,2cm]&-f_{\varepsilon}(k_c(h_1,f_{m_{\tau}},f_{m_{\beta}},v,e))D_hk_c(h_1,F_{Y|\bX}^{-1}(\tau|\cdot),F_{Y|\bX}^{-1}(\beta|\cdot),v,e)[\mathfrak{h}-h_1]
		\\[0,2cm]&=\quad f_{\varepsilon}(k_c(h_1,f_{m_{\tau}},f_{m_{\beta}},v,e))(k_c(\mathfrak{h},f_{m_{\tau}},f_{m_{\beta}},v,e)-k_c(h_1,f_{m_{\tau}},f_{m_{\beta}},v,e))
		\\[0,2cm]&\quad\quad+f_{\varepsilon}'(\tilde{k})(k_c(h_1,f_{m_{\tau}},f_{m_{\beta}},v,e)-k_c(\mathfrak{h},f_{m_{\tau}},f_{m_{\beta}},v,e))^2
		\\[0,2cm]&\quad\quad-f_{\varepsilon}(k_c(h_1,f_{m_{\tau}},f_{m_{\beta}},v,e))D_hk_c(h_1,F_{Y|\bX}^{-1}(\tau|\cdot),F_{Y|\bX}^{-1}(\beta|\cdot),v,e)[\mathfrak{h}-h_1]
		\end{align*}
		for some $\tilde{k}$ between $k_c(\mathfrak{h},f_{m_{\tau}},f_{m_{\beta}},v,e)$ and $k_c(h_1,f_{m_{\tau}},f_{m_{\beta}},v,e)$ it suffices to prove
		$$\big|k_c(\mathfrak{h},f_{m_{\tau}},f_{m_{\beta}},v,e)-k_c(h_1,f_{m_{\tau}},f_{m_{\beta}},v,e)-D_hk_c(s_0,v,e)[\mathfrak{h}-h_1]\big|\leq C||\mathfrak{h}-h_1||_{[z_a,z_b]}^{\frac{3}{2}}$$
		for an appropriate $C>0$ and uniformly in $c\in B_{\delta},v\in M_{\bX},e\in[e_a,e_b],f_{m_{\tau}},f_{m_{\beta}}$, such that $(\mathfrak{h},f_{m_{\tau}},f_{m_{\beta}})\in\tilde{\mathcal{H}}$. Analogous calculations for $D_{f_{m_{\tau}}}$ and $D_{f_{m_{\beta}}}$ yield the sufficient conditions
		\begin{align}
		&\big|k_c(h_1,f_{m_{\tau}},f_{m_{\beta}},v,e)-k_c(h_1,F_{Y|\bX}^{-1}(\tau|\cdot),f_{m_{\beta}},v,e)\nonumber
		\\[0,2cm]\quad&-D_{f_{m_{\tau}}}k_c(h_1,F_{Y|\bX}^{-1}(\tau|\cdot),f_{m_{\beta}},v,e)\big[f_{m_{\tau}}-F_{Y|\bX}^{-1}(\tau|\cdot)\big]\big|\nonumber
		\\[0,2cm]&\leq C||f_{m_{\tau}}-F_{Y|\bX}^{-1}(\tau|\cdot)||_{M_{\bX}}^2\label{proofasympBlemmac3sufficientalpha}
		\end{align}
		(uniformly in $c\in B_{\delta},v\in M_{\bX},e\in[e_a,e_b]$ and $f_{m_{\beta}}$, such that $(h_1,f_{m_{\tau}},f_{m_{\beta}})\in\tilde{\mathcal{H}}$) and
		\begin{align}
		&\big|k_c(h_1,F_{Y|\bX}^{-1}(\tau|\cdot),f_{m_{\beta}},v,e)-k_c(s_0,v,e)-D_{f_{m_{\beta}}}k_c(s_0,v,e)\big[f_{m_{\beta}}-F_{Y|\bX}^{-1}(\beta|\cdot)\big]\big|\nonumber
		\\[0,2cm]&\leq C||f_{m_{\beta}}-F_{Y|\bX}^{-1}(\beta|\cdot)||_{M_{\bX}}^2\label{proofasympBlemmac3sufficientbeta}
		\end{align}
		uniformly in $c\in B_{\delta},v\in M_{\bX},e\in[e_a,e_b]$ to handle the second and third summand, respectively.
		
		(\ref{Dhkc}) leads to
		\begin{align*}
		&k_c(\mathfrak{h},f_{m_{\tau}},f_{m_{\beta}},v,e)-k_c(h_1,f_{m_{\tau}},f_{m_{\beta}},v,e)-D_hk_c(s_0,v,e)[\mathfrak{h}-h_1]
		\\[0,2cm]&=\frac{1}{\sigma(v)}\Bigg(h_1(\mathfrak{h}^{-1}(z_c(\mathfrak{h},f_{m_{\tau}},f_{m_{\beta}})))^B-h_1(h_1^{-1}(z_c(h_1,f_{m_{\tau}},f_{m_{\beta}})))^B
		\\[0,2cm]&\quad-Bz_c(h_1,f_{m_{\tau}},f_{m_{\beta}})^{B-1}\bigg(\frac{\partial}{\partial t}z_c(f_{h,t},f_{m_{\tau}},f_{m_{\beta}})\bigg|_{t=0}+h_1(h_1^{-1}(z_c(h_1,f_{m_{\tau}},f_{m_{\beta}})))
		\\[0,2cm]&\quad-\mathfrak{h}(h_1^{-1}(z_c(h_1,f_{m_{\tau}},f_{m_{\beta}})))\bigg)\Bigg)
		\\[0,2cm]&=\frac{1}{\sigma(v)}\bigg(h_1(\mathfrak{h}^{-1}(z_c(\mathfrak{h},f_{m_{\tau}},f_{m_{\beta}})))^B-h_1(h_1^{-1}(z_c(\mathfrak{h},f_{m_{\tau}},f_{m_{\beta}})))^B+z_c(\mathfrak{h},f_{m_{\tau}},f_{m_{\beta}})^B
		\\[0,2cm]&\quad-z_c(h_1,f_{m_{\tau}},f_{m_{\beta}})^B-B(z_c(h_1,f_{m_{\tau}},f_{m_{\beta}}))^{B-1}\bigg(\frac{\partial}{\partial t}z_c(f_{h,t},f_{m_{\tau}},f_{m_{\beta}})\bigg|_{t=0}
		\\[0,2cm]&\quad+h_1(h_1^{-1}(z_c(h_1,f_{m_{\tau}},f_{m_{\beta}})))-\mathfrak{h}(h_1^{-1}(z_c(h_1,f_{m_{\tau}},f_{m_{\beta}})))\bigg)\bigg)
		\\[0,2cm]&=\frac{1}{\sigma(v)}\Big(h_1(\mathfrak{h}^{-1}(z_c(\mathfrak{h},f_{m_{\tau}},f_{m_{\beta}})))^B-h_1(h_1^{-1}(z_c(\mathfrak{h},f_{m_{\tau}},f_{m_{\beta}})))^B
		\\[0,2cm]&\quad-Bz_c(h_1,f_{m_{\tau}},f_{m_{\beta}})^{B-1}(h_1(h_1^{-1}(z_c(h_1,f_{m_{\tau}},f_{m_{\beta}})))-\mathfrak{h}(h_1^{-1}(z_c(h_1,f_{m_{\tau}},f_{m_{\beta}}))))\Big)
		\\[0,2cm]&\quad+\mathcal{O}\big(||\mathfrak{h}-h_1||_{[z_a,z_b]}^2\big),
		\end{align*}
		because
		\begin{align*}
		z_c(\mathfrak{h},f_{m_{\tau}},f_{m_{\beta}})^B-z_c(0)^B-Bz_c(0)^{B-1}\frac{\partial}{\partial t}z_c(t)\bigg|_{t=0}&=\frac{\frac{\partial^2}{\partial t^2}z_c(f_{h,t},f_{m_{\tau}},f_{m_{\beta}})\big|_{t=\tilde{t}}}{2}
		\\[0,2cm]\quad&=\mathcal{O}\big(||\mathfrak{h}-h_1||_{[z_a,z_b]}^2\big)
		\end{align*}
		for an appropriate $\tilde{t}$ in $(0,1)$. Apply Lemma \ref{lemmahdevinv} to obtain
		\begin{align*}
		&h_1(\mathfrak{h}^{-1}(z_c(\mathfrak{h},f_{m_{\tau}},f_{m_{\beta}})))^B-h_1(h_1^{-1}(z_c(\mathfrak{h},f_{m_{\tau}},f_{m_{\beta}})))^B
		\\[0,2cm]&-Bz_c(h_1,f_{m_{\tau}},f_{m_{\beta}})^{B-1}(h_1(h_1^{-1}(z_c(h_1,f_{m_{\tau}},f_{m_{\beta}})))-\mathfrak{h}(h_1^{-1}(z_c(h_1,f_{m_{\tau}},f_{m_{\beta}}))))
		\\[0,2cm]&=Bh_1(h_1^{-1}(z_c(\mathfrak{h},f_{m_{\tau}},f_{m_{\beta}})))^{B-1}(h_1(\mathfrak{h}^{-1}(z_c(\mathfrak{h},f_{m_{\tau}},f_{m_{\beta}})))-h_1(h_1^{-1}(z_c(\mathfrak{h},f_{m_{\tau}},f_{m_{\beta}}))))
		\\[0,2cm]&\quad-Bz_c(h_1,f_{m_{\tau}},f_{m_{\beta}})^{B-1}(h_1(h_1^{-1}(z_c(h_1,f_{m_{\tau}},f_{m_{\beta}})))-\mathfrak{h}(h_1^{-1}(z_c(h_1,f_{m_{\tau}},f_{m_{\beta}}))))
		\\[0,2cm]&\quad+\mathcal{O}(||\mathfrak{h}-h_1||^2_{[z_a,z_b]})
		\\[0,2cm]&=Bz_c(h_1,f_{m_{\tau}},f_{m_{\beta}})^{B-1}\Big(h_1(\mathfrak{h}^{-1}(z_c(\mathfrak{h},f_{m_{\tau}},f_{m_{\beta}})))-h_1(h_1^{-1}(z_c(\mathfrak{h},f_{m_{\tau}},f_{m_{\beta}})))
		\\[0,2cm]&\quad-h_1(h_1^{-1}(z_c(h_1,f_{m_{\tau}},f_{m_{\beta}})))+\mathfrak{h}(h_1^{-1}(z_c(h_1,f_{m_{\tau}},f_{m_{\beta}})))\Big)+\mathcal{O}(||\mathfrak{h}-h_1||^2_{[z_a,z_b]})
		\end{align*}
		as well as
		\begin{align*}
		&h_1(\mathfrak{h}^{-1}(z_c(\mathfrak{h},f_{m_{\tau}},f_{m_{\beta}})))-h_1(h_1^{-1}(z_c(\mathfrak{h},f_{m_{\tau}},f_{m_{\beta}})))-h_1(h_1^{-1}(z_c(h_1,f_{m_{\tau}},f_{m_{\beta}})))
		\\[0,2cm]&+\mathfrak{h}(h_1^{-1}(z_c(h_1,f_{m_{\tau}},f_{m_{\beta}})))
		\\[0,2cm]&=h_1(\mathfrak{h}^{-1}(z_c(\mathfrak{h},f_{m_{\tau}},f_{m_{\beta}})))-h_1(h_1^{-1}(z_c(h_1,f_{m_{\tau}},f_{m_{\beta}})))
		\\[0,2cm]&\quad+\mathfrak{h}(h_1^{-1}(z_c(h_1,f_{m_{\tau}},f_{m_{\beta}})))-\mathfrak{h}(\mathfrak{h}^{-1}(z_c(\mathfrak{h},f_{m_{\tau}},f_{m_{\beta}})))
		\\[0,2cm]&=h_1'(h_1^{-1}(z_c(0)))(\mathfrak{h}^{-1}(z_c(\mathfrak{h},f_{m_{\tau}},f_{m_{\beta}}))-h_1^{-1}(z_c(h_1,f_{m_{\tau}},f_{m_{\beta}})))
		\\[0,2cm]&\quad+\mathfrak{h}'(h_1^{-1}(z_c(0)))(h_1^{-1}(z_c(h_1,f_{m_{\tau}},f_{m_{\beta}}))-\mathfrak{h}^{-1}(z_c(\mathfrak{h},f_{m_{\tau}},f_{m_{\beta}})))+\mathcal{O}(||\mathfrak{h}-h_1||^2_{[z_a,z_b]})
		\\[0,2cm]&=\big(h_1'(h_1^{-1}(z_c(0)))-\mathfrak{h}'(h_1^{-1}(z_c(0)))\big)\big(\mathfrak{h}^{-1}(z_c(\mathfrak{h},f_{m_{\tau}},f_{m_{\beta}}))-h_1^{-1}(z_c(h_1,f_{m_{\tau}},f_{m_{\beta}}))\big)
		\\[0,2cm]&\quad+\mathcal{O}(||\mathfrak{h}-h_1||^2_{[z_a,z_b]})
		\\[0,2cm]&=\mathcal{O}(||\mathfrak{h}-h_1||^{\frac{3}{2}}_{[z_a,z_b]}).
		\end{align*}
		It remains to treat the second and the third summand. Recall that it is sufficient to prove the equations (\ref{proofasympBlemmac3sufficientalpha}) and (\ref{proofasympBlemmac3sufficientbeta}). For that purpose, notice that
		\begin{align*}
		&z_c(h_1,f_{m_{\tau}},f_{m_{\beta}})-z_c(h_1,F_{Y|\bX}^{-1}(\tau|\cdot),f_{m_{\beta}})
		\\[0,2cm]&=\big(h_1(f_{m_{\tau}}(v))^c+e(h_1(f_{m_{\beta}}(v))^c-h_1(f_{m_{\tau}}(v))^c)\big)^{\frac{1}{c}}
		\\[0,2cm]&\quad-\big(h_1(F_{Y|\bX}^{-1}(\tau|v))^c+e(h_1(f_{m_{\beta}}(v))^c-h_1(F_{Y|\bX}^{-1}(\tau|v))^c)\big)^{\frac{1}{c}}
		\\[0,2cm]&=\frac{1}{c}\big(h_1(F_{Y|\bX}^{-1}(\tau|v))^c+e(h_1(f_{m_{\beta}}(v))^c-h_1(F_{Y|\bX}^{-1}(\tau|v))^c)\big)^{\frac{1}{c}-1}
		\\[0,2cm]&\quad\quad ch_1(F_{Y|\bX}^{-1}(\tau|v))^{c-1}(1-e)h_1'(F_{Y|\bX}^{-1}(\tau|v))(f_{m_{\tau}}(v)-F_{Y|\bX}^{-1}(\tau|v))
		\\[0,2cm]&\quad+\mathcal{O}(||f_{m_{\tau}}-F_{Y|\bX}^{-1}(\tau|\cdot)||_{M_{\bX}}^2)
		\\[0,2cm]&=-\frac{(1-e)z_c(h_1,F_{Y|\bX}^{-1}(\tau|\cdot),f_{m_{\beta}})^{1-c}h_1(F_{Y|\bX}^{-1}(\tau|v))^c(f_{m_{\tau}}(v)-F_{Y|\bX}^{-1}(\tau|v))}{\lambda(F_{Y|\bX}^{-1}(\tau|v))}
		\\[0,2cm]&\quad+\mathcal{O}(||f_{m_{\tau}}-F_{Y|\bX}^{-1}(\tau|\cdot)||_{M_{\bX}}^2).
		\end{align*}
		Therefore,
		\begin{align*}
		&k_c(h_1,f_{m_{\tau}},f_{m_{\beta}},v,e)-k_c(h_1,F_{Y|\bX}^{-1}(\tau|\cdot),f_{m_{\beta}},v,e)-D_{f_{m_{\tau}}}k_c(s_0,v,e)\big[f_{m_{\tau}}-F_{Y|\bX}^{-1}(\tau|\cdot)\big]
		\\[0,2cm]&=\frac{1}{\sigma(v)}\bigg(z_c(h_1,f_{m_{\tau}},f_{m_{\beta}})^B-z_c(h_1,F_{Y|\bX}^{-1}(\tau|\cdot),f_{m_{\beta}})^B
		\\[0,2cm]&\quad-\frac{Bz_c(h_1,F_{Y|\bX}^{-1}(\tau|\cdot),f_{m_{\beta}})^{B-c}(1-e)h_1(F_{Y|\bX}^{-1}(\tau|v))^c(f_{m_{\tau}}(v)-F_{Y|\bX}^{-1}(\tau|v))}{\lambda(F_{Y|\bX}^{-1}(\tau|v))}\bigg)
		\\[0,2cm]&=\frac{1}{\sigma(v)}\bigg(Bz_c(h_1,F_{Y|\bX}^{-1}(\tau|\cdot),f_{m_{\beta}})^{B-1}(z_c(h_1,f_{m_{\tau}},f_{m_{\beta}})-z_c(h_1,F_{Y|\bX}^{-1}(\tau|\cdot),f_{m_{\beta}}))
		\\[0,2cm]&\quad-\frac{Bz_c(h_1,F_{Y|\bX}^{-1}(\tau|\cdot),f_{m_{\beta}})^{B-c}(1-e)h_1(F_{Y|\bX}^{-1}(\tau|v))^c(f_{m_{\tau}}(v)-F_{Y|\bX}^{-1}(\tau|v))}{\lambda(F_{Y|\bX}^{-1}(\tau|v))}\bigg)
		\\[0,2cm]&\quad+\mathcal{O}(||f_{m_{\tau}}-F_{Y|\bX}^{-1}(\tau|\cdot)||_{M_{\bX}}^2)
		\\[0,2cm]&=\mathcal{O}(||f_{m_{\tau}}-F_{Y|\bX}^{-1}(\tau|\cdot)||_{M_{\bX}}^2).
		\end{align*}
		Analogously,
		\begin{align*}
		&k_c(h_1,F_{Y|\bX}^{-1}(\tau|\cdot),f_{m_{\beta}},v,e)-k_c(s_0,v,e)-D_{f_{m_{\tau}}}k_c(s_0,v,e)\big[f_{m_{\beta}}-F_{Y|\bX}^{-1}(\beta|\cdot)\big]
		\\[0,2cm]&=\mathcal{O}(||f_{m_{\beta}}-F_{Y|\bX}^{-1}(\beta|\cdot)||_{M_{\bX}}^2).
		\end{align*}
		Hence, $(i)$ is proven.\\[0,2cm]
		\textbf{Proof of }$\mathbf{(ii)}$: Remember \ref{C4} and let $c\in B_{\delta_n}$. As before, one has
		\begin{align*}
		&||\Gamma_2(c,s_0)(\bx,e)[\hat{s}-s_0]-\Gamma_2(B,s_0)(\bx,e)[\hat{s}-s_0]||
		\\[0,2cm]&=||D_hG_{MD}(c,s_0)(\bx,e)[\bar{h}_1-h_1]-D_hG_{MD}(B,s_0)(\bx,e)[\bar{h}_1-h_1]
		\\[0,2cm]&\quad+D_{f_{m_{\tau}}}G_{MD}(c,s_0)(\bx,e)[\hat{f}_{m_{\tau}}-F_{Y|\bX}^{-1}(\tau|\cdot)\big]
		\\[0,2cm]&\quad-D_{f_{m_{\tau}}}G_{MD}(B,s_0)(\bx,e)[\hat{f}_{m_{\tau}}-F_{Y|\bX}^{-1}(\tau|\cdot)\big]
		\\[0,2cm]&\quad+D_{f_{m_{\beta}}}G_{MD}(c,s_0)(\bx,e)[\hat{f}_{m_{\beta}}-F_{Y|\bX}^{-1}(\beta|\cdot)\big]
		\\[0,2cm]&\quad-D_{f_{m_{\beta}}}G_{MD}(B,s_0)(\bx,e)[\hat{f}_{m_{\beta}}-F_{Y|\bX}^{-1}(\beta|\cdot)\big]||
		\\[0,2cm]&\leq||D_hG_{MD}(c,s_0)(\bx,e)[\bar{h}_1-h_1]-D_hG_{MD}(B,s_0)(\bx,e)[\bar{h}_1-h_1]||
		\\[0,2cm]&\quad+||D_{f_{m_{\tau}}}G_{MD}(c,s_0)(\bx,e)[\hat{f}_{m_{\tau}}-F_{Y|\bX}^{-1}(\tau|\cdot)\big]
		\\[0,2cm]&\quad-D_{f_{m_{\tau}}}G_{MD}(B,s_0)(\bx,e)[\hat{f}_{m_{\tau}}-F_{Y|\bX}^{-1}(\tau|\cdot)\big]||
		\\[0,2cm]&\quad+||D_{f_{m_{\beta}}}G_{MD}(c,s_0)(\bx,e)[\hat{f}_{m_{\beta}}-F_{Y|\bX}^{-1}(\beta|\cdot)\big]
		\\[0,2cm]&\quad-D_{f_{m_{\beta}}}G_{MD}(B,s_0)(\bx,e)[\hat{f}_{m_{\beta}}-F_{Y|\bX}^{-1}(\beta|\cdot)\big]||
		\end{align*}
		so that it is again sufficient to prove the condition for each of the summands. To treat the first summand, let $\mathfrak{h}\in\tilde{\mathcal{H}},c\in[B_1,B_2]$ and recall the definitions of $f_{h,t}$ and $z_c(t)$ from (\ref{deffht}) and (\ref{defzc}) as well as
		$$|z_c(0)-z_B(0)|=\mathcal{O}(|c-B|)\quad\textup{uniformly in }(\bx,e)$$
		and
		\begin{align*}
		&D_hG_{MD}(c,s_0)(\bx,e)[\mathfrak{h}-h_1]
		\\[0,2cm]&=\frac{1}{P(\bX\in M_{\bX})}\frac{\partial}{\partial t}\Big(\int_{M_{\bX}}F_{\varepsilon}(k_c(f_{h,t},F_{Y|\bX}^{-1}(\tau|\cdot),F_{Y|\bX}^{-1}(\beta|\cdot),\mathbf{w},e))I_{\{\mathbf{w}\leq \bx\}}f_{\bX}(\mathbf{w})\,d\mathbf{w}
		\\[0,2cm]&\quad-P(\bX\leq \bx|\bX\in M_{\bX})\int_{M_{\bX}}F_{\varepsilon}(k_c(f_{h,t},F_{Y|\bX}^{-1}(\tau|\cdot),F_{Y|\bX}^{-1}(\beta|\cdot),\mathbf{w},e))f_{\bX}(\mathbf{w})\,d\mathbf{w}\Big)\bigg|_{t=0}
		\\[0,2cm]&=\frac{1}{P(\bX\in M_{\bX})}\Big(\int_{M_{\bX}}f_{\varepsilon}(k_c(s_0,\mathbf{w},e))\big(I_{\{\mathbf{w}\leq \bx\}}-P(\bX\leq \bx|\bX\in M_{\bX})\big)
		\\[0,2cm]&\quad D_hk_c(s_0,\mathbf{w},e)[\mathfrak{h}-h_1]f_{\bX}(\mathbf{w})\,d\mathbf{w}\Big).
		\end{align*}
		At the beginning of the proof of this lemma, it was shown in (\ref{Dhkc}) that
		$$D_hk_c(s_0,\bx,e)[\mathfrak{h}-h_1]=\frac{Bz_c(0)^{B-1}\big(z_c(0)-\mathfrak{h}(h_1^{-1}(z_c(0)))+\frac{\partial}{\partial t}z_c(t)\big|_{t=0}\big)}{\sigma(\bx)}$$
		with $f_{m_{\tau}}=F_{Y|\bX}^{-1}(\tau|\cdot),f_{m_{\beta}}=F_{Y|\bX}^{-1}(\beta|\cdot)$ in $z_c(t)$), where
		\begin{align*}
		\frac{\partial}{\partial t}z_c(t)\bigg|_{t=0}&=\frac{1}{c}\big(h_1^c(F_{Y|\bX}^{-1}(\tau|\bx))+e(h_1^c(F_{Y|\bX}^{-1}(\beta|\bx))-h_1^c(F_{Y|\bX}^{-1}(\tau|\bx)))\big)^{\frac{1}{c}-1}
		\\[0,2cm]&\quad\big(ch_1^{c-1}(F_{Y|\bX}^{-1}(\tau|\bx))(\mathfrak{h}(F_{Y|\bX}^{-1}(\tau|\bx))-h_1(F_{Y|\bX}^{-1}(\tau|\bx))
		)
		\\[0,2cm]&\quad+e(ch_1^{c-1}(F_{Y|\bX}^{-1}(\beta|\bx))(\mathfrak{h}(F_{Y|\bX}^{-1}(\beta|\bx))-h_1(F_{Y|\bX}^{-1}(\beta|\bx))
		\\[0,2cm]&\quad\quad-ch_1^{c-1}(F_{Y|\bX}^{-1}(\tau|\bx))(\mathfrak{h}(F_{Y|\bX}^{-1}(\tau|\bx))-h_1(F_{Y|\bX}^{-1}(\tau|\bx)))\big).
		\end{align*}
		Hence,
		$$\underset{\bx\in M_{\bX},e\in[e_a,e_b]}{\sup}\,\bigg|\frac{\partial}{\partial t}z_c(t)\bigg|_{t=0}\bigg|=\mathcal{O}(||\mathfrak{h}-h_1||_{[z_a,z_b]})$$
		and
		$$\underset{\bx\in M_{\bX},e\in[e_a,e_b]}{\sup}\,\bigg|\frac{\partial}{\partial t}z_c(t)\bigg|_{t=0}-\frac{\partial}{\partial t}z_B(t)\bigg|_{t=0}\bigg|=\mathcal{O}\big(||\mathfrak{h}-h_1||_{[z_a,z_b]}|c-B|\big),$$
		so that
		\begin{align*}
		&D_hG_{MD}(c,s_0)(\bx,e)[\mathfrak{h}-h_1]-D_hG_{MD}(B,s_0)(\bx,e)[\mathfrak{h}-h_1]
		\\[0,2cm]&=\int_{M_{\bX}}\big(I_{\{\mathbf{w}\leq \bx\}}-P(\bX\leq \bx|\bX\in M_{\bX})\big)\Big(\varphi(c,\mathbf{w},e)\big(z_c(0)-\mathfrak{h}(h_1^{-1}(z_c(0)))\big)
		\\[0,2cm]&\quad\quad-\varphi(B,\mathbf{w},e)\big(z_B(0)-\mathfrak{h}(h_1^{-1}(z_B(0)))\big)\Big)\,d\mathbf{w}+o\big(|c-B|\big)
		\end{align*}
		for some continuously differentiable function $\varphi:[B_1,B_2]\times M_{\bX}\times[e_a,e_b]\rightarrow\mathbb{R}$. Due to
		\begin{align*}
		&\varphi(c,\mathbf{w},e)\big(h_1(h_1^{-1}(z_c(0)))-\mathfrak{h}(h_1^{-1}(z_c(0)))\big)
		\\[0,2cm]&-\varphi(B,\mathbf{w},e)\big(h_1(h_1^{-1}(z_B(0)))-\mathfrak{h}(h_1^{-1}(z_B(0)))\big)
		\\[0,2cm]&\quad=(\varphi(c,\mathbf{w},e)-\varphi(B,\mathbf{w},e))\big(h_1(h_1^{-1}(z_c(0)))-\mathfrak{h}(h_1^{-1}(z_c(0)))\big)
		\\[0,2cm]&\quad\quad+\varphi(B,\mathbf{w},e)\big(h_1(h_1^{-1}(z_c(0)))-h_1(h_1^{-1}(z_B(0)))+\mathfrak{h}(h_1^{-1}(z_B(0)))-\mathfrak{h}(h_1^{-1}(z_c(0)))\big)
		\\[0,2cm]&\quad=\varphi(B,\mathbf{w},e)\big(h_1'(h_1^{-1}(z_B(0)))(h_1^{-1}(z_c(0))-h_1^{-1}(z_B(0)))
		\\[0,2cm]&\quad\quad-\mathfrak{h}'(h_1^{-1}(z_B(0)))(h_1^{-1}(z_c(0))-h_1^{-1}(z_B(0)))\big)+o(|c-B|)
		\\[0,2cm]&\quad=\mathcal{O}(||\mathfrak{h}'-h_1'||_{[z_a,z_b]}|c-B|)+o(|c-B|)
		\\[0,2cm]&\quad=o(|c-B|),
		\end{align*}
		it holds that
		\begin{align*}
		||D_hG_{MD}(c,s_0)(\bx,e)[\bar{h}_1-h_1]-D_hG_{MD}(B,s_0)(\bx,e)[\bar{h}_1-h_1]||&=o_p(|c-B|).
		\end{align*}
		The second summand can be written as
		\begin{align*}
		&||D_{f_{m_{\tau}}}G_{MD}(c,s_0)(\bx,e)[\hat{F}_{Y|\bX}^{-1}(\tau|\cdot)-F_{Y|\bX}^{-1}(\tau|\cdot)\big]
		\\[0,2cm]&-D_{f_{m_{\tau}}}G_{MD}(B,s_0)(\bx,e)[\hat{F}_{Y|\bX}^{-1}(\tau|\cdot)-F_{Y|\bX}^{-1}(\tau|\cdot)\big]||
		\\[0,2cm]&=\bigg|\bigg|\frac{1}{P(\bX\in M_{\bX})}\int_{M_{\bX}}\big(I_{\{\mathbf{w}\leq\cdot\}}-P(\bX\leq\cdot|\bX\in M_{\bX})\big)
		\\[0,2cm]&\quad\quad\big(f_{\varepsilon}(k_c(s_0,\mathbf{w},.))D_{f_{m_{\tau}}}k_c(s_0,\mathbf{w},.)[\hat{F}_{Y|\bX}^{-1}(\tau|\cdot)-F_{Y|\bX}^{-1}(\tau|\cdot)\big]
		\\[0,2cm]&\quad\quad-f_{\varepsilon}(k_B(s_0,\mathbf{w},.))D_{f_{m_{\tau}}}k_B(s_0,\mathbf{w},.)[\hat{F}_{Y|\bX}^{-1}(\tau|\cdot)-F_{Y|\bX}^{-1}(\tau|\cdot)\big]\big)f_{\bX}(\mathbf{w})\,d\mathbf{w}\bigg|\bigg|.
		\end{align*}
		Thus, it is sufficient to prove
		\begin{align*}
		&f_{\varepsilon}(k_c(s_0,\mathbf{w},e))D_{f_{m_{\tau}}}k_c(s_0,\mathbf{w},e)[\hat{F}_{Y|\bX}^{-1}(\tau|\cdot)-F_{Y|\bX}^{-1}(\tau|\cdot)\big]
		\\[0,2cm]&-f_{\varepsilon}(k_B(s_0,\mathbf{w},e))D_{f_{m_{\tau}}}k_B(s_0,\mathbf{w},e)[\hat{F}_{Y|\bX}^{-1}(\tau|\cdot)-F_{Y|\bX}^{-1}(\tau|\cdot)\big]
		\\[0,2cm]&=\big(f_{\varepsilon}(k_c(s_0,\mathbf{w},e))-f_{\varepsilon}(k_B(s_0,\mathbf{w},e))\big)D_{f_{m_{\tau}}}k_c(s_0,\mathbf{w},e)[\hat{F}_{Y|\bX}^{-1}(\tau|\cdot)-F_{Y|\bX}^{-1}(\tau|\cdot)\big]
		\\[0,2cm]&\quad+f_{\varepsilon}(k_B(s_0,\mathbf{w},e))\big(D_{f_{m_{\tau}}}k_c(s_0,\mathbf{w},e)[\hat{F}_{Y|\bX}^{-1}(\tau|\cdot)-F_{Y|\bX}^{-1}(\tau|\cdot)\big]
		\\[0,2cm]&\quad\quad-D_{f_{m_{\tau}}}k_B(s_0,\mathbf{w},e)[\hat{F}_{Y|\bX}^{-1}(\tau|\cdot)-F_{Y|\bX}^{-1}(\tau|\cdot)\big]\big)
		\\[0,2cm]&=o_p(|c-B|)+\mathcal{O}_p\big(n^{-\frac{1}{2}}\big)
		\end{align*}
		uniformly in $(\mathbf{w},e)\in M_{\bX}\times[e_a,e_b]$. By condition \ref{C4}, $\hat{f}_{m_{\tau}}(\bx)-F_{Y|\bX}^{-1}(\tau|\bx)=o_p(1)$ uniformly in $\bx\in M_{\bX}$, so that for an appropriate $\tilde{c}$ between $c$ and $B$
		\begin{align*}
		&\big(f_{\varepsilon}(k_c(s_0,\mathbf{w},e))-f_{\varepsilon}(k_B(s_0,\mathbf{w},e))\big)D_{f_{m_{\tau}}}k_c(s_0,\mathbf{w},e)[\hat{f}_{m_{\tau}}-F_{Y|\bX}^{-1}(\tau|\cdot)\big]
		\\[0,2cm]&=f_{\varepsilon}'(k_{\tilde{c}}(s_0,\mathbf{w},e))\frac{\partial}{\partial c}k_c(s_0,\mathbf{w},e)\bigg|_{c=\tilde{c}}(c-B)o_p(1)
		\\[0,2cm]&=o_p(|c-B|).
		\end{align*}
		On the other hand, the remaining term can be rewritten via
		\begin{align*}
		&D_{f_{m_{\tau}}}k_c(s_0,\mathbf{w},e)[[\hat{f}_{m_{\tau}}-F_{Y|\bX}^{-1}(\tau|\cdot)\big]-D_{f_{m_{\tau}}}k_B(s_0,\mathbf{w},.)[\hat{f}_{m_{\tau}}-F_{Y|\bX}^{-1}(\tau|\cdot)\big]
		\\[0,2cm]&=\frac{B(e-1)(\hat{f}_{m_{\tau}}(\mathbf{w})-F_{Y|\bX}^{-1}(\tau|\mathbf{w}))}{\sigma(\mathbf{w})\lambda(F_{Y|\bX}^{-1}(\tau|\mathbf{w}))}(\psi(B,\mathbf{w},e)-\psi(c,\mathbf{w},e)),
		\end{align*}
		where
		$$\psi(c,\mathbf{w},e)=\big(h_c(F_{Y|\bX}^{-1}(\tau|\mathbf{w}))+e\big(h_c(F_{Y|\bX}^{-1}(\beta|\mathbf{w}))-h_c(F_{Y|\bX}^{-1}(\tau|\mathbf{w}))\big)\big)^{\frac{B}{c}-1}h_c(F_{Y|\bX}^{-1}(\tau|\bx)).$$
		Due to
		\begin{align*}
		&\frac{\partial}{\partial c}\psi(c,\mathbf{w},e)
		\\[0,2cm]&=\frac{\partial}{\partial c}\big(h_c(F_{Y|\bX}^{-1}(\tau|\mathbf{w}))+e\big(h_c(F_{Y|\bX}^{-1}(\beta|\mathbf{w}))-h_c(F_{Y|\bX}^{-1}(\tau|\mathbf{w}))\big)\big)^{\frac{B}{c}-1}h_c(F_{Y|\bX}^{-1}(\tau|\bx))
		\\[0,2cm]&=-\frac{B}{c^2}\log\big(h_c(F_{Y|\bX}^{-1}(\tau|\mathbf{w}))+e\big(h_c(F_{Y|\bX}^{-1}(\beta|\mathbf{w}))-h_c(F_{Y|\bX}^{-1}(\tau|\mathbf{w}))\big)\big)\psi(c,\mathbf{w},e)
		\\[0,2cm]&\quad+\bigg(\frac{B}{c}-1\bigg)h_c(F_{Y|\bX}^{-1}(\tau|\bx))
		\\[0,2cm]&\quad\quad\Big(h_c(F_{Y|\bX}^{-1}(\tau|\mathbf{w}))+e\big(h_c(F_{Y|\bX}^{-1}(\beta|\mathbf{w}))-h_c(F_{Y|\bX}^{-1}(\tau|\mathbf{w}))\big)\Big)^{\frac{B}{c}-2}
		\\[0,2cm]&\quad\quad\bigg(\log\big(h_1(F_{Y|\bX}^{-1}(\tau|\mathbf{w}))\big)h_c(F_{Y|\bX}^{-1}(\tau|\mathbf{w}))+e\Big(\log\big(h_1(F_{Y|\bX}^{-1}(\beta|\mathbf{w}))\big)h_c(F_{Y|\bX}^{-1}(\beta|\mathbf{w}))
		\\[0,2cm]&\quad\quad-\log\big(h_1(F_{Y|\bX}^{-1}(\tau|\mathbf{w}))\big)h_c(F_{Y|\bX}^{-1}(\tau|\mathbf{w}))\Big)\bigg)+\psi(c,\mathbf{w},e)\log\big(h_1(F_{Y|\bX}^{-1}(\beta|\mathbf{w}))\big),
		\end{align*}
		the derivative of $\psi$ with respect to $c$ is uniformly bounded in $(\mathbf{w},e)\in M_{\bX}\times[e_a,e_b]$. Hence,
		\begin{align*}
		&D_{f_{m_{\tau}}}k_c(s_0,\mathbf{w},e)[\hat{f}_{m_{\tau}}-F_{Y|\bX}^{-1}(\tau|\cdot)\big]-D_{f_{m_{\tau}}}k_B(s_0,\mathbf{w},e)[\hat{f}_{m_{\tau}}-F_{Y|\bX}^{-1}(\tau|\cdot)\big]
		\\[0,2cm]&=o_p(|c-B|)
		\end{align*}
		uniformly in $(\mathbf{w},e)\in M_{\bX}\times[e_a,e_b]$. The same reasoning can be applied for
		$$D_{f_{m_{\beta}}}G_{MD}(c,s_0)(\bx,e)[\hat{f}_{m_{\beta}}-F_{Y|\bX}^{-1}(\beta|\cdot)\big],$$
		which completes the proof of Lemma \ref{lemmac3}.

		\begin{lemma}\label{lemmac5}
			For all sequences $\delta_n\searrow0$ it holds that
			$$\underset{||c-B||\leq\delta_n,||s-s_0||\leq\delta_n}{\sup}\,||G_{nMD}(c,s)-G_{MD}(c,s)-G_{nMD}(B,s_0)||_2=o_p(n^{-\frac{1}{2}}),$$
			that is, \ref{C5} is valid.
		\end{lemma}\noindent
		\textbf{Proof:} In  a moment, it will be shown that the process
		$$G_n(c,s,\bx,e)=G_{nMD}(c,s)(\bx,e)-G_{MD}(c,s)(\bx,e),$$
		as a process in $c\in[B_1,B_2],s\in\tilde{\mathcal{H}},\bx\in M_{\bX},e\in[e_a,e_b]$ is Donsker. Then, Corollary 2.3.12 of \citet{vdVW1996} yields
		$$\underset{\overset{\scriptstyle c,\tilde{c}\in[B_1,B_2],s,\tilde{s}\in\tilde{\mathcal{H}},\bx\in M_{\bX},e\in[e_a,e_b]}{||s-\tilde{s}||_{\mathcal{H}}<\delta_n,|c-\tilde{c}|<\delta_n}}{\sup}\,\sqrt{n}|G_n(c,s,\bx,e)-G_n(\tilde{c},\tilde{s},\bx,e)|=o_p(1).$$
		Due to $G_{MD}(B,s_0)(\bx,e)=0$ for all $\bx\in M_{\bX},e\in[e_a,e_b]$ the assertion then follows from the compactness of $M_{\bX}$ and $[e_a,e_b]$.
		
		First, define the function class
		$$\mathcal{F}=\{(\bX,\varepsilon)\mapsto I_{\{\bX\in M_{\bX}\}}I_{\{\varepsilon\leq k_c(\mathfrak{h},f_{m_{\tau}},f_{m_{\beta}},\bX,e)\}}:s\in\tilde{\mathcal{H}},c\in[B_1,B_2],e\in[e_a,e_b]\}.$$
		Due to the definition of $\tilde{\mathcal{H}}$ in (\ref{defHcaltilde}) and the compactness of $[B_1,B_2],M_{\bX},[e_a,e_b]$, there exists a compact set $\mathcal{K}$ such that
		$$k_c(s,\bx,e)\in\mathcal{K},\quad\textup{for all }s\in\tilde{\mathcal{H}},\bx\in M_{\bX},e\in[e_a,e_b].$$
		Consider $s,\tilde{s}\in\tilde{\mathcal{H}},c,\tilde{c}\in[B_1,B_2],e,\tilde{e}\in[e_a,e_b]$. For some $s^*,c^*$ and $e^*$ between $s$ and $\tilde{s}$, $c$ and $\tilde{c}$ and $e$ and $\tilde{e}$, respectively, as well as some $C>0$ the $L^2(P^{Y,\bX})$-distance can be bounded by
		\begin{align*}
		&||I_{\{\cdot\in M_{\bX}\}}I_{\{\cdot\leq k_c(s,\cdot,e)\}}-I_{\{\cdot\in M_{\bX}\}}I_{\{\cdot\leq k_{\tilde{c}}(\tilde{s},\cdot,\tilde{e})\}}||_2
		\\[0,2cm]&=E\big[I_{\{\bX\in M_{\bX}\}}\big(I_{\{\varepsilon\leq k_c(s,\bX,e)\}}-I_{\{\varepsilon\leq k_{\tilde{c}}(\tilde{s},\bX,\tilde{e})\}}\big)^2\big]^{\frac{1}{2}}
		\\[0,2cm]&=\bigg(\int_{M_{\bX}}|F_{\varepsilon}(k_c(s,\mathbf{w},e))-F_{\varepsilon}(k_{\tilde{c}}(\tilde{s},\mathbf{w},\tilde{e}))|f_{\bX}(\mathbf{w})\,d\mathbf{w}\bigg)^{\frac{1}{2}}
		\\[0,2cm]&\leq\underset{e\in\mathcal{K}}{\sup}\,|f_{\varepsilon}(e)|\bigg(\int_{M_{\bX}}|f_{\bX}(\mathbf{w})|\,d\mathbf{w}\bigg)^{\frac{1}{2}}\,\underset{\mathbf{w}\in M_{\bX}}{\sup}\,\big|k_c(s,\mathbf{w},e)-k_{\tilde{c}}(\tilde{s},\mathbf{w},\tilde{e})\big|^{\frac{1}{2}}
		\\[0,2cm]&\leq\underset{e\in\mathcal{K}}{\sup}\,|f_{\varepsilon}(e)|\bigg(\int_{M_{\bX}}|f_{\bX}(\mathbf{w})|\,d\mathbf{w}\bigg)^\frac{1}{2}\,\underset{\mathbf{w}\in M_{\bX}}{\sup}\,\Big|D_hk_c(s^*,\mathbf{w},e^*)[\tilde{\mathfrak{h}}-\mathfrak{h}]
		\\[0,2cm]&\quad\quad+D_{f_{m_{\tau}}}k_c(s^*,\mathbf{w},e^*)[\tilde{f}_{m_{\tau}}-f_{m_{\tau}}]+D_{f_{m_{\beta}}}k_c(s^*,\mathbf{w},e^*)[\tilde{f}_{m_{\beta}}-f_{m_{\beta}}]
		\\[0,2cm]&\quad\quad+D_{e}k_c(s^*,\mathbf{w},e^*)[\tilde{e}-e]+D_{c}k_{c^*}(s^*,\mathbf{w},e^*)[\tilde{c}-c]\Big|^{\frac{1}{2}}.
		\end{align*}
		Similar to the the proof of (\ref{Dhkc}), one can show with
		\begin{align*}
		s^*&=(\mathfrak{h}^*,f_{m_{\tau}}^*,f_{m_{\beta}}^*),
		\\[0,2cm]f_{h,t}&=\mathfrak{h}^*+t(\tilde{\mathfrak{h}}-\mathfrak{h}),
		\\[0,2cm]\tilde{z}_c(t)&=\big(f_{h,t}^c(f_{m_{\tau}}^*(\bx))+e(f_{h,t}^c(f_{m_{\beta}}^*(\bx))-f_{h,t}^c(f_{m_{\tau}}^*(\bx)))\big)^{\frac{1}{c}}
		\end{align*}
		that
		\begin{align*}
		&D_hk_c(s^*,\mathbf{w},e)[\tilde{\mathfrak{h}}-\mathfrak{h}]
		\\[0,2cm]&=\frac{Bh_1\big((h^*)^{-1}(\tilde{z}_c(0))\big)^{B-1}h_1'\big((h^*)^{-1}(z_c(0))\big)\Big(\frac{\partial}{\partial t}\tilde{z}_c(t)\big|_{t=0}-\big(\tilde{\mathfrak{h}}-\mathfrak{h}\big)\big((h^*)^{-1}(z_c(0))\big)\Big)}{(h^*)'\big(h_1^{-1}(z_c(0))\big)\sigma(\mathbf{w})}
		\end{align*}
		and that
		$$\underset{c\in[B_1,B_2],\mathbf{w}\in M_{\bX},e\in[e_a,e_b]}{\sup}\,|D_hk_c(s^*,\mathbf{w},e)[\tilde{\mathfrak{h}}-\mathfrak{h}]|\leq \tilde{C}||\tilde{\mathfrak{h}}-\mathfrak{h}||_{[z_a,z_b]}$$
		for an appropriate constant $\tilde{C}>0$ and all $s,\tilde{s}\in\tilde{\mathcal{H}}$ and $s^*$ between $s$ and $\tilde{s}$. A similar reasoning for $D_{f_{m_{\tau}}}k_{c^*}(s^*,\mathbf{w},e^*)[\tilde{f}_{m_{\tau}}-f_{m_{\tau}}],...,D_{c^*}k_{c^*}(s^*,\mathbf{w},e^*)[\tilde{c}-c]$ leads to
		\begin{align}
		&\underset{\mathbf{w}\in M_{\bX}}{\sup}\,\big|k_c(s,\mathbf{w},e)-k_{\tilde{c}}(\tilde{s},\mathbf{w},\tilde{e})\big|\nonumber
		\\[0,2cm]&\quad\quad\leq \bar{C}\Big(||\tilde{\mathfrak{h}}-\mathfrak{h}||_{[z_a,z_b]}+||\tilde{f}_{m_{\tau}}-f_{m_{\tau}}||_{M_{\bX}}+||\tilde{f}_{m_{\beta}}-f_{m_{\beta}}||_{M_{\bX}}+|\tilde{e}-e|+|\tilde{c}-c|\Big)\label{auxbracket}
		\end{align}
		for some appropriate constant $\bar{C}>0$, which is independent of $c^*,s^*,e^*$. This will be used in the following to define brackets for $\mathcal{F}$. Let $\xi,\eta>0$ and consider $\xi$-brackets for $c\in[B_1,B_2],e\in[e_a,e_b]$ and $\mathfrak{h},f_{m_{\tau}},f_{m_{\beta}}$ such that $s\in\tilde{\mathcal{H}}$. Construct $\eta$-brackets for $\mathcal{F}$ as follows. Let
		$$\xi=\xi(\eta)=\frac{\eta^2}{10\bar{C}\,\underset{e\in\mathcal{K}}{\sup}\,f_{\varepsilon}(e)^2\,\int|f_{\bX}(\mathbf{w})|\,d\mathbf{w}}$$
		with $\bar{C}$ from (\ref{auxbracket}). For each combination of the $\xi$-brackets take representatives $\bar{\mathfrak{h}},\bar{f}_{m_{\tau}},\bar{f}_{m_{\beta}}$, $\bar{c},\bar{e}$ within these brackets and define
		$$l(\bX,\varepsilon)=I\{\bX\in M_{\bX}\}I\Bigg\{\varepsilon\leq k_{\bar{c}}(\bar{\mathfrak{h}},\bar{f}_{m_{\tau}},\bar{f}_{m_{\beta}},\bX,\bar{e})-\frac{\eta^2}{2\,\underset{e\in\mathcal{K}}{\sup}\,f_{\varepsilon}(e)^2\,\int|f_{\bX}(\mathbf{w})|\,d\mathbf{w}}\Bigg\}$$
		and
		$$u(\bX,\varepsilon)=I\{\bX\in M_{\bX}\}I\Bigg\{\varepsilon\leq k_{\bar{c}}(\bar{\mathfrak{h}},\bar{f}_{m_{\tau}},\bar{f}_{m_{\beta}},\bX,\bar{e})+\frac{\eta^2}{2\,\underset{e\in\mathcal{K}}{\sup}\,f_{\varepsilon}(e)^2\,\int|f_{\bX}(\mathbf{w})|\,d\mathbf{w}}\Bigg\}.$$
		Then, $||u-l||_2\leq\eta$ by the same reasoning as above and equation (\ref{auxbracket}) ensures that each combination of the $\xi(\eta)$-brackets for $\mathfrak{h},f_{m_{\tau}},f_{m_{\beta}}$ such that $s\in\tilde{\mathcal{H}}$ and $c\in[B_1,B_2],e\in[e_a,e_b]$ is covered by its corresponding $[l,u]$-bracket.
		
		Since $\tilde{\mathcal{H}}\subseteq C_{R_h}^{\gamma_h}([z_a,z_b])\times C_{R_{f_{m_{\tau}}}}^{\gamma_{f_{m_{\tau}}}}(M_{\bX})\times C_{R_{f_{m_{\tau}}}}^{\gamma_{f_{m_{\beta}}}}(M_{\bX})$ one has for all $\eta>0$
		\begin{align*}
		\mathcal{N}_{[\,]}(\eta,\mathcal{F},L^2(P))&\leq\mathcal{N}_{[\,]}\big(\xi(\eta),C_{R_h}^{\gamma_h}([z_a,z_b]),||.||_{[z_a,z_b]}\big)\mathcal{N}_{[\,]}\big(\xi(\eta),C_{R_{f_{m_{\tau}}}}^{\gamma_{f_{m_{\tau}}}}(M_{\bX}),||.||_{M_{\bX}}\big)
		\\[0,2cm]&\quad\quad\mathcal{N}_{[\,]}\big(\xi(\eta),C_{R_{f_{m_{\beta}}}}^{\gamma_{f_{m_{\beta}}}}(M_{\bX}),||.||_{M_{\bX}}\big)\mathcal{N}_{[\,]}\big(\xi(\eta),[e_a,e_b],|.|\big)
		\\[0,2cm]&\quad\quad\mathcal{N}_{[\,]}\big(\xi(\eta),[B_1,B_2],|.|\big).
		\end{align*}
		According to Theorem 2.7.1 of \citet{vdVW1996}, one has
		\begin{align*}
		\log\bigg(\mathcal{N}_{[\,]}\Big(\xi,C_{R_h}^{\gamma_h}([z_a,z_b]),||.||_{[z_a,z_b]}\Big)\bigg)&\leq C_h\xi^{-\frac{1}{\gamma_h}},
		\\[0,2cm]\log\bigg(\mathcal{N}_{[\,]}\Big(\xi,C_{R_{f_{m_{\tau}}}}^{\gamma_{f_{m_{\tau}}}}(M_{\bX}),||.||_{M_{\bX}}\Big)\bigg)&\leq C_{f_{m_{\tau}}}\xi^{-\frac{\dX}{\gamma_{f_{m_{\tau}}}}},
		\\[0,2cm]\log\bigg(\mathcal{N}_{[\,]}\Big(\xi,C_{R_{f_{m_{\beta}}}}^{\gamma_{f_{m_{\beta}}}}(M_{\bX}),||.||_{M_{\bX}}\Big)\bigg)&\leq C_{f_{m_{\beta}}}\xi^{-\frac{\dX}{\gamma_{f_{m_{\beta}}}}}
		\end{align*}
		for some appropriate constants $C_{h},C_{f_{m_{\tau}}},C_{f_{m_{\beta}}}>0$. Note that
		$$\gamma_h>1,\quad\gamma_{f_{m_{\tau}}}>\dX\quad\textup{and}\quad\gamma_{f_{m_{\beta}}}>\dX$$
		by definition of $\tilde{\mathcal{H}}$ in (\ref{defHcaltilde}). Hence, for some $C>0$
		\begin{align*}
		&\int_0^{1}\sqrt{\log\big(\mathcal{N}_{[\,]}(\eta,\mathcal{F},L^2(P))\big)}\,d\eta
		\\[0,2cm]&\leq\int_0^1\sqrt{\log\Big(\mathcal{N}_{[\,]}\big(\xi(\eta),C_{R_h}^{\gamma_h}([z_a,z_b]),||.||_{[z_a,z_b]}\big)\Big)}\,d\eta
		\\[0,2cm]&\quad+\int_0^1\sqrt{\log\Big(\mathcal{N}_{[\,]}\big(\xi(\eta),C_{R_{f_{m_{\tau}}}}^{\gamma_{f_{m_{\tau}}}}(M_{\bX}),||.||_{M_{\bX}}\big)\Big)}\,d\eta
		\\[0,2cm]&\quad+\int_0^1\sqrt{\log\Big(\mathcal{N}_{[\,]}\big(\xi(\eta),C_{R_{f_{m_{\beta}}}}^{\gamma_{f_{m_{\beta}}}}(M_{\bX}),||.||_{M_{\bX}}\big)\Big)}\,d\eta
		\\[0,2cm]&\quad+\int_0^1\sqrt{\log\Big(\mathcal{N}_{[\,]}\big(\xi(\eta),[e_a,e_b],|.|\big)\Big)}\,d\eta
		\\[0,2cm]&\quad+\int_0^1\sqrt{\log\Big(\mathcal{N}_{[\,]}\big(\xi(\eta),[B_1,B_2],|.|\big)\Big)}\,d\eta
		\\[0,2cm]&\leq C\int_0^1\Bigg(\bigg(\frac{1}{\eta}\bigg)^{\frac{1}{\gamma_h}}+\bigg(\frac{1}{\eta}\bigg)^{\frac{\dX}{\gamma_{f_{m_{\tau}}}}}+\bigg(\frac{1}{\eta}\bigg)^{\frac{\dX}{\gamma_{f_{m_{\beta}}}}}+\log\bigg(\frac{1}{\eta^2}\bigg)+\log\bigg(\frac{1}{\eta^2}\bigg)\Bigg)\,d\eta
		\\[0,2cm]&<\infty,
		\end{align*}
		so that the function class $\mathcal{F}$ is Donsker. Of course, the function class $\{\bX\mapsto I_{\{\bX\leq \bx\}}:\bx\in M_{\bX}\}$ is Donsker and by the same reasoning as before it can be shown that the class
		$$\tilde{\mathcal{F}}=\{(\bX,\varepsilon)\mapsto I_{\{\bX\leq \bx\}}I_{\{\bX\in M_{\bX}\}}I_{\{\varepsilon\leq k_c(s,\bX,e)\}}:s\in\tilde{\mathcal{H}},c\in[B_1,B_2],e\in[e_a,e_b],\bx\in M_{\bX}\}$$
		is Donsker as well. Finally, it can be shown that
		\begin{align*}
		&G_{nMD}(c,s)(\bx,e)-G_{MD}(c,s)(\bx,e)
		\\[0,2cm]&=\bigg(\frac{2P(\bX\leq \bx|\bX\in M_{\bX})P(\varepsilon\leq k_c(s,\bX,e)|\bX\in M_{\bX})}{P(\bX\in M_{\bX})}
		\\[0,2cm]&=-\frac{P(\bX\leq \bx,\varepsilon\leq k_c(s,\bX,e)|\bX\in M_{\bX})}{P(\bX\in M_{\bX})}\bigg)\bigg(\frac{1}{n}\sum_{i=1}^n(I_{\{\bX_i\in M_{\bX}\}}-P(\bX\in M_{\bX}))\bigg)
		\\[0,2cm]&\quad+\frac{1}{P(\bX\in M_{\bX})}\frac{1}{n}\sum_{i=1}^n\Big(I_{\{\bX_i\leq \bx\}}I_{\{\bX_i\in M_{\bX}\}}I_{\{\varepsilon_i\leq k_c(\mathfrak{h},f_{m_{\tau}},f_{m_{\beta}},\bX_i,e)\}}
		\\[0,2cm]&\quad\quad-P(\bX\in M_{\bX},\bX\leq \bx,\varepsilon\leq k_c(\mathfrak{h},f_{m_{\tau}},f_{m_{\beta}},\bX,e))
		\\[0,2cm]&\quad\quad-I_{\{\bX_i\leq \bx\}}I_{\{\bX_i\in M_{\bX}\}}+P(\bX\in M_{\bX},\bX\leq \bx)-I_{\{\bX_i\in M_{\bX}\}}I_{\{\varepsilon_i\leq k_c(\mathfrak{h},f_{m_{\tau}},f_{m_{\beta}},\bX_i,e)\}}
		\\[0,2cm]&\quad\quad+P(\bX\in M_{\bX},\varepsilon\leq k_c(\mathfrak{h},f_{m_{\tau}},f_{m_{\beta}},\bX,e))\Big)+o_p\bigg(\frac{1}{\sqrt{n}}\bigg),
		\end{align*}
		so that because of Corollary 2.3.1 of \citet{vdVW1996},
		$$\underset{\overset{\scriptstyle c,\tilde{c}\in[B_1,B_2],s,\tilde{s}\in\tilde{\mathcal{H}},\bx\in M_{\bX},e\in[e_a,e_b]}{||s-\tilde{s}||_{\mathcal{H}}<\delta_n,|c-\tilde{c}|<\delta_n}}{\sup}\,\sqrt{n}|G_n(c,s,\bx,e)-G_n(\tilde{c},\tilde{s},\bx,e)|=o_p(1).$$
		\hfill$\square$

		\begin{lemma}\label{lemmac6}
			There exists a real valued function $\psi_{\Gamma_2}$ with $E[\psi_{\Gamma_2}(Y,\bX)]=o\big(n^{-\frac{1}{2}}\big)$ and $\sigma_A^2:=E[\psi_{\Gamma_2}(Y,\bX)^2]\in(0,\infty)$ such that
			\begin{align*}
			&\sqrt{n}\int_{M_{\bX}}\int_{[e_a,e_b]}\Gamma_1(B,s_0)(\bx,e)\Big(G_{nMD}(B,s_0)(\bx,e)+\Gamma_2(B,s_0)(\bx,e)[\hat{s}-s_0]\Big)\,de\,d\bx
			\\[0,2cm]&=\quad\frac{1}{\sqrt{n}}\sum_{i=1}^n\psi_{\Gamma_2}(Y_i,\bX_i)+o_p(1)
			\\[0,2cm]&\overset{\mathcal{D}}{\rightarrow}\mathcal{N}(0,\sigma_A^2).
			\end{align*}
			Moreover, it holds that
			\begin{equation}\label{Gamma2convergence}
			||\Gamma_2(B,s_0)[\hat{s}-s_0]||_2=\mathcal{O}_p\bigg(\frac{1}{\sqrt{n}}\bigg).
			\end{equation}
		\end{lemma}\noindent
		\textbf{Proof:} First, the left hand side will be rewritten such that
		\begin{align}\label{proofasympBlemmac6expint}
		&\sqrt{n}\int_{M_{\bX}}\int_{[e_a,e_b]}\Gamma_1(B,s_0)(\bx,e)\big(G_{nMD}(B,s_0)(\bx,e)+\Gamma_2(B,s_0)(\bx,e)[\hat{s}-s_0]\big)\,de\,d\bx\nonumber
		\\[0,2cm]&=\frac{1}{\sqrt{n}}\sum_{i=1}^n(\psi_{\Gamma_2}(Y_i,\bX_i)-E[\psi_{\Gamma_2}(Y,\bX)])+o_p(1)
		\end{align}
		for some appropriate function $\psi_{\Gamma_2}:\mathbb{R}^{\dX+1}\rightarrow\mathbb{R}$. Afterwards, the usual Central Limit Theorem can be applied to obtain the desired convergence. Third, it will be shown that $||\Gamma_2(B,s_0)[\hat{s}-s_0]||_2=\mathcal{O}_p\big(n^{-\frac{1}{2}}\big)$.
		
		For this purpose, note that it was shown in equation (\ref{sepGamma2}) that
		\begin{align*}
		\Gamma_2(B,s_0)(\bx,e)[\hat{s}-s_0]&=D_sG_{MD}(B,s_0)(\bx,e)[\hat{s}-s_0]
		\\[0,2cm]&=D_hG_{MD}(B,s_0)(\bx,e)[\hat{h}-h_1]
		\\[0,2cm]&\quad+D_{f_{m_{\tau}}}G_{MD}(B,s_0)(\bx,e)[\hat{F}_{Y|\bX}^{-1}(\tau|\cdot)-F_{Y|\bX}^{-1}(\tau|\cdot)\big]
		\\[0,2cm]&\quad+D_{f_{m_{\beta}}}G_{MD}(B,s_0)(\bx,e)[F_{Y|\bX}^{-1}(\beta|\cdot)-F_{Y|\bX}^{-1}(\beta|\cdot)\big].
		\end{align*}
		Hence, there are actually four terms that have to be fitted to Expression (\ref{proofasympBlemmac6expint}). Straightforward calculations lead to
		\begin{align*}
		&\sqrt{n}\int_{M_{\bX}}\int_{[e_a,e_b]}\Gamma_1(B,s_0)(\bx,e)G_{nMD}(B,s_0)(\bx,e)\,de\,d\bx
		\\[0,2cm]&=\frac{1}{\sqrt{n}}\sum_{i=1}^n(\psi_1(Y_i,\bX_i)-E[\psi_1(Y_i,\bX_i)])+o_p(1),
		\end{align*}
		where
		\begin{align*}
		&\psi_1(Y,\bX)
		\\[0,2cm]&=I_{\{\bX\in M_{\bX}\}}\int_{M_{\bX}}\int_{[e_a,e_b]}\Bigg(\frac{P(\bX\leq \bx|\bX\in M_{\bX})P(\varepsilon\leq F_{\varepsilon}^{-1}(\tau)+e(F_{\varepsilon}^{-1}(\beta)-F_{\varepsilon}^{-1}(\tau)))}{P(\bX\in M_{\bX})}
		\\[0,2cm]&\quad\quad+\frac{1}{P(\bX\in M_{\bX})}\bigg(I_{\{\bX\leq \bx\}}I_{\big\{\frac{h_B(Y)-g(\bX)}{\sigma(\bX)}\leq F_{\varepsilon}^{-1}(\tau)+e(F_{\varepsilon}^{-1}(\beta)-F_{\varepsilon}^{-1}(\tau))\big\}}-I_{\{\bX\leq \bx\}}
		\\[0,2cm]&\quad\quad-I_{\big\{\frac{h_B(Y)-g(\bX)}{\sigma(\bX)}\leq F_{\varepsilon}^{-1}(\tau)+e(F_{\varepsilon}^{-1}(\beta)-F_{\varepsilon}^{-1}(\tau))\big\}}\bigg)\Bigg)\Gamma_1(B,s_0)(\bx,e)\,de\,d\bx.
		\end{align*}
		
		Recall from the proof of Lemma \ref{lemmac3} that
		\begin{align}
		&D_hG_{MD}(B,s_0)(\bx,e)[\hat{h}-h_1]\nonumber
		\\[0,2cm]&=\frac{1}{P(\bX\in M_{\bX})}\bigg(\int_{M_{\bX}}f_{\varepsilon}(k_B(s_0,\mathbf{w},e))\big(I_{\{\mathbf{w}\leq \bx\}}-P(\bX\leq \bx|\bX\in M_{\bX})\big)\nonumber
		\\[0,2cm]&\quad D_hk_B(s_0,\mathbf{w},e)[\hat{h}-h_1]f_{\bX}(\mathbf{w})\,d\mathbf{w}\bigg),\label{DhGMD}
		\end{align}
		where
		$$D_hk_B(s_0,\mathbf{w},e)[\hat{h}-h_1]=\frac{Bz_B(0)^{B-1}\big(z_B(0)-\hat{h}(h_1^{-1}(z_B(0)))+\frac{\partial}{\partial t}z_B(t)\big|_{t=0}\big)}{\sigma(\mathbf{w})},$$
		with $z_B$ as in (\ref{defzc}) and $f_{m_{\tau}}=F_{Y|\bX}^{-1}(\tau|\cdot),f_{m_{\beta}}=F_{Y|\bX}^{-1}(\beta|\cdot)$. Furthermore, it was shown in the proof of Lemma \ref{asympintegral} that
		\begin{align*}
		\hat{h}_1(y)-h_1(y)&=\exp\bigg(-\int_{y_1}^y\frac{1}{\hat{\lambda}(u)}\,du\bigg)-\exp\bigg(-\int_{y_1}^y\frac{1}{\lambda(u)}\,du\bigg)
		\\[0,2cm]&=-h_1(y)\int_{y_1}^y\bigg(\frac{1}{\hat{\lambda}(u)}-\frac{1}{\lambda(u)}\bigg)\,du+o_p\bigg(\frac{1}{\sqrt{n}}\bigg)
		\\[0,2cm]&=\frac{1}{n}\sum_{i=1}^n-h_1(y)\eta_i(y)+o_p\bigg(\frac{1}{\sqrt{n}}\bigg),
		\end{align*}
		with $\eta_i$ as in (\ref{defeta}). Remember $E[\eta_1(y)]=0$ uniformly in $[z_a,z_b]$, so that
		\begin{equation}\label{kBhath-h}
		D_hk_B(s_0,\mathbf{w},e)[\hat{h}-h_1]=\frac{1}{n}\sum_{i=1}^n\tilde{\psi}_2(Y_i,\bX_i,\mathbf{w},e)+o_p\bigg(\frac{1}{\sqrt{n}}\bigg)
		\end{equation}
		holds for an appropriate function $\tilde{\psi}_2$, which is centred and uniformly bounded in $(Y,\bX,\mathbf{w},e)\in\mathbb{R}^{\dX+1}\times M_{\bX}\times[e_a,e_b]$. Hence,
		\begin{align*}
		&\sqrt{n}\int_{M_{\bX}}\int_{[e_a,e_b]}\Gamma_1(B,s_0)(\bx,e)D_hG_{MD}(B,s_0)(\bx,e)[\hat{h}-h_1]\,de\,d\bx
		\\[0,2cm]&=\frac{1}{nP(\bX\in M_{\bX})}\sum_{i=1}^n\int_{M_{\bX}}\int_{[e_a,e_b]}\int_{M_{\bX}}\Gamma_1(B,s_0)(\bx,e)f_{\varepsilon}(k_B(s_0,\mathbf{w},e))
		\\[0,2cm]&\quad\quad\big(I_{\{\mathbf{w}\leq \bx\}}-P(\bX\leq \bx|\bX\in M_{\bX})\big)\tilde{\psi}_2(Y_i,\bX_i,\mathbf{w},e)f_{\bX}(\mathbf{w})\,d\mathbf{w}\,de\,d\bx+o_p(1)
		\\[0,2cm]&=\frac{1}{\sqrt{n}}\sum_{i=1}^n\psi_2(Y_i,\bX_i)+o_p(1)
		\end{align*}
		for
		\begin{align*}
		\psi_2(Y_i,\bX_i)&=\frac{1}{P(\bX\in M_{\bX})}\int_{M_{\bX}}\int_{[e_a,e_b]}\int_{M_{\bX}}\Gamma_1(B,s_0)(\bx,e)f_{\varepsilon}(k_B(s_0,\mathbf{w},e))
		\\[0,2cm]&\quad\quad\big(I_{\{\mathbf{w}\leq \bx\}}-P(\bX\leq \bx|\bX\in M_{\bX})\big)\tilde{\psi}_2(Y_i,\bX_i,\mathbf{w},e)f_{\bX}(\mathbf{w})\,d\mathbf{w}\,de\,d\bx.
		\end{align*}
		
		It remains to rewrite $D_{f_{m_{\iota}}}G_{MD}(B,s_0)[\hat{F}_{Y|\bX}^{-1}(\iota|\cdot)-F_{Y|\bX}^{-1}(\iota|\cdot)\big]$ for $\iota\in\{\tau,\beta\}$. In the proof of Lemma \ref{lemmac5}, it was shown that
		\begin{align*}
		&D_{f_{m_{\tau}}}k_B(s_0,\bx,e)[\hat{F}_{Y|\bX}^{-1}(\tau|\cdot)-F_{Y|\bX}^{-1}(\tau|\cdot)\big]
		\\[0,2cm]&=-\frac{B(1-e)h_{B}(F_{Y|\bX}^{-1}(\tau|\bx))(\hat{F}_{Y|\bX}^{-1}(\tau|\bx)-F_{Y|\bX}^{-1}(\tau|\bx))}{\sigma(\bx)\lambda(F_{Y|\bX}^{-1}(\tau|\bx))}
		\end{align*}
		and
		$$D_{f_{\beta}}k_B(s_0,\bx,e)[\hat{F}_{Y|\bX}^{-1}(\beta|\cdot)-F_{Y|\bX}^{-1}(\beta|\cdot)\big]=-\frac{Beh_{B}(F_{Y|\bX}^{-1}(\beta|\bx))(\hat{F}_{Y|\bX}^{-1}(\beta|\bx)-F_{Y|\bX}^{-1}(\beta|\bx))}{\sigma(\bx)\lambda(F_{Y|\bX}^{-1}(\beta|\bx))}.$$
		Equation (\ref{FYXhatinverse}) yields
		\begin{align*}
		\hat{F}_{Y|\bX}^{-1}(\iota|\bx)-F_{Y|\bX}^{-1}(\iota|\bx)&=\frac{1}{f_{Y|\bX}(F_{Y|\bX}^{-1}(\iota|\bx)|\bx)f_{\bX}(\bx)}\frac{1}{n}\sum_{i=1}^n\bK_{h_x}(\bx-\bX_i)
		\\[0,2cm]&\quad\quad\Bigg(\mathcal{K}_{h_y}(F_{Y|\bX}^{-1}(\iota|\bx)-Y_i)-\frac{p(F_{Y|\bX}^{-1}(\iota|\bx),\bx)}{f_{\bX}(\bx)}\Bigg)+o_p\bigg(\frac{1}{\sqrt{n}}\bigg)
		\end{align*}
		for $\iota\in\{\tau,\beta\}$. Note that the order of the remaining term, that is obtained in the proof there, is actually $o_p(\hat{F}_{Y|\bX}^{-1}(\iota|\bx)-F_{Y|\bX}^{-1}(\iota|\bx))$, but this order can be extended to
		$$\mathcal{O}_p((\hat{F}_{Y|\bX}^{-1}(\iota|\bx)-F_{Y|\bX}^{-1}(\iota|\bx))^2)=o_p(n^{-\frac{1}{2}})$$
		by using the Lagrange form of the remainder of the Taylor expansion. Due to \ref{B3}, one has
		\begin{equation}\label{expectmathcalKp}
		E\bigg[\bK_{h_x}(\bx-\bX_1)\bigg(\mathcal{K}_{h_y}(F_{Y|\bX}^{-1}(\iota|\bx)-Y_1)-\frac{p(F_{Y|\bX}^{-1}(\iota|\bx),\bx)}{f_{\bX}(\bx)}\bigg)\bigg]=o_p\bigg(\frac{1}{\sqrt{n}}\bigg)
		\end{equation}
		uniformly in $\bx\in M_{\bX}$ and $\iota\in\{\tau,\beta\}$. So far, a representation
		\begin{align*}
		&\sqrt{n}\int_{M_{\bX}}\int_{[e_a,e_b]}\Gamma_1(B,s_0)(\bx,e)\big(G_{nMD}(B,s_0)(\bx,e)+\Gamma_2(B,s_0)(\bx,e)[\hat{s}-s_0]\big)\,de\,d\bx
		\\[0,2cm]&=\frac{1}{\sqrt{n}}\sum_{i=1}^n\big(\psi_1(Y_i,\bX_i)+\psi_2(Y_i,\bX_i)+\psi_{3,n}(Y_i,\bX_i)+\psi_{4,n}(Y_i,\bX_i)
		\\[0,2cm]&\quad-E[\psi_1(Y,\bX)-\psi_2(Y,\bX)-\psi_{3,n}(Y,\bX)-\psi_{4,n}(Y,\bX)]\big)+o_p(1),
		\end{align*}
		was found, where
		\begin{align*}
		\psi_{3,n}(Y_i,\bX_i)&=-\int_{M_{\bX}}\int_{[e_a,e_b]}\int_{M_{\bX}}\Gamma_1(B,s_0)(\bx,e)
		\big(I_{\{\mathbf{w}\leq \bx\}}-P(\bX\leq \bx|\bX\in M_{\bX})\big)
		\\[0,2cm]&\quad\quad\frac{B(1-e)h_{B}(F_{Y|\bX}^{-1}(\tau|\mathbf{w}))f_{\varepsilon}(k_B(s_0,\mathbf{w},e))}{\sigma(\mathbf{w})\lambda(F_{Y|\bX}^{-1}(\tau|\mathbf{w}))f_{Y|\bX}(F_{Y|\bX}^{-1}(\tau|\mathbf{w})|\mathbf{w})}
		\\[0,2cm]&\quad\quad\bK_{h_x}(\mathbf{w}-\bX_i)\bigg(\mathcal{K}_{h_y}(F_{Y|\bX}^{-1}(\tau|\mathbf{w})-Y_i)-\frac{p(F_{Y|\bX}^{-1}(\tau|\mathbf{w}),\mathbf{w})}{f_{\bX}(\mathbf{w})}\bigg)\,d\mathbf{w}\,de\,d\bx
		\end{align*}
		and
		\begin{align*}
		\psi_{4,n}(Y_i,\bX_i)&=-\int_{M_{\bX}}\int_{[e_a,e_b]}\int_{M_{\bX}}\Gamma_1(B,s_0)(\bx,e)
		\big(I_{\{\mathbf{w}\leq \bx\}}-P(\bX\leq \bx|\bX\in M_{\bX})\big)
		\\[0,2cm]&\quad\quad\frac{Beh_{B}(F_{Y|\bX}^{-1}(\beta|\mathbf{w}))f_{\varepsilon}(k_B(s_0,\mathbf{w},e))}{\sigma(\mathbf{w})\lambda(F_{Y|\bX}^{-1}(\beta|\mathbf{w}))f_{Y|\bX}(F_{Y|\bX}^{-1}(\beta|\mathbf{w})|\mathbf{w})}
		\\[0,2cm]&\quad\quad\bK_{h_x}(\mathbf{w}-\bX_i)\bigg(\mathcal{K}_{h_y}(F_{Y|\bX}^{-1}(\beta|\mathbf{w})-Y_i)-\frac{p(F_{Y|\bX}^{-1}(\beta|\mathbf{w}),\mathbf{w})}{f_{\bX}(\mathbf{w})}\bigg)\,d\mathbf{w}\,de\,d\bx
		\end{align*}
		depend on $n$. To fit this expression to equation (\ref{proofasympBlemmac6expint}) it suffices to replace $\psi_{3,n}$ and $\psi_{4,n}$ with some functions $\psi_3$ and $\psi_4$ (independent of $n$ and with finite second moments), respectively, such that
		\begin{equation}\label{psi3psi4}
		E[(\psi_{3,n}(Y,\bX)-\psi_3(Y,\bX))^2]=o(1)\quad\textup{and}\quad E[(\psi_{4,n}(Y,\bX)-\psi_4(Y,\bX))^2]=o(1),
		\end{equation}
		since it was already shown that $E[\psi_1(Y,\bX)^2],E[\psi_2(Y,\bX)^2]<\infty$. It can be shown that equation (\ref{psi3psi4}) is fulfilled for $\psi_3$ and $\psi_4$, where
		\begin{align}
		\psi_{\tau}(\bx,e)&:=-\frac{B(1-e)h_{B}(F_{Y|\bX}^{-1}(\tau|\bx))f_{\varepsilon}(k_B(s_0,\bx,e))}{\sigma(\bx)\lambda(F_{Y|\bX}^{-1}(\tau|\bx))f_{Y|\bX}(F_{Y|\bX}^{-1}(\tau|\bx)|\bx)},\label{psitau}
		\\[0,2cm]\psi_{\beta}(\bx,e)&:=-\frac{Beh_{B}(F_{Y|\bX}^{-1}(\beta|\bx))f_{\varepsilon}(k_B(s_0,\bx,e))}{\sigma(\bx)\lambda(F_{Y|\bX}^{-1}(\beta|\bx))f_{Y|\bX}(F_{Y|\bX}^{-1}(\beta|\bx)|\bx)},\nonumber
		\\[0,2cm]\psi_3(Y,\bX)&:=\bigg(I_{\{Y\leq F_{Y|\bX}^{-1}(\tau|\bX)\}}-\frac{p(F_{Y|\bX}^{-1}(\tau|\bX),\bX)}{f_{\bX}(\bX)}\bigg)\nonumber
		\\[0,2cm]&\quad\quad\int_{M_{\bX}}\int_{[e_a,e_b]}\Gamma_1(B,s_0)(\bx,e)\big(I_{\{\bX\leq \bx\}}-P(\bX\leq \bx|\bX\in M_{\bX})\big)\psi_{\tau}(\bX,e)\,d\bx\,de,\nonumber
		\\[0,2cm]\psi_4(Y,\bX)&:=\bigg(I_{\{Y\leq F_{Y|\bX}^{-1}(\beta|\bX)\}}-\frac{p(F_{Y|\bX}^{-1}(\beta|\bX),\bX)}{f_{\bX}(\bX)}\bigg)\nonumber
		\\[0,2cm]&\quad\quad\int_{M_{\bX}}\int_{[e_a,e_b]}\Gamma_1(B,s_0)(\bx,e)\big(I_{\{\bX\leq \bx\}}-P(\bX\leq \bx|\bX\in M_{\bX})\big)\psi_{\beta}(\bX,e)\,d\bx\,de.\nonumber
		\end{align}
		In total, this leads to
		\begin{align*}
		&\sqrt{n}\int_{M_{\bX}}\int_{[e_a,e_b]}\Gamma_1(B,s_0)(\bx,e)\big(G_{nMD}(B,s_0)(\bx,e)+\Gamma_2(B,s_0)(\bx,e)[\hat{s}-s_0]\big)\,de\,d\bx
		\\[0,2cm]&=\frac{1}{\sqrt{n}}\sum_{i=1}^n(\psi_{\Gamma_2}(Y_i,\bX_i)-E[\psi(Y,\bX)])+o_p(1)
		\end{align*}
		for
		$$\psi_{\Gamma_2}(Y,\bX)=\psi_1(Y,\bX)+\psi_2(Y,\bX)+\psi_3(Y,\bX)+\psi_4(Y,\bX).$$
		One has $E[\psi_{\Gamma_2}(Y,\bX)]=0$, so that the Central Limit Theorem implies
		\begin{align*}
		&\sqrt{n}\int_{M_{\bX}}\int_{[e_a,e_b]}\Gamma_1(B,s_0)(\bx,e)\big(G_{nMD}(B,s_0)(\bx,e)+\Gamma_2(B,s_0)(\bx,e)[\hat{s}-s_0]\big)\,de\,d\bx
		\\[0,2cm]&\overset{\mathcal{D}}{\rightarrow}\mathcal{N}(0,\sigma_A^2)
		\end{align*}
		with $\sigma_A^2=\operatorname{Var}(\psi(Y,\bX))$.
		
		To prove Lemma \ref{lemmac6}, it remains to prove (\ref{Gamma2convergence}). In the following, the complexity of the dominating term in $||\Gamma_2(B,s_0)[\hat{s}-s_0]||_2$ will be reduced stepwise. First, apply (\ref{sepGamma2}) to obtain
		\begin{align*}
		&||\Gamma_2(B,s_0)[\hat{s}-s_0]||_2
		\\[0,2cm]&\leq\big|\big|D_hG_{MD}(c,s_0)[\hat{h}_1-h_1]\big|\big|_2+\big|\big|D_{f_{m_{\tau}}}G_{MD}(B,s_0)\big[\hat{F}_{Y|\bX}^{-1}(\tau|\cdot)-F_{Y|\bX}^{-1}(\tau|\cdot)\big]\big|\big|_2
		\\[0,2cm]&\quad+\big|\big|D_{f_{m_{\beta}}}G_{MD}(B,s_0)\big[\hat{F}_{Y|\bX}^{-1}(\beta|\cdot)-F_{Y|\bX}^{-1}(\beta|\cdot)\big]\big|\big|_2
		\\[0,2cm]&=\big|\big|D_{f_{m_{\tau}}}G_{MD}(B,s_0)\big[\hat{F}_{Y|\bX}^{-1}(\tau|\cdot)-F_{Y|\bX}^{-1}(\tau|\cdot)\big]\big|\big|_2
		\\[0,2cm]&\quad+\big|\big|D_{f_{m_{\beta}}}G_{MD}(B,s_0)\big[\hat{F}_{Y|\bX}^{-1}(\beta|\cdot)-F_{Y|\bX}^{-1}(\beta|\cdot)\big]\big|\big|_2+\mathcal{O}_p\bigg(\frac{1}{\sqrt{n}}\bigg),
		\end{align*}
		where the last equation follows from the equations (\ref{DhGMD}) and (\ref{kBhath-h}). Both of the terms
		$$\big|\big|D_{f_{m_{\tau}}}G_{MD}(B,s_0)\big[\hat{F}_{Y|\bX}^{-1}(\tau|\cdot)-F_{Y|\bX}^{-1}(\tau|\cdot)\big]\big|\big|_2$$
		and
		$$\big|\big|D_{f_{m_{\beta}}}G_{MD}(B,s_0)\big[\hat{F}_{Y|\bX}^{-1}(\beta|\cdot)-F_{Y|\bX}^{-1}(\beta|\cdot)\big]\big|\big|_2$$
		can be treated similarly to each other, so that only the first term is considered in the following. Recall
		\begin{align*}
		&D_{f_{m_{\tau}}}G_{MD}(B,s_0)(\bx,e)\big[\hat{F}_{Y|\bX}^{-1}(\tau|\cdot)-F_{Y|\bX}^{-1}(\tau|\cdot)\big]
		\\[0,2cm]&=\frac{1}{P(\bX\in M_{\bX})}\Big(\int_{M_{\bX}}f_{\varepsilon}(k_B(s_0,\mathbf{w},e))\big(I_{\{\mathbf{w}\leq \bx\}}-P(\bX\leq \bx|\bX\in M_{\bX})\big)
		\\[0,2cm]&\quad\quad D_{f_{m_{\tau}}}k_B(s_0,\mathbf{w},e)\big[\hat{F}_{Y|\bX}^{-1}(\tau|\cdot)-F_{Y|\bX}^{-1}(\tau|\cdot)\big]f_{\bX}(\mathbf{w})\,d\mathbf{w}\Big)
		\\[0,2cm]&=-\frac{1}{P(\bX\in M_{\bX})}\int_{M_{\bX}}f_{\varepsilon}(k_B(s_0,\mathbf{w},e))\big(I_{\{\mathbf{w}\leq \bx\}}-P(\bX\leq \bx|\bX\in M_{\bX})\big)
		\\[0,2cm]&\quad\quad \frac{B(1-e)h(F_{Y|\bX}^{-1}(\tau|\mathbf{w}))(\hat{F}_{Y|\bX}^{-1}(\tau|\mathbf{w})-F_{Y|\bX}^{-1}(\tau|\mathbf{w}))}{\sigma(\mathbf{w})\lambda(F_{Y|\bX}^{-1}(\tau|\mathbf{w}))}f_{\bX}(\mathbf{w})\,d\mathbf{w}.
		\end{align*}
		Inserting equation (\ref{FYXhatinverse}) leads to
		\begin{align*}
		&D_{f_{m_{\tau}}}G_{MD}(B,s_0)(\bx,e)\big[\hat{F}_{Y|\bX}^{-1}(\tau|\cdot)-F_{Y|\bX}^{-1}(\tau|\cdot)\big]
		\\[0,2cm]&=-\frac{1}{nP(\bX\in M_{\bX})}\sum_{i=1}^n\int_{M_{\bX}}\big(I_{\{\mathbf{w}\leq \bx\}}-P(\bX\leq \bx|\bX\in M_{\bX})\big)\psi_{\tau}(\mathbf{w},e)
		\\[0,2cm]&\quad\quad\Bigg(\mathcal{K}_{h_y}(F_{Y|\bX}^{-1}(\tau|\mathbf{w})-Y_i)-\frac{p(F_{Y|\bX}^{-1}(\tau|\mathbf{w}),\mathbf{w})}{f_{\bX}(\mathbf{w})}\Bigg)\bK_{h_x}(\mathbf{w}-\bX_i)f_{\bX}(\mathbf{w})\,d\mathbf{w}+o_p\bigg(\frac{1}{\sqrt{n}}\bigg)
		\\[0,2cm]&=-\frac{1}{nP(\bX\in M_{\bX})}\sum_{i=1}^n\int_{M_{\bX}}\big(I_{\{\bX_i+h_x\mathbf{w}\leq \bx\}}-P(\bX\leq \bx|\bX\in M_{\bX})\big)\psi_{\tau}(\bX_i+h_x\mathbf{w},e)
		\\[0,2cm]&\quad\quad\Bigg(\mathcal{K}_{h_y}(F_{Y|\bX}^{-1}(\tau|\bX_i+h_x\mathbf{w})-Y_i)-\frac{p(F_{Y|\bX}^{-1}(\tau|\bX_i+h_x\mathbf{w}),\bX_i+h_x\mathbf{w})}{f_{\bX}(\bX_i+h_x\mathbf{w})}\Bigg)
		\\[0,2cm]&\quad\quad\bK(\mathbf{w})f_{\bX}(\bX_i+h_x\mathbf{w})\,d\mathbf{w}+o_p\bigg(\frac{1}{\sqrt{n}}\bigg)
		\\[0,2cm]&=:\frac{1}{n}\sum_{i=1}^nZ_{i,h_x,\tau}(\bx,e)+o_p\bigg(\frac{1}{\sqrt{n}}\bigg)
		\end{align*}
		uniformly in $\bx\in M_{\bX},e\in[e_a,e_b]$, where $\psi_{\tau}$ was defined in (\ref{psitau}). Equation (\ref{expectmathcalKp}) yields $E[Z_{1,h_x,\tau}(\bx,e)]=o\big(n^{-\frac{1}{2}}\big)$ uniformly in $\bx\in M_{\bX},e\in[e_a,e_b]$, so that
		$$E\bigg[\bigg|\bigg|\frac{1}{\sqrt{n}}\sum_{i=1}^nZ_{i,h_x,\tau}\bigg|\bigg|_2^2\bigg]=E[Z_{1,h_x,\tau}^2]+o(1)\leq C$$
		for some sufficiently large $C>0$. Therefore, it holds that
		$$||\Gamma_2(B,s_0)[\hat{s}-s_0]||_2=\mathcal{O}_p\bigg(\frac{1}{\sqrt{n}}\bigg).$$
		\hfill$\square$

		\subsubsection{Proof of Lemma \ref{asympB}}
		Similar to \citet{LSvK2008}, define
		$$L_n(\bx,e)=G_{nMD}(B,s_0)(\bx,e)-G_{MD}(B,s_0)(\bx,e)$$
		as well as
		$$\mathcal{L}_n(c)(\bx,e)=L_n(\bx,e)+\Gamma_{1}(B,s_0)(\bx,e)(c-B)+\Gamma_2(B,s_0)(\bx,e)[\hat{s}-s_0].$$
		In the proof of Lemma \ref{lemmac5}, it was shown that
		\begin{equation}\label{convLn}
		||L_n||_2=||G_{nMD}(B,s_0)||_2=\mathcal{O}_p\bigg(\frac{1}{\sqrt{n}}\bigg).
		\end{equation}
		Then, one has for all sequences $\delta_n\searrow0$
		\begin{align}
		&||G_{nMD}(c,\hat{s})-\mathcal{L}_n(c)||_2\nonumber
		\\[0,2cm]&\quad\ =||G_{nMD}(c,\hat{s})-G_{nMD}(B,s_0)+G_{MD}(B,s_0)\nonumber
		\\[0,2cm]&\quad\quad\ -\Gamma_1(B,s_0)(c-B)-\Gamma_2(B,s_0)[\hat{s}-s_0]||_2\nonumber
		\\[0,2cm]&\quad\overset{\ref{C1}}{\leq}||G_{nMD}(c,\hat{s})-G_{MD}(c,\hat{s})-G_{nMD}(B,s_0)||_2\nonumber
		\\[0,2cm]&\quad\quad\ +||G_{MD}(c,\hat{s})-\Gamma_1(B,s_0)(c-B)-\Gamma_2(B,s_0)[\hat{s}-s_0]||_2\nonumber
		\\[0,2cm]&\quad\overset{\ref{C5}}{=}||G_{MD}(c,\hat{s})-\Gamma_1(B,s_0)(c-B)-\Gamma_2(B,s_0)[\hat{s}-s_0]||_2+o_p\bigg(\frac{1}{\sqrt{n}}\bigg)\nonumber
		\\[0,2cm]&\quad\ \leq||G_{MD}(c,\hat{s})-G_{MD}(c,s_0)-\Gamma_2(B,s_0)[\hat{s}-s_0]||_2\nonumber
		\\[0,2cm]&\quad\quad\ +||G_{MD}(c,s_0)-\Gamma_1(B,s_0)(c-B)||_2+o_p\bigg(\frac{1}{\sqrt{n}}\bigg)\nonumber
		\\[0,2cm]&\quad\overset{\ref{C3}}{\leq}||G_{MD}(c,\hat{s})-G_{MD}(c,s_0)-\Gamma_2(c,s_0)[\hat{s}-s_0]||_2\nonumber
		\\[0,2cm]&\quad\quad\ +||G_{MD}(c,s_0)-\Gamma_1(B,s_0)(c-B)||_2+o_p(|c-B|)+o_p\bigg(\frac{1}{\sqrt{n}}\bigg)\nonumber
		\\[0,2cm]&\overset{\ref{C3}+\ref{C4}}{\leq}||G_{MD}(c,s_0)-G_{MD}(B,s_0)-\Gamma_1(B,s_0)(c-B)||_2+o_p(|c-B|)+o_p\bigg(\frac{1}{\sqrt{n}}\bigg)\nonumber
		\\[0,2cm]&\quad\overset{\ref{C2}}{=}o_p(|c-B|)+o_p\bigg(\frac{1}{\sqrt{n}}\bigg)\label{GnMDL}
		\end{align}
		uniformly in $c\in B_{\delta_n}$. Denote the minimizer of $c\mapsto||\mathcal{L}_n(c)||_2$ by $\bar{B}$. Then, $\bar{B}$ can be calculated explicitly by solving
		\begin{align*}
		&\frac{\partial}{\partial c}||\mathcal{L}_n(c)||_2^2
		\\[0,2cm]&=\frac{\partial}{\partial c}\bigg(||L_n||_2^2+||\Gamma_{1}(B,s_0)||_2^2(c-B)^2+||\Gamma_2(B,s_0)[\hat{s}-s_0]||_2^2
		\\[0,2cm]&\quad+2\int_{M_{\bX}}\int_{[e_a,e_b]}\Gamma_{1}(B,s_0)(\bx,e)(L_n(\bx,e)+\Gamma_2(B,s_0)(\bx,e)[\hat{s}-s_0])\,de\,d\bx(c-B)
		\\[0,2cm]&\quad+2\int_{M_{\bX}}\int_{[e_a,e_b]}L_n(\bx,e)\Gamma_2(B,s_0)(\bx,e)[\hat{s}-s_0]\,de\,d\bx\bigg)
		\\[0,2cm]&=2||\Gamma_{1}(B,s_0)||_2^2(c-B)
		\\[0,2cm]&\quad+2\int_{M_{\bX}}\int_{[e_a,e_b]}\Gamma_{1}(B,s_0)(\bx,e)(L_n(\bx,e)+\Gamma_2(B,s_0)(\bx,e)[\hat{s}-s_0])\,de\,d\bx
		\\[0,2cm]&=0.
		\end{align*}
		Therefore, one has
		\begin{align}
		\bar{B}&\ \,=B-\frac{\int_{M_{\bX}}\int_{[e_a,e_b]}\Gamma_{1}(B,s_0)(\bx,e)(L_n(\bx,e)+\Gamma_2(B,s_0)(\bx,e)[\hat{s}-s_0])\,de\,d\bx}{||\Gamma_{1}(B,s_0)||_2^2}\nonumber
		\\[0,2cm]&\overset{\ref{C6}}{=}B+\mathcal{O}_p\bigg(\frac{1}{\sqrt{n}}\bigg).\label{BBbar}
		\end{align}
		One has
		$$||\mathcal{L}_n(\bar{B})||_2\leq||\mathcal{L}_n(B)||_2\leq||L_n(\bx,e)||_2+||\Gamma_2(B,s_0)[\hat{s}-s_0]||_2\overset{(\ref{Gamma2convergence})+(\ref{convLn})}{=}\mathcal{O}_p\bigg(\frac{1}{\sqrt{n}}\bigg),$$
		$$||G_{nMD}(\hat{B},\hat{s})||_2\leq||G_{nMD}(B,\hat{s})||_2\overset{(\ref{GnMDL})}{\leq}||\mathcal{L}_n(B)||_2+o_p\bigg(\frac{1}{\sqrt{n}}\bigg)=\mathcal{O}_p\bigg(\frac{1}{\sqrt{n}}\bigg).$$
		For all $c\in B_{\delta_n}$ a Taylor expansion yields for some $B^*$ between $c$ and $\bar{B}$
		\begin{align}
		||\mathcal{L}_n(c)||_2^2&=||\mathcal{L}_n(\bar{B})||_2^2+\frac{\partial}{\partial c}||\mathcal{L}_n(c)||_2^2\bigg|_{c=\bar{B}}(c-\bar{B})+\frac{\frac{\partial^2}{\partial c^2}||\mathcal{L}_n(c)||_2^2\big|_{c=B^*}}{2}(c-\bar{B})^2\nonumber
		\\[0,2cm]&=||\mathcal{L}_n(\bar{B})||_2^2+||\Gamma_1(B,s_0)||_2^2(c-\bar{B})^2.\label{LcLBbar}
		\end{align}
		These assertions in turn can be used to obtain
		\begin{align*}
		&||G_{nMD}(\hat{B},\hat{s})||_2^2
		\\[0,2cm]&\ \ \leq||G_{nMD}(\bar{B},\hat{s})||_2^2
		\\[0,2cm]&\overset{(\ref{GnMDL})}{=}\big(||\mathcal{L}_n(\bar{B})||_2+o_p(|\bar{B}-B|)+o_p\big(n^{-\frac{1}{2}}\big)\big)^2
		\\[0,2cm]&\ \ =||\mathcal{L}_n(\bar{B})||_2^2+||\mathcal{L}_n(\bar{B})||_2o_p\big(n^{-\frac{1}{2}}\big)+o_p(n^{-1})
		\\[0,2cm]&\overset{(\ref{LcLBbar})}{=}||\mathcal{L}_n(\hat{B})||_2^2-||\Gamma_1(B,s_0)||_2^2(\hat{B}-\bar{B})^2+o_p(n^{-1})
		\\[0,2cm]&\overset{(\ref{GnMDL})}{=}\big(||G_{nMD}(\hat{B},\hat{s})||_2+o_p(|\hat{B}-B|)+o_p\big(n^{-\frac{1}{2}}\big)\big)^2-||\Gamma_1(B,s_0)||_2^2(\hat{B}-\bar{B})^2+o_p(n^{-1})
		\\[0,2cm]&\overset{(\ref{BBbar})}{=}\big(||G_{nMD}(\hat{B},\hat{s})||_2+o_p(|\hat{B}-\bar{B}|)+o_p\big(n^{-\frac{1}{2}}\big)\big)^2-||\Gamma_1(B,s_0)||_2^2(\hat{B}-\bar{B})^2+o_p(n^{-1})
		\\[0,2cm]&\ \ =||G_{nMD}(\hat{B},\hat{s})||_2^2-||\Gamma_1(B,s_0)||_2^2(\hat{B}-\bar{B})^2+o_p\big(n^{-\frac{1}{2}}|\hat{B}-\bar{B}|\big)+o_p(|\hat{B}-\bar{B}|^2)
		\\[0,2cm]&\ \ \quad+o_p(n^{-1}).
		\end{align*}
		Thus,
		$$||\Gamma_1(B,s_0)||_2^2(\hat{B}-\bar{B})^2=o_p\big(n^{-\frac{1}{2}}|\hat{B}-\bar{B}|\big)+o_p(|\hat{B}-\bar{B}|^2)+o_p(n^{-1})$$
		and consequently $\hat{B}-\bar{B}=o_p(n^{-\frac{1}{2}})$. Finally, \ref{C6} yields
		\begin{align}
		&\sqrt{n}(\hat{B}-B)\nonumber
		\\[0,2cm]&=\sqrt{n}(\bar{B}-B)+o_p(1)\nonumber
		\\[0,2cm]&=-\frac{\sqrt{n}\int_{M_{\bX}}\int_{[e_a,e_b]}\Gamma_{1}(B,s_0)(\bx,e)(G_{nMD}(\bx,e)+\Gamma_2(B,s_0)(\bx,e)[\hat{s}-s_0])\,de\,d\bx}{||\Gamma_{1}(B,s_0)||_2^2}+o_p(1)\nonumber
		\\[0,2cm]&\overset{\mathcal{D}}{\rightarrow}\mathcal{N}\bigg(0,\frac{\sigma_A^2}{||\Gamma_{1}(B,s_0)||_2^4}\bigg).\label{convBhat}
		\end{align}
		\hfill$\square$

		\subsection{Proof of Lemma \ref{asympaltB}}\label{proofasympaltB}
		Recall the definition of $\tilde{B}$ from equation (\ref{defestBtilde})
		$$\tilde{B}=-\frac{\partial}{\partial y}\hat{\lambda}(y)\Big|_{y=\hat{y}_0}.$$
		First, a Taylor expansion leads to
		$$\tilde{B}=-\frac{\partial}{\partial y}\hat{\lambda}(y)\Big|_{y=\hat{y}_0}=-\frac{\partial}{\partial y}\hat{\lambda}(y)\Big|_{y=y_0}-\frac{\partial^2}{\partial y^2}\hat{\lambda}(y)\Big|_{y=y^*}(\hat{y}_0-y_0)$$
		for some $y^*$ between $\hat{y}_0$ and $y_0$. Similarly to the result of \citet{Han2008}, it can be shown that $\frac{\partial^2}{\partial y^2}\hat{\lambda}$ converges uniformly on compact sets to $\frac{\partial^2}{\partial y^2}{\lambda}$, which in turn is continuous under \ref{B1}--\ref{B5}, so that Lemma \ref{asympy0alpha2} implies
		$$\tilde{B}=-\frac{\partial}{\partial y}\hat{\lambda}(y)\Big|_{y=y_0}-\frac{\partial^2}{\partial y^2}\lambda(y)\Big|_{y=y_0}(\hat{y}_0-y_0)=-\frac{\partial}{\partial y}\hat{\lambda}(y)\Big|_{y=y_0}+\mathcal{O}_p\bigg(\frac{1}{\sqrt{nh_y}}\bigg).$$
		Therefore, it suffices to treat the first summand. The assertion follows from some tedious calculations together with a Lindeberg-Feller-Theorem. See \citet{Klo2019} for details.
		\hfill$\square$

	\bibliographystyle{plainnat}
	\bibliography{References}

\end{appendix}

\end{document}